\documentclass{amsart}
\usepackage[active]{srcltx}
\usepackage{amsmath}
\usepackage{amssymb}
\usepackage{epsfig}
\usepackage{multicol}
\usepackage{hyperref}
\newtheorem{Proposition}{Proposition}
  \newtheorem{Remark}[Proposition]{Remark}
  \newtheorem{Corollary}[Proposition]{Corollary}
  \newtheorem{Lemma}[Proposition]{Lemma}
  
  \newtheorem{Theorem}{Theorem}
 
 \newtheorem{Definition}[Proposition]{Definition}
 \newtheorem{Note}[Proposition]{Note}
 
 \newtheorem{Assumption}[Proposition]{Assumption}

\usepackage{multicol}
\makeatletter
  {\end{multicols}\if@restonecol\onecolumn\else\clearpage\fi}
\makeatother

\newcommand {\z}{{\noindent}}
\def\blackslug{\hbox{\hskip 1pt \vrule width 4pt height 8pt depth 1.5pt
\hskip 1pt}}
\def\qed{\quad\blackslug\lower 8.5pt\null\par}
 \def\HH{\mathbb{H}}
\def\CC{\mathbb{C}}
 \def\RR{\mathbb{R}}
 \def\NN{\mathbb{N}}
\def\ZZ{\mathbb{Z}}

\def\Re{\mathrm{Re}}
\def\Im{\mathrm{Im}}
 \def\bchi{\mbox{\raisebox{.4ex}{\begin{large}$\chi$\end{large}}}}
 \renewcommand{\thefootnote}{(\arabic{footnote})}

  \makeindex
\begin{document}

\title[Ionization of Coulomb systems in $\RR^3$]{Ionization 
of Coulomb systems in $\RR^3$ by time periodic forcings of
  arbitrary size} 

 \author[O. Costin, J. L.  Lebowitz  and S.  Tanveer] {O. Costin$^1$, J. L.  Lebowitz$^{2,3}$  and S.  Tanveer$^1$ }
\thanks{$1$.  Department of Mathematics, Ohio State University.}  \thanks{$2$.
  Department of Mathematics, Rutgers University } \thanks{$3$.  Department of
  Physics, Rutgers University.}  
\maketitle

\bigskip
 
\centerline{\today}

\begin{abstract}
  
  We analyze the long time behavior of solutions of the Schr\"odinger
  equation $i\psi_t=(-\Delta-b/r+V(t,x))\psi$, $x\in\RR^3$, $r=|x|$,
  describing a Coulomb system subjected to a spatially compactly
  supported 
  time periodic potential $V(t,x)=V(t+2\pi/\omega,x)$ with
  zero time average.

We  show that, for any $V(t,x)$ of the form $2\Omega(r)\sin (\omega t-\theta)$, with $\Omega(r)$ nonzero on its support, Floquet bound states do not
  exist. This implies that the system ionizes, {\em i.e.}  $P(t,K)=\int_K|\psi(t,x)|^2dx\to 0$ as
  $t\to\infty$ for any compact set $K\subset\RR^3$. 
 Furthermore,
if the initial state is compactly supported and  has only finitely many 
spherical harmonic modes, then $P(t,K)$ decays like $t^{-5/3}$ as 
  $t \rightarrow \infty $. 

To prove these statements, we develop  a rigorous WKB theory 
for  infinite systems of ordinary differential equations.  \bigskip

\z {\bf Keywords:} Ionization, time periodic Schr\"odinger equation, resonances, Floquet theory, Coulomb potential.

\end{abstract}
\tableofcontents

\section{Introduction and overview of results}

The long time behavior of solutions of the Schr\"odinger equation of a system
with both discrete and continuous spectrum subjected to a time periodic
potential is a longstanding problem.  Powerful results  have been obtained under various
assumptions on the potentials, see
\cite{Bourgain4,Bourgain1,Bourgain2,Bourgain3,MG,Moller, RT,Schlag1,Skibsted},
and references therein. In particular, there are conditional  results on the ionization of the
Hydrogen atom, subjected to an external time-harmonic dipole field
$V(t,x)=E\cdot x \cos\omega t$ if $E$ is sufficiently small, see
\cite{YajimaCMP1,YajimaCMP}. In addition,  M\"oller and  Skibsted  proved the equivalence of absence of point spectrum and ionization
for a large class of such  systems subject to periodic fields  \cite{Moller}. There are also detailed results about the
behavior of the wave function for systems subjected to general time
periodic potentials, decaying faster than $r^{-2}$, under the
additional assumption of absence of point spectrum of the Floquet operator, see \cite{Yajima}.

None of these results however prove or disprove ionization of
Coulomb-- bound particles subject to time-periodic forcing of fixed
amplitude and zero average.  In fact, such results have only recently been
obtained even for simple model systems, see
\cite{{JPA},{CMP221},{CRM},{Stucchio},{JSPLieb},[6a]} and references
cited there.  For a periodic dipole field of nonzero average ionization
was proved in \cite{Moller2} (we note that the time averaged Hamiltonian has no bound states in this seetting).

What experiments and simplified models show is that the behavior of
systems with both discrete and continuous spectrum, subject to
time-periodic fields of arbitrary strength, can be very
complicated. For amplitudes where perturbation theory is not
applicable (such fields are becoming of increasing practical importance in technology), qualitative departures from the behavior at small fields
are observed. There are even situations, see {\em e.g.}
\cite{CMP221}, where for small enough fields ionization occurs for all
initial states while for larger fields there exist localized
time--{\em quasiperiodic} solutions of the Schr\"odinger equation,
{\em i.e.}  Floquet bound states.  Though these situations are rather
exceptional, constructive methods of analysis are required to determine
the outcome in specific settings.

In this paper we prove ionization for Coulomb systems with very
special (non-dipole) type of forcings
of arbitrary
magnitude. This is equivalent to establishing the absence of point and 
singular continuous spectrum of the corresponding Floquet operators. 
We also obtain the large time behavior of the wave function. The
time decay  of the wave function,  for compactly supported initial conditions,
is 
of order $t^{-5/6}$. This differs from the $t^{-3/2}$ or,
exceptionally, $t^{-1/2}$ power law found for
shorter range reference potentials, see \cite{Yajima,JSPLieb}.
  The nonperturbative methods include the
development of rigorous WKB techniques for infinite systems of ODEs.  
\subsection{The Coulomb Hamiltonian}\label{sH} In units such that
$\hbar^2/2m=1$, the Coulomb quantum Hamiltonian of a Hydrogen atom (more
generally a Rydberg atom) is \begin{equation}
  \label{eq:eqH}
  H_C=-\Delta-\frac{b}{r}
\end{equation}
where $b>0$, $r=|x|$, $x\in\RR^3$ and $\Delta$ is the Laplacian. It is
well known, see e.g. \cite{Kato}, that $H_C$ is self-adjoint 
  on the Sobolev space $H^2 (\mathbb{R}^3)=D(-\Delta)$,  the domain of $-\Delta$
 (cf. also \cite{Kato}, p. 303). 
The
spectrum of $H_C$ consists of
isolated eigenvalues $E_n=-b^2/4n^2$, with multiplicity $n^2$, and an absolutely
continuous part, $[0,\infty)$.  

\z \subsection{Setting}\label{setting} Our starting point is the time evolution of the wave function $\psi(t,x)$ of the Hydrogen atom  described
by the Schr\"odinger
equation 
\begin{multline}
\label{sch-eq}
 i\psi_t=H_C\psi +V(t,x)\psi;\ \  \ \ \psi(0,x)=\psi_0(x)\in H^2(\RR^3) \\
{\rm where} ~~V(t, x) = \sum_{j\in \ZZ} \Omega_j (x) e^{i j \omega t}  \text{\ is real valued and $\Omega_0\equiv 0$}
\end{multline}
The operator $H_C+V(t,x)$ satisfies the assumptions of Theorem X.71,
p. 290, in \cite{Reed-Simon} v.2.; Theorem X.70, p. 285 also applies
in our setting. Thus, for any $t$, $\psi(t,\cdot)\in H^2
(\mathbb{R}^3)$, and the unitary propagator $U(t)$ for (\ref{sch-eq})
is strongly differentiable in $t$; see \S\ref{lapts1} for a short proof
in our case.

\begin{Assumption}\label{AA1}
  The $ \Omega_j (x) $, $j\in\ZZ$ 
are smooth inside a  common compact support, 
chosen  without loss of generality
to be the ball ${\sf B}_1 \subset \mathbb{R}^3 $ of radius $1$,
and 
$ \sum_{j \in \mathbb{Z}} (1+|j|) \| \Omega_j \|_{L^\infty ({\sf B}_1 ) } 
< \infty $.
\end{Assumption}

 \z \subsection{ Ionization} We say that the system  {\em ionizes} if the probability
to find the particle in any compact set vanishes for large $t$, {\em i.e.}, 
for any $a > 0 $ we have 
\begin{equation}
    \label{eq:ioniz}
P(t,{\sf B_{a}})=    \int_{{\sf B_{a}}}|\psi(t,x)|^2dx\to 0\ \ \ {as} \ \ \ t\to\infty
  \end{equation}
where ${\sf B_{a}}=\{x:|x|<a\}$. To prove ionization, it clearly 
suffices to prove (\ref{eq:ioniz}) for all  $a > 1$.

\z 
A simple way in which ionization may fail  is the existence
of a solution of the Schr\"odinger equation in the form
\begin{multline}
  \label{eq:persol}
 \psi(t,x)= e^{i\phi t} v (t,x)\ \text{with
$\phi\in\RR$}\text{ and $v \in L^2([0,2\pi/\omega] \times \RR^3)$ time-periodic.}
\end{multline}
Substitution in (\ref{sch-eq})  leads to the 
equation:
\begin{equation}
  \label{eq:floquet1}
  Kv =\phi v 
\end{equation}
where 
\begin{equation}
  \label{eq:floquet2}
 K= i\frac{\partial}{\partial t}-\Big(-\Delta-br^{-1}+V(t,x)\Big)
\end{equation}
is the Floquet operator,  densely
defined on $L^2([0,2\pi/\omega] \times \RR^3 )$;  
$0\ne v\in L^2$  implies by definition that $\phi\in \sigma_p(K)$, the
point spectrum of $K$.

Somewhat surprisingly, in all studied systems, $\sigma_p(K)\ne
\emptyset$ is in fact the {\em only} possibility for ionization to
fail. As we will show this is also true for (\ref{sch-eq}).  The proof
of ionization also implies that $K$ does not have any singular
continuous spectrum. This turns out to be a consequence of the
existence of an underlying compact operator formulation, the operator
being closely related to $K$.  Generic ionization is then expected
since $L^2$ solutions of the Schr\"odinger equation of the special
form (\ref{eq:persol}) are unlikely.  We prove that for $V(t,x)=2
\Omega(r)\sin(\omega t -\theta) $, $\Omega>0$ on $[0,1]$ and
sufficiently smooth, they do not exist.

\subsection{Laplace space formulation}\label{St1}
For $\psi\in H^2(\RR^3)$, 
the Laplace transform
\begin{displaymath}
  \label{eq:Lap}
  \hat{\psi}(p,\cdot):=\int_0^{\infty}\psi(t,\cdot)e^{-pt}dt
\end{displaymath}
exists for $p\in\mathbb{H}$,  \index{$p $} \index{bfhh@$\mathbb{H} $} 
the right half complex plane, and the map 
$p\to \hat \psi(p,\cdot)$ is $H^2$ valued
analytic in $\Re ~p > 0$. 
The Laplace transform converts the asymptotic
problem (\ref{eq:ioniz}) into an analytical one.

To improve the decay in $p$ of the Laplace transform, it is convenient to 
write
\begin{equation}
\label{psidecomp}
\psi (t,x) = \psi_0 (x)  e^{-t} + y (t,x)
\end{equation} 
Now, $y(t,x)$ satisfies
\begin{equation}
\label{sch-eq-y}
 i y_t - H_C y - V(t, x) y = 
e^{-t} \left [ i \psi_0 + H_C \psi_0 + V(t,x) \psi_0 
\right ] \equiv -y^0  (t,x) ; \ \ \ \ \ y(0,x)= 0 
\end{equation}
Standard arguments (see Appendix \ref{lapts}) show that
the $t-$Laplace transform of $y$,  $\hat{y}$ is in $H^2$ and satisfies
\begin{equation}
  \label{eq:S-lt-y}
(H_C-ip)\hat{y}(p,x)= {\hat y}^{0} (p, x) -\sum_{j\in\ZZ}
\Omega_j(x)\hat{y}(p-i j\omega,x)
\end{equation}
where \index{${\hat y}^{0} $}
\begin{equation}
\label{eq:haty0}
{\hat y}^{0} (p, x) =
-\frac{1}{1+p}\Bigg( i \psi_0 + H_C \psi_0\Bigg)
- \sum_{j\in\ZZ}  \frac{\Omega_j (x) \psi_0 (x) }{1+p-i j \omega}  
\end{equation} 
 
\subsection{The homogeneous equation and   the PDE-difference equations}  The
homogeneous system associated to
(\ref{eq:S-lt-y}) is 
\begin{equation}
  \label{hmdiff00}
  (-\Delta-b/r-ip) w (p,x)=-\sum_{j\in\ZZ}\Omega_j(x)w (p-ij\omega,x)
\end{equation} 
\begin{Note}\label{N1} {\rm (i) Clearly, (\ref{eq:S-lt-y}) and
    (\ref{hmdiff00}) couple two values of $p$ only if $(p_1-p_2)\in
    i\omega\ZZ$, and are effectively infinite systems of partial
    differential equations.  
    Setting \index{$p_1$}
\begin{equation}
  \label{defnp}
  p=p_1+in\omega, \text{ with}
  \ p_1\in \CC \text{\,mod\,} (i \omega )
\end{equation}
 we denote
$$y_n (p_1, x)=  \hat{y} (p_1+i n\omega,x), \\ y^{0}_n (p_1, x)
= {\hat y}^0 (p_1 + i n \omega, x) , \\  
w_n (p_1, x)= w (p_1+in\omega,x)$$ 
Eqns.  (\ref{eq:S-lt-y}) and (\ref{hmdiff00}) 
now become \index{yn@$y_n$} \index{$y_n^0$} \index{$w_n$} 
\begin{eqnarray}
  \label{eq:S-lt2}
(H_C-ip_1+n\omega) y_n=y^0_n -\sum_{j\in\ZZ}\Omega_j(x)
y_{n-j}  \\\label{eqvv}
(H_C-ip_1+n\omega) w_n=-\sum_{j\in\ZZ} \Omega_j (x) w_{n-j}
\end{eqnarray}}
\end{Note}
\begin{Note}\label{NN2}{\rm 
  Seen as a differential difference equation, the solution
  $\hat{y} (p,x)$ is then a vector
$\{ y_n (p_1,x) \}_{n\in\ZZ}$ and the whole problem depends only
parametrically on $p_1$. We have
\begin{equation}
  \label{eq:identif}
  y_n (p_1+i\omega,x)=y_{n+1}(p_1,x)
\end{equation}
and the analysis can be restricted to 
$$\mathbb{S}_0=\{p\in\overline{\mathbb{H}}:\Im\,p\in[0,\omega)\}$$
 where
$\overline{\mathbb{H}}$ \index{hbar@$\overline{\mathbb{H}}$} is the closure of ${\mathbb{H}}$.
There is
  arbitrariness in the choice of $ \mathbb{S}_0$ and,  to see
  analyticity in $p_1$ on $\partial\mathbb{S}_0$, it is convenient  to
  allow $p_1 \in \overline{\mathbb{H}}$, using (\ref{eq:identif}) to
  identify different strips of width $\omega$.}
 \end{Note}

  \section{Main results}

\begin{Theorem}\label{T3}
   Assume $V(t,x)=2 \Omega(r)\sin (\omega t -\theta)$, with $\Omega(r)=0$ for $r>1$, 
  $\Omega(r)>0$ for $r\le 1$ and $\Omega(r)\in C^\infty[0,1]$.
   Then
  $\sigma_p(K)=\emptyset$ and ionization always occurs. Furthermore,
  if $\psi_0 (x)$ is compactly supported and has only finitely many spherical harmonics, then
  $P(t,{\sf B_a})=O(t^{-5/3})$.
  \end{Theorem}
For the proof, given in \S\ref{PfT3}, \S\ref{asymptvn} and
\S\ref{Fredh},  we develop a relatively general rigorous WKB theory for infinite systems of differential equations. This yields  the asymptotic behavior of $w_n$
as $n\to-\infty$. The argument relies 
on Theorems \ref{cdi} and \ref{cdi(ii)} below.

\begin{Remark}
  { \rm The condition $\Omega(1^-)\ne 0$ simplifies the arguments but
    these could accommodate an algebraically vanishing $\Omega$.  
   (We also note  that  some one-dimensional models with rough $\Omega$ 
such as a  $\delta$ mass show failure of ionization, 
    see \cite{CMP221}, \cite{[6a]}
    and \cite{RCL}.) 

  }
\end{Remark}

\bigskip

 We will later derive 
   equivalent systems of integral equations, (\ref{eq:S-lt0230w}),  allowing
for a compact
   operator reformulation of the problem.

 \begin{Theorem}\label{cdi}
   In the setting \S\ref{setting},
assuming spherical
symmetry  in $x$ of the  forcing $V(t, x)$, ionization occurs {\em iff} for all
$p_1 \in\overline{\HH}$, (\ref{eqvv}) has only zero $H^2$ solutions
 decaying in $n$\footnote {For precise conditions,  see \S \ref{defC}
below and  the integral form (\ref{eq:S-lt0230w}).}.  This is
   true {\em iff} $\sigma_p(K)=\emptyset$.
\end{Theorem}

This extends results about absence of singular continuous spectrum of $K$,
\cite{Yajima}, to this class of systems, with Coulombic potential and
nonanalytic forcing.

 The proof is
given in \S\ref{Ab} and \S\ref{pcdi}. 

\bigskip
{\bf Properties of Floquet bound states for general compactly supported 
$V(t,x)$}. 
\addtocounter{footnote}{-1}
\begin{Theorem}\label{cdi(ii)}
  If there exists an $H^2$ nonzero solution $w$
of  (\ref{eqvv}) decaying in $n$, \footnotemark,
  then it has the
  further property
  \begin{equation}
    \label{I2}
   w_n=\bchi_{\mathsf B_1} w_n\ \ \text{for all}\ \ n<0
\end{equation}
with $\bchi_{A}$ the characteristic function of the set $A$.
\end{Theorem}
\z The general idea  of the proof is explained in \S~\ref{Sec7} and 
the details are given in \S~\ref{iproof}. 
\begin{Note}\label{N07}{\rm (i) 
{{
The Sobolev embedding theorem implies} that $w_n$ is continuous in
      $x$. From (\ref{eqvv}), $w_n$ is piecewise $C^2$,
      implying continuity of $\nabla w_n$ up to $\partial {\sf B}_1$.}

\noindent (ii)
Equation (\ref{I2}) makes the  second order system (\ref{eqvv})
    formally overdetermined since the regularity of $w$ in $x$ imposes both
    Dirichlet and Neumann conditions on $\partial {\sf B}_1$ for $n<0$. 
    Nontrivial
    solutions are not, in general, expected to exist.}\end{Note}

\section{Proofs}

{\bf Outline of the ideas.} 
 As in our previous work \cite{UAB}--\cite{JSPLieb}, 
summarized in \cite{JMP} 
on simpler systems,
we rely on a modified Fredholm theory to prove a 
dichotomy:
there are bound Floquet states, or the system gets ionized. 
Mathematically the Coulomb potential introduces a number of
substantial difficulties compared to the
potentials considered before (for references, see
  e.g. \cite{JSPLieb}), due to its singular behavior at the origin
and, more importantly, its very slow decay at infinity.

The slow decay translates into potential-specific corrections at
infinity, and standard general methods to show compactness in weighted
spaces of the Floquet resolvent, such as those in \cite{Yajima} or our
previous ones do not apply. Instead, the asymptotic behavior in the
far field of the resolvent has to be calculated in detail. 
The accumulation of eigenvalues of increasing multiplicity at the top
of the discrete spectrum of $H_C$ produces an essential
singularity at zero of the Floquet resolvent with a local expansion
of the form $\sum_{j,l} p^{l/2}e^{i j A (-ip)^{-1/2}}$ for
small $p$ when $\Re \, p \ge 0$,  where
$A= \pi b/2 $.  For sufficiently rapidly decaying potentials
  the exponentials would be absent. Their presence  clearly makes the analysis at $p=0$ of the Floquet
  resolvent  more delicate and is responsible for the change
  in the large time asymptotic behavior of the wave function, from
  $t^{-3/2}$ to $t^{-5/6}$.

  We introduce an extended parameter $X=(p^{1/2},
  e^{i A(-ip)^{-1/2}}$) and prove analyticity of the solution ${\hat
    y} $ in $X$, whose $p$-counterpart is $p$ small, 
$\Re\,p \ge 0$, and  similarly in regions near the special points 
  $p\in i \omega \mathbb{Z}$. We reformulate the problem in terms of an
  integral operator $\mathfrak{C}$, defined in \S\ref{defC}, closely
  related to the Floquet resolvent, shown to be compact in a suitable
  space and analytic in a variable corresponding to $X$.  Then, by the
  Fredholm alternative, $(I-\mathfrak{C})^{-1}$ is meromorphic, and
  in fact analytic in $X$, 
  since we show absence  of eigenvalues
of $I-\mathfrak{C}$ for {any $p \in \overline{\mathbb{H}}$. 
\subsection{The Hilbert space $\mathcal{H}$} \label{Hsp}  
\index{v@$\mathcal{H}$} Let $\mathcal{H}$ be  the Hilbert
space of sequences $Y = \{y_n\}_{n\in\ZZ},\,y_n\in L^2({\sf B}_a)$, with
$a>1$, and  with
$$\|Y\|^2:=\|Y\|^2_a
=\sum_{n\in\ZZ}(1+|n|)^{4/3} \|y_n\|^2_{L^2({\sf B_a})}<\infty$$
  \begin{Note}\label{reg1}
The properties of 
$(I-\mathfrak{C})^{-1}$ as $\Re \, p_1 \rightarrow 0^+$ 
ensure that 
$Y (p_1, \cdot) \in \mathcal{H}\cup H^2$ and is locally integrable in 
$p_1$  along $i\RR$.
  \end{Note}

  We then extend the stationary phase method to such a setting,
  cf. \S\ref{Newsp}, to evaluate, asymptotically for large $t$,
the inverse Laplace transform of $\hat{y}$ on $i \mathbb{R}$ and  
obtain the ionization result and time decay
estimates.

To show ionization we then
have to rule out the existence of a point spectrum of the Floquet
operator, that is the existence of nontrivial solutions of 
(\ref{eqvv}). We use the general criterion in Theorem \ref{cdi(ii)}
to  show that, if 
there exists a nonzero solution to (\ref{eqvv}), then a 
subsequence of $\{ w_n \}_{n\in-\NN} $ would be  singular at $x=0$, in contradiction with Note \ref{N07}, (i).

To find the behavior of solutions for large $n$, we develop a
WKB theory for infinite systems of ODEs and find the asymptotic
behavior  of $w_n$ in $n$ in detail. 
The formal WKB calculation of the behavior
is straightforward
algebra, relatively easy even in much more general settings, see
\S\ref{Hr}. Justifying the procedure is however delicate, and a
good part of the paper is devoted to that; cf. \S\ref{prf91},
\S\ref{Addtl}.

 The procedure of introducing an enlarged set of parameters with
  respect to which the solution is regular, when this does not hold in
  the original parameter, should also  be applicable to other problems where
  complicated singularities arise.
\bigskip

\subsection{ Proof of Theorem~\ref{cdi}} \label{Ab} We show that
$\hat{y}$ has a limit in  $L^1_{loc}$  on $\partial \mathbb{H}=i\RR$
   where it is smooth
  except for  possible poles and a discrete set of
  essential (but $L^1$) singularities. 
  Poles are present iff
  the integral form (\ref{eq:S-lt0230w}) of
  (\ref{hmdiff00}) 
  has nontrivial solutions in $\mathcal{H}$.  
  There is sufficient decay in $p$
  at infinity, so that, when poles are absent,  the Riemann-Lebesgue lemma
  applies, implying that $y$ decays as $t\to\infty$, proving
  ionization --since $\psi_0 (x) e^{-t}$ obviously goes to zero in this
  limit.  More detailed analysis of the resolvent reveals the nature
  of the essential singularity at $p = i \omega
  \mathbb{Z}$. Stationary phase analysis shows a $t^{-5/6} $ decay of
  the wave function if the initial condition is spatially compactly supported
  and contains only a finite number
  of spherical harmonics.

\begin{Proposition}\label{PP3}
  Ionization holds for every $\psi_0\in L^2$ iff it holds  for
  any $\psi_0$ in a set densely spanning
  $L^2$.
 
\end{Proposition}
\begin{proof}
  We make use of the standard triangle inequality argument to estimate 
 $U(t)\psi_0$, where $U(t)$ is the 
unitary operator associated to the Schr\"odinger evolution (\ref{sch-eq}).
\end{proof}

We choose $\psi_0$ in  a dense set $C_c^{\infty}(\RR^3)$, the smooth, compactly
supported functions in $\RR^3$.
Define as usual the angular momentum operators
$$-\mathbf{L}^2=\frac{\partial^2}{\partial\theta^2}+
\frac{1}{\tan\theta}\frac{\partial}{\partial\theta}+\frac{1}{\sin^2\theta}\frac{\partial^2}{\partial\phi^2}$$
  and
  $$L_z=-i\frac{\partial}{\partial\phi}$$
   Let  \index{$ P_{l,m} $}
 $ P_{l,m}$ be the orthogonal projector on $\{\phi:
  \mathbf{L}^2\phi=l(l+1)\phi,\ {L}_z\phi =m\phi\}$ for some $m\in\ZZ$,
  $|m|\le l\in \NN\cup \{0\}$. 

  Since $\sum_{l,m}P_{l,m}=I$,  we can now assume 
  without loss of generality that  $ \psi_0\in
  P_{l,m}\Big(C_c^{\infty}(\RR^3)\big)$ if $l$ and $m$ are arbitrary.
  Likewise, if $P (t, {\sf B_a} )$ decays 
  like $t^{-5/3}$ when $\psi_0\in
P_{l,m}\Big(C_c^{\infty}(\RR^3)\big)$, then the same decay rate clearly
holds for any $\psi_0$ given by a 
finite linear combination over $(l, m)$
(but not, 
in general, for any $\psi_0 \in L^2 (\mathbb{R}^3)) $.

\z {\bf Further notations.} An 
index of commonly used notations together with
their pages of definition
is given at the end of the paper. 
As usual we write 
$\mathbb{D}_{\epsilon}=\{z:|z|<\epsilon\}$, 
$\mathbb{D}=\mathbb{D}_1$  and we denote 
$\mathbb{D}^+_{\epsilon}=\mathbb{D}_{\epsilon}\cap
\left \{ z: \arg z \in \left ( -\pi/4, \pi/4 \right ) 
\right \} $. We also let
 $\mathcal{I}_{\epsilon} = i [-\epsilon, 0]$,
$\mathbb{H}^{+c}=\{p+c:p\in\HH\}$, $\ell_{\alpha}=
\{p:\Re(p) \ge 0,\Im (p)=\alpha \}$
and for a set $A$, $A_{\setminus \ell_\alpha}=A\setminus \ell_\alpha$. 
We denote 
$\mathcal{D}=\mathbb{H}\cup i\RR^+$, 
$\mathcal{O}(\mathcal{D})$ will denote {\em some} small open 
neighborhood
of $\mathcal{D}$. 
\index{$\mathbb{D}_{\epsilon}$} \index{$\mathcal{O}(\mathcal{D})$} \index{$\mathcal{I}_{\epsilon}$} \index{$\mathbb{H}^{+c}$} \index{$\mathbb{D}_1 $} 
 \index{$\mathbb{D}^+_{\epsilon} $}
\subsection{Step 1. Compact operator reformulation}\label{S3.2}
To investigate the analytic properties of $\hat{\psi}$ it
is convenient to introduce a new operator $\mathcal{A}_\beta$ which is a complex perturbation 
of $H_C$, having no real eigenvalues. More precisely, define \index{$\mathcal{A}_{\beta} $}
\begin{equation}
  \label{eq:defA}
\mathcal{A}_{\beta} :=
  H_C-i\beta(p){\bchi_{{\sf B}_a}}(r)-ip; \ \ (\text{with the understanding that}\ \  \mathcal{A}_0=H_C-ip)
\end{equation}  where $a>1$ and
\begin{equation}
  \label{eq:S-lt0239}
\beta={\beta}(p)=\left\{\begin{array}{ll}c>0 \text{ if } \Im\, p\in [-\epsilon, p_c]\text{ and } \Re\,p\ge 0
\\ 0\ \text{otherwise}
\end{array}\right.
\end{equation}
\z  \index{$\beta $} Here $\epsilon < \omega/2$ is small  
  as required in  Proposition \ref{P131} below,  and 
we choose $p_c$ so that  \index{$ p_c$}
$p_c/\omega \notin \ZZ$ and $ p_c>-E_0=b^2/4$, the ground state 
energy of the unperturbed 
  atom.  
Clearly $\mathcal{A}_\beta$ is defined on $D(H_0)$ and
$\mathcal{A}_\beta^*=\mathcal{A}_{-\beta}+ip^* + ip $. 
We  rewrite
(\ref{eq:S-lt2}) and (\ref{eqvv}) in the equivalent form
\begin{eqnarray}
  \label{eq:S-lt222}\ \ \ \ \ 
\left [H_C-ip_1+n\omega-i\beta(p){\bchi_{{\sf B}_a}}(r) \right ] y_n=
y^0_n -i\beta(p) {\bchi_{{\sf B}_a}}(r)y_n-\sum_{j\in\ZZ}\Omega_j(x)
y_{n-j}  \\\label{eqvvvv}\ \ \ \ \ \ \
\left [H_C-ip_1+n\omega-i\beta(p){\bchi_{{\sf B}_a}}(r) 
\right ] w_n
=-i\beta(p) {\bchi_{{\sf B}_a}}(r)w_n-\sum_{j\in\ZZ} \Omega_j (x) w_{n-j}
\end{eqnarray}
We show next that $\mathcal{A}_{\beta}^{-1}$ is analytic
in $p\in\HH\setminus \{ \ell_{p_c} \cup \ell_{-\epsilon} \} $, 
and sufficiently regular on
$i\RR$. Since the parameter $p_c$ is artificial, the non-analyticity at 
$\ell_{p_c} \cup \ell_{-\epsilon}$ of $\mathcal{A}_\beta^{-1}$
is not reflected in the actual solution ${\hat y} $, as discussed in Note~\ref{ancon}.

\begin{Proposition}\label{P4} There exists
 an open 
 neighborhood $\mathcal{O} $ 
of $\mathcal{D}\setminus \left \{ \ell_{p_c} \cup \ell_{-\epsilon} \right \}$,
not containing the  origin $0$, 
such that the operator  \index{$\mathfrak{R}_\beta $}
$\mathfrak{R}_\beta=\mathcal{A}_\beta^{-1}$ exists and is
analytic in $p\in \mathcal{O}\setminus (\ell_{p_c} \cup
\ell_{-\epsilon}) $. Furthermore, for any $p$ for which $\mathfrak{R}_\beta$ exists, we have
$\mathfrak{R}_\beta : L^2 (\mathbb{R}^3 ) \rightarrow H^2 \left (\mathbb{R}^3 \right )$.
\end{Proposition}
\z The proof is given in \S\ref{Sp4}.

\subsection{Restriction to a ball ${\sf B}$; Definition of $\mathfrak{C}$}
\label{defC}
To study ionization, we only need to know $y(t,x)$ for $x$ in a
fixed (but arbitrary) ball ${\sf B}_a \supset {\sf B}_1$. 
Henceforth, to simplify the notation,
we write
${\sf B}_a={\sf B}$. 
 We shall
therefore need to study the properties of ${\bchi_{\sf
    B}}\mathfrak{R}_\beta{\bchi_{\sf B}}.$ This sandwiched operator
(which preserves information about $L^2(\RR^3)$ through built-in
boundary conditions on $\partial {\sf B}$) is the one that we shall
most often use below. We recall that $p = p_1 + i n \omega$ and
$\hat{y} (p_1 + i n \omega, x) = {y}_n (p_1, x)$.  Since
$\Omega_j( x)={\bchi_{\sf B}}\Omega_j(x)$, 
(\ref{eq:S-lt2}) implies 
that for $x\in {\sf B}$,
\begin{equation}
  \label{eq:S-lt0230}
y_n (p_1,x)=
{\bchi_{\sf B}}\mathfrak{R}_{{\beta}} y^0_n + 
{\bchi_{\sf B}}
\mathfrak{R}_\beta{\bchi_{\sf B}}\Big[-i{{\beta}}
y_n (p_1, x)-\sum_{j\in\ZZ}\Omega_j(x)
y_{n-j} (p_1, x) \Big],
\end{equation}
where we may assume that ${\sf B}$ contains the support of $\psi_0 (x)$,  
and
therefore of $y^0_n$.
Note that $\mathfrak{R}_\beta$ depends on $n$ through
$p = i n \omega + p_1$.
Corresponding to (\ref{eq:S-lt0230}), we obtain the
homogeneous system:
\begin{equation}
  \label{eq:S-lt0230w}
w_n (p_1,x)=
{\bchi_{\sf B}}
\mathfrak{R}_\beta{\bchi_{\sf B}}\Big[-i{{\beta}}
w_n (p_1, x)-\sum_{j\in\ZZ}\Omega_j(x)
w_{n-j} (p_1, x) \Big],
\end{equation}

The elements of $\mathcal{H}$ will be denoted by capital letters \index{$Y $}
\index{$ \mathfrak{T}$}\index{$\mathfrak{C} $}
{\em e.g.}
$\left \{ y_n \right \}_{n \in \ZZ}=:Y $, $\left \{ y^0_n \right \}_{n
\in \ZZ} =: Y_0 $. We define the operators $\mathfrak{T}$ on 
$L^2({\sf B})$ by
$$\left \{ \mathfrak{T} Y_0 \right \}_n =  
{\bchi_{\sf B}} 
\mathfrak{R}_{\beta} y^0_n $$ 
and 
$\mathfrak{C} $ on $\mathcal{H}$ by
$$ \left \{ \mathfrak{C} Y \right \}_n =
{\bchi_{\sf B}}
\mathfrak{R}_\beta{\bchi_{\sf B}}\Big[ -i {\beta}
y_n (p_1, x)-\sum_{j\in\ZZ}\Omega_j(x)
y_{n-j} (p_1, x) \Big] $$
Then,  we rewrite 
(\ref{eq:S-lt0230}) in the form
\begin{equation}
  \label{eq:intform}
   {Y}=\mathfrak{T} {Y}_0+\mathfrak{C}{Y}
 \end{equation}
 \begin{Note}\label{ancon}{\rm 
   We shall see that for any $\beta$ satisfying
  (\ref{eq:S-lt0239}), eq.
   (\ref{eq:intform}) has a unique 
solution in $\mathbb{H}$ (call it now $Y^{(\beta)}$). 
Thus, away from the artificial
cuts, all these solutions coincide (since the domain 
of $Y^{(0)}$ corresponding to  $c=0$, contains all of the others). 
Hence, wherever some $Y^{(\beta)}$ 
has analytic continuation,
so will $Y^{(0)}$. 
}
 \end{Note}
 The homogeneous system corresponding to (\ref{eq:intform}) is given by
\begin{equation}
\label{eq:chomo}
w = \mathfrak{C} w 
\end{equation} 

\begin{Note}\label{clarif}{We have shown (cf. \S\ref{St1} and 
\S\ref{lapts}) that
      $\hat{\psi}$ and the Laplace transform ${\hat y} (p, \cdot) =
      \mathcal{L} \left (\psi - \psi_0 e^{-t} \right ) $ \index{$\mathcal{L}$} exist
      for $\Re \, p > 0$.  The corresponding $Y = \left \{ y_n
      \right \}$,  $y_n = {\hat y} (in \omega + p_1,\cdot) $,
      restricted to ${\sf B}$, will therefore satisfy
      (\ref{eq:intform}) for $\beta =0$ when $\Re\, p_1 > 0$.  It will be
      shown that (\ref{eq:intform}) has a unique solution $Y \in
      \mathcal{H} $ for any $\Re \, p_1 \ge 0 $, and that $Y$ is
      analytic in $p_1 \in \mathbb{H}$ and has an $L^1_{loc}$ limit on
      $i\RR$, with sufficient decay in $n$.  
     The implied decay and regularity properties of 
     ${\hat y} (p, \cdot)$ on $i\RR$ show that $
      \mathcal{L}^{-1} {\hat y} + e^{-t} \psi_0 $ 
      (the integration 
contour taken to be $i\RR$)
equals $\psi$ for $x\in {\sf B}$.}
\end{Note}

\begin{Proposition}[Asymptotic behavior of
${\bchi_{\sf B}}\mathfrak{R}_{\beta}{\bchi_{\sf B}}$]\label{P12}
If $\Re\, p=0$ and   $| \Im\, p | \to\infty$  (see Note \ref{NN2}), then 
$\|{\bchi_{\sf B}} \mathfrak{R}_{\beta}{\bchi_{\sf B}}\|$ 
$=O(|p|^{-1/2} )$ (recall that
$\beta=0$ if $|\Im ~p|$ is large).
Moreover, for any $\epsilon > 0$, 
$\bchi_B \mathfrak{R}_{\beta} \bchi_B $ is analytic in $p$
in an open set containing $ -i (\epsilon, \infty)$.
\end{Proposition}

\z  For $\Im ~p \rightarrow +\infty$, the $|p|^{-1/2}$ decay rate follows from
the spectral theorem since we are outside the spectrum, while for $\Im ~p \rightarrow -\infty$,
the  rate is obtained using  Mourre estimates \cite{Jensen1989}, Theorem 6.1. 
The rest
of the proof involving analyticity is given in \S\ref{P12P1} and 
relies on an explicit representation of 
the resolvent for $H_C$, see \S\ref{Green}.
Using spherical symmetry,
the explicit Green's function could be avoided, but in view of possible
future
generalizations to non-spherical $V(t,x)$, 
we prefer this more delicate approach.

\begin{Lemma}
\label{lemY0C}
We have  $\mathfrak{T} Y_0 \in \mathcal{H}$. The operators $\mathfrak{S}:=Y\to \{\sum_{j\in\ZZ}\Omega_j(x)
y_{n-j} (p_1, x)\}_{n\in\ZZ}$ \index{$\mathfrak{S}$} and  
$\mathfrak{C}$ are bounded in $\mathcal{H}$.
\end{Lemma}
\begin{proof} 
We note from Proposition \ref{P12} that 
$\mathfrak{R}_\beta = O\left (|p|^{-1/2} \right )$ for large
$p$, {\it i.e.} $O\left (|n|^{-1/2} \right ) $
for large $|n|$, since $p = i n \omega + p_1 $. 
Therefore, from the expression of $\hat{y}_n^0$ in (\ref{eq:haty0}),  
\begin{multline*}
\sum_{n \in \ZZ} (1+|n|)^{4/3} \| 
\left \{ \mathfrak{T} Y_0 \right \}_n \|^2 \\
\le  
\sum_{n \in \ZZ} (1+|n|)^{4/3} \left ( \| \bchi_{\sf B} \mathfrak{R}_{\beta} 
\bchi_{\sf B} \|  
\left [ \frac{\| \psi_0 \|_{L^2} + \| \Delta \psi_0 \|_{L^2} }{ 
|1+p_1 + i n \omega|} \right ] \right . \\
\left . + \| \bchi_{\sf B} \mathfrak{R}_{\beta} 
\bchi_{\sf B} \| \| \psi_0 \|_{L^2} 
\sum_{j \in \ZZ} \frac{\| \Omega_j \|_{L^2}}{
| 1 + p_1 + i (n-j) \omega |}  \right )^2  \\
\le C \left [ 1 + \sum_{n \in \ZZ} (1+|n|)^{6/7} 
\left ( \sum_{j \in \ZZ} \frac{\| \Omega_j \|_{L^2}}{(1+|n-j|)} \right )^2 
\right ]
\end{multline*}
Using
(7.12) and (7.13) in  \cite{JSPLieb} with $\gamma = 6/7$, the above
is finite since 
$$\sum_{j\in \ZZ} (1+|j|)^{3/7} 
\| \Omega_j \|_{L^2} \le \sum_{j \in \ZZ} (1+|j|) \| \Omega_j \|_{L^2}
< \infty$$ The proof that $\mathfrak{S}$ is the same as that of Lemma
27 of \cite{JSPLieb}, with
$\gamma=\frac{4}{3}$,  replacing absolute values by norms in $x$. Since
$\mathfrak{R}_{\beta}$ is uniformly bounded (in the
operator norm) and acts diagonally in $n$, $\mathfrak{C}$ is bounded too.
\end{proof}

\begin{Lemma}
\label{lemCanal}
Both $\mathfrak{T} Y_0$ and the operator $\mathfrak{C}$ are 
analytic in $p_1$ for
$$p_1 \in \mathcal{O}\left (\overline{\mathbb{H}} 
\setminus \left(\{ \ell_{p_c} + i \omega \ZZ \} \cup \{ \ell_{-\epsilon} + i 
\omega \ZZ \} \cup
\{ \mathcal{I}_\epsilon + i \omega {\ZZ}\} \right)\right )$$
\end{Lemma}
\begin{proof} Propositions 
\ref{P4} and \ref{P12} imply that  $\mathfrak{R}_\beta$ is
 analytic in $ p \in \mathcal{O} \setminus \{ \ell_{p_c}
  \cup \ell_{-\epsilon} \} $ and in an open set containing
$-i(\epsilon,\infty)$. Analyticity of $\mathfrak{T} Y_0 $ and
$\mathfrak{C}$ follow from their definition (we
note
$\mathfrak{C}$ is a
norm limit of analytic operators: 
its restrictions  to the subspaces 
with nonzero  components for $|n|\le N$ only). 
\end{proof}

\begin{Remark}
\label{remanal}
As shown later, 
$(I-\mathfrak{C})$ is invertible.
Since the solution $Y$ cannot depend on the arbitrary
parameters $\epsilon$ and $p_c$ (see Note \ref{ancon}), 
the non-analyticity of $\mathfrak{C}$
and $\mathfrak{T} Y_0$ for 
$$ p_1 \in \{ \ell_{p_c} + i \omega \ZZ \} \cup
\{ \ell_{-\epsilon}+i\omega \ZZ \} \cup \{ 
\mathcal{I}_\epsilon + i \omega \ZZ \} $$ 
is not
reflected in $Y$.
\end{Remark}

\begin{Proposition}\label{P010}
  For $\Re\, p_1>0$ large enough, (\ref{eq:intform}) has a unique solution
  in $\mathcal{H}$. The inverse Laplace transform in $p$ of 
   ${\hat y} (p, x) =: y_n (p_1, x)$, where
  $p= i n \omega + p_1$, 
solves the initial value problem (\ref{sch-eq-y})  in  ${\sf B}$
(see Note \ref{clarif}).
\end{Proposition}
\z The proof is given in \S\ref{010}.

\subsection{Step 2. Regularity of 
$\mathfrak{R}_{\beta,l,m}$ at $p = 0$ and of $\mathfrak{C}_{l,m}$ at $p_1=0$} 
$ $ \index{$ \mathfrak{R}_{\beta,l,m}$}\index{$ \mathfrak{C}_{l,m}$}

Define $\mathfrak{R}_{\beta, l,m} = P_{l,m} \mathfrak{R}_\beta $ \index{$ \mathfrak{C}_{l,m}$}
and $\mathfrak{C}_{l,m} = P_{l,m} \mathfrak{C}$.

\begin{Note}[Compactness versus regularity of $\mathfrak{R}_{\beta, l, m}$]
\label{N3}
  {\rm The term $-i\beta{\bchi_{\sf B}}$ was introduced in \S\ref{Ab} to ensure
    that $\mathfrak{R}_{\beta, l, m}$ is bounded in
    $\overline{\mathbb{H}}$. Since $-i{\beta}{\bchi_{\sf B}}$ is localized in
    $x$, the shifts in the poles created by the  point
spectrum of
    $H_C$ are smaller as $p\to 0$ (the size of the orbitals of the
    Hydrogen atom grows when the energy approaches zero.) The resulting
    integral operators have an essential singularity at $p=0$. The factor
    ${\bchi_{\sf B}}$ is needed to ensure compactness,  
    simplifying the analysis.}
\end{Note}
The poles of the resolvent $\mathfrak{R}_{\beta, l,m}$
    accumulate at $p=0$ from $-\mathbb{H}$, 
     along a curve tangent to the positive imaginary $p$-axis 
     {(see Note \ref{NoteAcal} in \S \ref{Acalc})}.  
    As a result, while being uniformly bounded, 
    $\mathfrak{R}_{\beta, l,m}$ is not
    continuous along the imaginary $p$ line at zero 
    but oscillates without
    limit. Boundedness of 
    ${\bchi_{\sf B}}\mathfrak{R}_{\beta, l, m}{\bchi_{\sf B}}$ (which is not
    difficult to prove)
   does not ensure boundedness of the solution $Y$.  
    However, we  do have analyticity in an extended, two-dimensional, parameter.
    Let  $\lambda:=\sqrt{-ip}$ (with the usual branch of the square
root,  $\Im\lambda<0$ if $p\in\mathbb{H}$) and let
    $X:=(p^{1/2}, Z)$ \index{$ X$} with 
    \begin{equation}
      \label{eq:eqZ}
      Z= e^{\frac{i\pi b}{2\lambda}}
    \end{equation}\index{$ Z$}
    (The dependence of $Z$ on $\lambda$ reflects
    the actual behavior of the solution.) The resolvent is analytic
    in $X$ and a useful Fredholm alternative can be applied.
 
For any $a > 1$  we can choose a $c$ in (\ref{eq:S-lt0239})
  (see \S\ref{pp131} below) such that the following
  statement holds.
\begin{Proposition}[Analyticity in $X$]\label{P131}
  ${\bchi_{\sf B}} \mathfrak{R}_{{\beta,l,m}}{\bchi_{\sf B}} $ 
   is  analytic with respect to $X$
   on the compact 
set  $\overline{\mathbb{D}^+_\epsilon}\times \overline{\mathbb{D}}$. \footnote{ As usual, by analyticity
in a compact set, we mean analyticity in some open set containing the compact set. Analyticity in 
  $\overline{\mathbb{D}^+_\epsilon}\times \overline{\mathbb{D}}$ 
 of course, implies that  ${\bchi_{\sf B}} \mathfrak{R}_{{\beta,l,m}}{\bchi_{\sf B}} $  is 
given by a convergent double series in $p^{1/2}$ and $e^{\frac{i\pi b}{2\lambda}}$.}
\end{Proposition}
\z The proof is given in \S \ref{pp131}.

As a corollary, we have the following regularity
property of $\mathfrak{C}_{l,m}$ and $\mathfrak{T}_{l,m} Y_0 $.
Let
\begin{equation}
  \label{eq:Se}
  \mathbb{S}_{-\epsilon}=\{p:\Im\,p\in[-\epsilon,\omega -\epsilon), \ \ \Re \, p \ge 0 \}
\end{equation}\index{$\mathbb{S}_{-\epsilon} $}
\begin{Corollary}
\label{CorP131}
For $p_1 \in \mathbb{S}_{-\epsilon}$, 
define $X_1  = \left (p_1^{1/2}, Z_1 \right )$, \index{$X_1 $} where 
$Z_1= \exp \left [ i\pi b/\left ( 2\sqrt{-i p_1} \right ) \right ]$.\index{$Z_1 $}
Then, $\mathfrak{T}_{l,m} Y_0$ and  $\mathfrak{C}_{l,m}$ 
are analytic in $X_1$ on the 
compact set
$\overline{\mathbb{D}^+_\epsilon}\times \overline{\mathbb{D}} $.
\end{Corollary}
\begin{proof} Note first that Propositions \ref{P4} and \ref{P12} and 
the relative arbitrariness in the choice of $p_c$ and $\epsilon$ imply that
$\bchi_B \mathfrak{R}_\beta \bchi_B $ is analytic in $p$ in a neighborhood of
$p = i n \omega$, $n \in \mathbb{Z} \setminus \{0 \}$.
Since for large $|n|$,
$ \bchi_B \mathfrak{R}_\beta \bchi_B = \bchi_B \mathfrak{R}_0 \bchi_B $, its
expression as an 
integral operator involving $G$ in (\ref{eq:exprr}) (see Note (\ref{NoteG})
as well) 
implies
a lower bound of the analyticity radius independent of
$n$.
For sufficiently small $\epsilon$ for  
any $n \in \mathbb{Z}$, including $n=0$, then, 
analyticity of $\bchi_B \mathfrak{R}_\beta \bchi_B $ 
in the expanded variable 
$$\left ( \sqrt{p-in\omega},  
\exp \left [ \frac{i \pi b}{(2 \sqrt{ -i (p-in\omega) } )} \right ] \right )$$ 
follows
 in the domain  
$$\left\{|p-in \omega|^{1/2} \le \epsilon;
\left | \exp \left [ \frac{i \pi b}{(2 \sqrt{ -i (p-in\omega) } )} \right ] 
\right | \le 1\right\}$$
(since Proposition 
\ref{P131}
gives analyticity 
$X \in \overline{\mathbb{D}^+_\epsilon} \times\mathbb{D}$.)
Analyticity of $\mathfrak{C}_{l,m}$ in
$X_1  
\in \overline{\mathbb{D}^+_\epsilon}\times \overline{\mathbb{D}} 
$ 
now follows 
since $\mathfrak{C}_{l,m}$ is the 
    norm limit of analytic operators (the restrictions  of
    $\mathfrak{C}_{l,m}$ 
    to the subspaces of $\mathcal{H}$ with zero
    components for $|n|>N$). (See Proposition \ref{P12} for the necessary 
 estimates of decay in $N$.) The analyticity of $\mathfrak{T}_{l,m}
Y_0$ follows from its definition. 
\end{proof}

\subsection{Compactness}
\begin{Proposition}\label{Ccomp}
  $\mathfrak{C}_{l,m}$ is compact in $\mathcal{H}$ (cf. Note~\ref{N1}) for
$p_1 \in \overline{\mathbb{H}}$. 
\end{Proposition}
\z The proof is given in \S \ref{020}.

\subsection{Step 3. The Fredholm alternative}
We can now formulate the ionization condition using the Fredholm
alternative.

\begin{Proposition}\label{P200}
  If (\ref{eq:chomo})  has no
  nontrivial solution in $\mathcal{H}$ for  
  $p_1 \in \mathbb{S}_{-\epsilon}$, 
  then $\left (I-\mathfrak{C}_{l,m} \right )^{-1}$
  exists and 
  the system ionizes  (cf. (\ref{eq:ioniz})).
\end{Proposition}

The first part is simply  the Fredholm alternative. Ionization
  follows from the following proposition. We recall that
 ${\hat y} (p_1 + in \omega , x) = y_n$ are the components of $Y$.

\begin{Proposition}\label{P20}
  Assume (\ref{eq:chomo}) has no nontrivial solution when $p_1 \in
  \overline{\mathbb{H}}$. Then, for $\psi_0 \in P_{l,m} \left (
    C_0^{\infty} (\mathbb{R}^3) \right ) $, the solution $Y \in
  \mathcal{H}$ to (\ref{eq:intform}) is analytic in $p_1\in
  \overline{\HH} \setminus \{ i \omega \ZZ \} $ and analytic with
  respect to $X_1$ in $\overline{\mathbb{D}^+_\epsilon}\times
  \overline{\mathbb{D}}$.  In particular $Y$ is bounded at $p_1=0$.
  These properties imply sufficient regularity and decay of
    $\hat{y} (p, x)$ so that the integration contour in
    $\mathcal{L}^{-1}\hat{y}$  can be taken to be $i\RR$. By
      the Riemann-Lebesgue lemma, $P(t, {\sf B} ) \rightarrow 0$ as
      $t\rightarrow \infty$.
\end{Proposition}
\z The proof is given in \S\ref{Fredh}.

\subsection{End of proof of Theorem~\ref{cdi}} \label{pcdi} It only remains to 
make the connection with Floquet theory. This is 
done in \S~\ref{Fl3}.

\subsection{Proof of Theorem~\ref{cdi(ii)}.}\label{Sec7}
Eq. (\ref{eqvv}),  restricted to $\sf B$,
follows from the homogeneous system
$w=\mathfrak{C} w$.  Multiplying (\ref{eqvv}) by
$\overline{w_n}$, summing over $n$, and integrating over
 ${\sf B}_{\tilde a}$, where  ${\tilde a} \in (1, a]$, we are lead to
a nonnegative definite quantity involving
$w_n|_{\partial{\sf B}_{\tilde a}} $ being zero for $n < 0$.
Details are given in \S~\ref{iproof}.

\subsection{Proof of ionization for spherically symmetric
 $\Omega$, Theorem~\ref{T3}}\label{PfT3}
We consider the case $V(t,x)=2 \Omega(r)\sin\omega t$, corresponding
to $\Omega_1 = - i \Omega$ and $\Omega_{-1} = i \Omega$ ($\Omega$ is
real valued). The proof in the slightly more general case
$2 \Omega (r) \sin \left (\omega t - \theta \right )$ amounts to 
replacing $t$ by $t -{\theta}/{\omega}$ and $\psi_0 (x)$ by
$\psi \left ({\theta}/{\omega}, x \right )$ in our proof.  Recall $\mathfrak{C}_{l,m} = P_{l,m}\mathfrak{C}$
\footnote{As discussed, it suffices to show ionization on a dense
  subset of initial conditions.}.  We obtain by projection of
(\ref{eq:intform}) to $P_{l,m}\Big(L^2(\RR^3)\Big)$,
\begin{equation}
  \label{eq:intform4}
   Y=Y_0+\mathfrak{C}_{lm}Y
\end{equation}
The homogeneous equation associated
to (\ref{eq:intform4}) is
\begin{equation}\label{eqhomogen}
   w=\mathfrak{C}_{lm} w,\ w \in\mathcal{H}
\end{equation}
 The Fredholm alternative applies and (\ref{eq:intform4}) has a
unique solution in
$\mathcal{H}$ iff (\ref{eqhomogen}) implies $w=0$.
\begin{Note}{\rm 
  By separation of variables in spherical coordinates, we see that
  $\mathfrak{C}_{l,m}$ can be defined in the same way as
  $\mathfrak{C}$, replacing $\mathcal{A}_{\beta}$ in (\ref{eq:defA})
  by
\begin{equation}
  \label{eq:S-lt232}
\mathcal{A}_{\beta,r}=-\frac{d^2}{dr^2}-\frac{2}{r}\frac{d}{dr}+\frac{l(l+1)}{r^2}-\frac{b}{r}-ip_1+n\omega-i\beta{\bchi_{\sf B}}
\end{equation} \index{$ \mathcal{A}_{\beta,r}$}
and the associated differential-difference systems are obtained
by replacing $-\Delta-b/r$ with $\displaystyle -\frac{d^2}{dr^2}-\frac{2}{r}\frac{d}{dr}+\frac{l(l+1)}{r^2}-\frac{b}{r}$.
}
\end{Note}
Clearly  if there exists a nontrivial solution $w\in\mathcal{H}$ of
  (\ref{eqhomogen}), then, again by  elliptic regularity (see
Proposition \ref{P4}), 
$ v $  defined by \index{$ v $}\index{$\mathcal{Y}_{l,m} $} 
$w=\mathcal{Y}_{l,m} v (r) $, where $\mathcal{Y}_{l,m}$ are the spherical harmonics, is a nontrivial solution to
\begin{equation}
  \label{eqhomogen.1}
  \mathcal{A}_{\beta, r} v_n 
= - i \Omega 
\left ( v_{n+1} - v_{n-1} \right ) - i \beta \bchi_B v_n ~\ \,~~~~{\rm implying}~~   
\mathcal{A}_{0, r} v_n = - i \Omega \left ( v_{n+1} - v_{n-1} \right ) 
\end{equation}
\begin{Proposition}\label{startpoint}
  If $v$ satisfies (\ref{eqhomogen.1}), then there exists $n\ge 0$ such
  that either ${\bf (i)}$ $v_{n}(1)\ne 0$; or ${\bf (ii)}$ $v_{n} (1)
  =0$, but $v_{n}'(1)\ne 0$; let $n_0$ be the smallest such $n$.  By
  homogeneity, we can assume that $v_{n_0}(1)=1 $ in case ${\bf (i)}$
  and $v_{n_0}'(1)=-\sqrt{\Omega (1)}$ in case ${\bf (ii)}$ 
(we use the positivity of $\Omega$).
\end{Proposition}
\z The proof is given in \S\ref{Lemhom1}.

\z {\bf Definition:} 
We define $\tau $  \index{$ \tau$} to be  $0$ or $1$ in case
\textbf{(i)} and $1$ in case \textbf{(ii)} respectively. 

\smallskip

\subsection{Asymptotic behavior of  $v_n$ in
(\ref{eqhomogen.1}) as $n\to -\infty$}\label{asymptvn}
\z In view of Proposition~\ref{P20} we see that (\ref{eqhomogen.1})  holds
the necessary ionization information.
\subsubsection{Notation} \label{eqxi}
Let \index{$ \mathfrak{s}$}
\begin{equation}
\label{defmk}
\mathfrak{s}(r): =  \int_{r}^1 \sqrt{\Omega (\rho)} d\rho \ \ (r\in (0,1))
\end{equation}
 By assumption $\Omega>0$ is smooth
and then so is $
\mathfrak{s} $.
Let \index{$\alpha $} \index{$\zeta $}
\begin{equation}
  \label{eq:defquants}
n _0-k , \, \,  \Omega_0=\Omega(0), \,\, \Omega'_0=\Omega'(0), \,
\, {\mathfrak{s}}_{_0}=\mathfrak{s}(0), \,\, 
\alpha = \frac{2 \sqrt{\Omega_0}}{{\mathfrak{s}}_{_0}}, \, \, \zeta = 
\alpha k r 
\end{equation}
Denote \index{$H_0 (\zeta) $}
$$H_0 (\zeta): =
\sqrt{\frac{2}{\pi}} e^{\zeta} \zeta^{1/2} K_{l+1/2} (\zeta);\ \ 
G_0 (\zeta) = \sqrt{\frac{\pi}{2}} e^{\zeta} \zeta^{1/2}
I_{l+1/2} (\zeta)$$ where $K_{l+1/2}$ and $I_{l+1/2}$ are  the modified Bessel functions of order
$l+1/2$.  It
  follows that for small $\zeta$, $$H_0 (\zeta) ~\sim 2^{-l}
  \zeta^{-l}(2 l)!/l!$$ Let   $\check{H}
(\zeta;k,l)$ be the unique solution of the integral equation \index{$\check{H} (\zeta;k,l) $}
\begin{multline}\label{inteq1} 
\check{H} (\zeta;k,l) = G_0 (\zeta) \int_{0}^\zeta e^{-2s} G_0 (s)  
\mathcal{R} (H_0 + k^{-1}{\check{H}}) (s)ds \\- H_0 (\zeta) \int_{k \alpha}^\zeta
e^{-2s} H_0 (s)  
\mathcal{R} (H_0 + k^{-1}{\check{H}})  (s) ds
\end{multline}
for $\zeta \in [0, k \alpha]$, where 
the operator $\mathcal{R}$ is defined by
$$ \Big(\mathcal{R} f\Big)(\zeta) = 
2 \left( -\frac{\omega}{2\alpha^2} 
+\frac{\Omega^\prime_0(1+2\zeta)}{4\alpha \Omega_0} +\frac{\tau}{2} \right)f^\prime
-  \frac{b f}{\alpha \zeta} 
$$ 
 Define \index{$H(\zeta) $}  $$H(\zeta)=H(\zeta; k,l):= H_0
(\zeta) + k^{-1} \check{H} (\zeta; k,l)$$ It can be checked that $H 
$ satisfies
\begin{equation}
\label{S7.5}
H^{\prime \prime}=
2\left(1-\frac{\omega}{2k\alpha^2} 
+\frac{\Omega^\prime_0(1+2\zeta)}{4 k \alpha \Omega_0} +\frac{\tau}{2k} \right) H^\prime
+ \left(\frac{l (l+1)}{\zeta^2}- 
  \frac{b }{\alpha
    \zeta k} \right) H
\end{equation}
with the following asymptotic condition\footnote{As is common,
the notation $O\left (a, b \right ) \equiv O \left ( |a| + |b| \right )$;
similarly $O\left (a, b, c \right ) \equiv O \left ( |a| + |b| +|c| \right )$}
\begin{equation}
\label{S7.6}
H(\zeta) \sim 1 + \frac{l(l+1)}{2 \zeta} +
\frac{b}{2 k \alpha} \log \zeta +
O \left ( \frac{\log \zeta}{k\zeta}, \frac{1}{\zeta^2} \right ) (\zeta,
k\to\infty, \zeta\leqslant  k\alpha)
\end{equation} 

\begin{Remark}\label{R46}
(i)   $\displaystyle H(\zeta;k,l) \sim 
  \sqrt{\frac{2}{\pi}} e^{\zeta} \zeta^{1/2} K_{l+1/2} (\zeta)(1+o(1))\ \ \text{as} \ k\to\infty$

(ii)   From the expression (\ref{inteq1}) for $\check{H}$ it is seen that
  as $\zeta \rightarrow 0$ we have
  $\check{H} (\zeta;k,l) \sim {\rm const.}\zeta^{-l+1}$ for $l \ne 1$ and $\check{H}
  \sim {\rm const.} +{\rm const~}  \zeta \log \zeta$ for $l=1$.  
   For $\tau =$  0 or 1, $\check{H}$
  is less singular than $H_0$ at $\zeta = 0$.  
\end{Remark}
Define \index{$ m_k (r)$}
\begin{equation}
\label{S7.3}
m_k (r) =
\frac{\mathfrak{s}^{2 k+\tau}\Omega^{\frac{1}{4}} (1)H (\alpha k r )}{(2 k+\tau)! \Omega^{\frac{1}{4}}(r)H(\alpha k)}  \exp \left [ \frac{\omega }{4}
  \int_1^r \frac{\mathfrak{s} (s)ds }{\sqrt{\Omega(s)}} \right ]
:= \frac{\mathfrak{s}^{2 k+\tau}}{(2 k+\tau)!}F_k (r)
\end{equation}
\begin{Note}\label{Nuk}{\rm 
 From standard  properties of the modified Bessel
function $K_{l+1/2}$, it follows  that for large enough $k$, $H(\alpha k r)$ is
continuous and nonzero
for $r \in (0, 1]$ and that as $r \rightarrow 0$, $H(\alpha k r) $
is singular as $r^{-l}$. Therefore for any $k$ sufficiently large 
$u_k \equiv r^l m_k $ has a finite limit nonzero limit
as $r \rightarrow 0^+$.}
\end{Note}
\begin{Definition}
 With $n_0$ as in Proposition ~\ref{startpoint},
we define $h_k (r)$ by \index{$ h_k$}
\begin{equation}
\label{DefRk}
w_{n_0-k} = \frac{i^k}{r} m_k (r) h_k (r) \mathcal{Y}_{l,m}  
\end{equation}  
\end{Definition}

\begin{Theorem}[Behavior as $k\to +\infty$ 
({\it i.e.} $n \rightarrow -\infty$) ]
\label{Rnasympt}
 For any sufficiently large $k $,
$u_k := r^l m_k (r) $ is continuous in $r \in (0, 1]$ and
$u_k (r) \to const \ne 0$ as $r \rightarrow 0^+$. 
Furthermore, if there is a  nontrivial
solution to 
(\ref{eqhomogen.1}), 
then
there exists a subsequence $k_j\to \infty $ such that for any $r \in [0,1]$,
\begin{equation}
\lim_{j\rightarrow \infty} h_{k,j} (r) = 1  
\end{equation}
\end{Theorem}
\z  The first part follows  simply from  Note \ref{Nuk}. The rest of the proof 
is given in \S\ref{prf91}.

\begin{Proposition}~\label{P210}
There is no nonzero solution of (\ref{eqhomogen}) in $\mathcal{H}$.
\end{Proposition} \z Indeed, Theorem \ref{Rnasympt} shows that   otherwise
$(r^{l+1}v_n )(0)=m_n\ne 0$ for a subsequence of $n<0$. This implies
that the corresponding $w_{n} (x) \sim m_n  r^{-l-1} \mathcal{Y}_{l,m} $ for 
$r=|x| \rightarrow 0$. 
 This singularity is incompatible with $w_n \in H^2$, 
(see Proposition \ref{P4}).
Thus there is no admissible solution of the
homogeneous system and the first part of Theorem \ref{T3}
follows from Theorem~\ref{cdi} (i).  See
also the remarks in \S\ref{improvements}. The result on the decay rate
follows from the type of essential singularities of ${\hat y} $ for 
$p \in i \omega \ZZ$; see \S \ref{Fredh}.  

\section{Proofs of intermediate steps}
\subsection{Proof of Proposition~\ref{P4}}\label{Sp4} 

As mentioned in \S\ref{S3.2}, $\mathcal{A}_\beta$ and
$\mathcal{A}_\beta^*$ are adjoints of each-other. They
are furthermore densely defined and hence (see, {\em e.g} \cite{Kato},
Theorems 5.28 and 5.29, p. 168), closed.  Once we show that
$\mathcal{A}_{\beta}(p)$ is invertible in $\mathcal{D}$, analyticity
of $\mathcal{R}_\beta$ in $\mathcal{O}\setminus ({\ell_{p_c} \cup
  \ell_{-\epsilon}})$ follows (the spectrum of the closed operator
$H_C-i\beta{\bchi_{\sf B}}(r)$ is a closed set). Analyticity holds
wherever $\mathcal{A}_{\beta}$ is analytic, \cite{Reed-Simon}, vol.
1, Theorem VIII.2, p. 254).

\z {\bf(1)} {\em Eigenvalues}. We first show below that no $ip\in
i\mathcal{D}$ is an eigenvalue of $H_C-i\beta{\bchi_{\sf
    B}}(r)$. Assume we had $\mathcal{A}_\beta\psi =0$. If
  $\beta=0$, $\mathcal{A}_0\psi =0$ implies $ip\in\sigma_p(H_C)$,
  but, by construction, these values of $p$ correspond to the region
  where $\beta\ne 0$. So we can assume $\beta>0$. Then
\begin{equation}
  \label{eq:eqsp}
  \langle \psi,(-\Delta-ip-br^{-1})\psi\rangle +\langle \psi,-i\beta{\bchi_{\sf B}}\psi\rangle =0
\end{equation}
 Taking the imaginary part of
(\ref{eq:eqsp}) we get
\begin{equation}
  \label{eq:eqsym}
\Re\,p\,\, \langle \psi,\psi\rangle +c\langle \psi,{\bchi_{\sf B}}\psi\rangle=\Re\,p\,\, \langle \psi,\psi\rangle +c\langle {\bchi_{\sf B}}\psi,{\bchi_{\sf B}}\psi\rangle=0
\end{equation}
If $\Re\,p>0$ this immediately implies $\psi=0$. If $\Re\,p=0$  we get
 ${\bchi_{\sf B}}\psi=0$. But
${\bchi_{\sf B}}\psi=0$ implies $0=\mathcal{A}_\beta\psi=\mathcal{A}_0\psi$.  In
spherical coordinates the equation $\mathcal{A}_0\psi=0$ becomes a system
of ordinary differential equations
\begin{equation}
  \label{eq:spherical}
 \left( -\frac{d^2}{dr^2}-\frac{2}{r}\frac{d}{dr}+\frac{l(l+1)}{r^2}-br^{-1}-ip\right)\psi_{n,l,m}=0
\end{equation}
Since ${\bchi_{\sf B}}\psi=0$, the solution of (\ref{eq:spherical}) vanishes identically
on $[0,a]$; but then, by standard arguments the solution is identically
zero.  

If $\Im\,p \notin [-\epsilon, p_c ] $, with $p \in \mathcal{D}$, 
then $\mathcal{A}_\beta=\mathcal{A}_0$, and we are, by
construction, outside the spectrum of $\mathcal{A}_0$,
and thus
$\mathcal{A}_0 \psi =0$ implies $\psi=0$.  

{\bf (2)} {\em The range of $\mathcal{A}_{\beta}$ is dense.}  Indeed, the
opposite would imply \footnote{\cite{Kato}, p. 267}
$\mathrm{Ker}(\mathcal{A}^*_{\beta})\ne 0$, which leads to the same
contradiction as in Step 1 (note that $\mathcal{A}^*_{\beta}$ is simply
$\mathcal{A}_{\beta}$  with the signs of $\beta$ and $\Re\,p$ changed at the
same time).

{\bf (3)} {\em For any $p\in\mathcal{D}$ there is an $\epsilon>0$
  such that $\|\mathcal{A}_{\beta}\psi\|>\epsilon\|\psi\|$}.

 (a) If $\Re \,p> 0$ and $\|\psi\|=1$, then 
\begin{equation}
  \|\mathcal{A}_{\beta}\psi\|\geqslant
  \left|\langle\mathcal{A}_{\beta}\psi,\psi\rangle\right|=
  \left|\langle\mathcal{A}_0\psi,\psi\rangle-i c \langle{\bchi_{\sf B}}\psi,{\bchi_{\sf B}}\psi\rangle\right|\geqslant
  \left|\Re\, p\langle\psi,\psi\rangle\right| \ge \Re\, p
\end{equation}

 (b) Let now $\Re\,p=0$, and assume $\Im\,p$ is between two 
 eigenvalues of $-H_C$, the distance to the nearest being
 $\delta>0$. To get a
contradiction, assume that $\|\psi_j\|=1$ and
$\|\mathcal{A}_{\beta}\psi_j\|=\epsilon_j\to 0$.  Then
\begin{equation}
 \epsilon_j= \|\mathcal{A}_{\beta}\psi_j\|\geqslant
\left|\langle\mathcal{A}_{\beta}\psi_j,\psi_j\rangle\right|=
\left|\langle\mathcal{A}_0\psi_j,\psi_j\rangle-i c \langle{\bchi_{\sf B}}\psi_j,{\bchi_{\sf B}}\psi_j\rangle\right|\geqslant
c\left|\langle{\bchi_{\sf B}}\psi_j,{\bchi_{\sf B}}\psi_j\rangle\right|\to 0
\end{equation}
thus ${\bchi_{\sf B}}\psi_j\to 0$, and by the definition of
$\mathcal{A}_{\beta}$ and $\mathcal{A}_0$ we get
\begin{equation}\label{CD2}
\|\mathcal{A}_0\psi_j\|\to 0
\end{equation}
which is impossible, since our assumption and (\ref{CD2}) imply
noninvertibility of $H_C-ip$ while $ip$ is
outside the spectrum of $H_C$.

 (c) In the last case, $\Re\,p=0,\Im\,p\in\sigma_p(-H_C)$; 
then if we assume there is a
sequence $\psi_j$, $\|\psi_j\|=1$ such that $\|\mathcal{A}_{\beta}\psi_j\|\to 0$
as $j\to \infty$ we get
\begin{equation}\label{LB3}
  \|\mathcal{A}_{\beta}\psi_j\|\geqslant
  |\langle\mathcal{A}_{\beta}\psi_j,\psi_j\rangle|=
  \left|\langle\mathcal{A}_0\psi_j,\psi_j\rangle-i c 
\langle{\bchi_{\sf B}}\psi_j,{\bchi_{\sf B}}\psi_j\rangle\right|
\geqslant\left|c \langle{\bchi_{\sf B}}\psi_j,{\bchi_{\sf B}}
\psi_j\rangle\right|\to
  0
\end{equation}
Since $\| \mathcal{A}_0 \psi_j \| \le \| \mathcal{A}_\beta \psi_j \|
+ c \| {\bchi_{\sf B}} \psi_j \| $, (\ref{LB3}) implies 
$\|\mathcal{A}_0\psi_j\|\to 0$. On the other hand, with $P$ the
orthogonal projection on the finite dimensional eigenspace of $H_C$
corresponding to the eigenvalue $ip$, we have
$\mathcal{A}_0 P=0\Rightarrow \mathcal{A}_0=\mathcal{A}_0(I-P)$ and then
 since $\mathcal{A}_0 \psi_j \rightarrow 0$,
\begin{equation}
  \label{eq:inva}
  \|\mathcal{A}_0(I-P)\psi_j\|\to 0
\end{equation}
But by definition $\mathcal{A}_0$ is invertible on $(I-P)L^2(\RR^3)$ and
(\ref{eq:inva}) then implies $\|(I-P)\psi_j\|\to 0$, i.e.  $P\psi_j
- \psi_j\rightarrow 0$. Since $\|\psi_j \|=1$, $\|P\psi_j\|\to 1$.  Then
$P\psi_j$ is a bounded sequence in the finite dimensional space $P
L^2(\RR^3)$, hence we can extract a convergent subsequence, which we may
without loss of generality assume to be $P\psi_j$ itself,
$P\psi_j\to\psi,\,\|\psi\|=1$, and also $\psi_j \rightarrow P \psi_j
\rightarrow \psi$, thus $P\psi=\psi$.  Therefore, $\mathcal{A}_0 \psi =
\mathcal{A}_0 P \psi = 0 $.  Also, since multiplication by 
$ c {\bchi_{\sf B}}$ is a
bounded operator we have $ c {\bchi_{\sf B}} \psi_j \rightarrow  c 
{\bchi_{\sf B}} \psi =0$, since
$ c {\bchi_{\sf B}} \psi_j \rightarrow 0$.  
Therefore, $\| A_\beta \psi \| \le \|
\mathcal{A}_0 \psi \| + \| c {\bchi_{\sf B}} \psi \| = 0$ in contradiction to
the absence of eigenvalues.

\z {\bf(4)} {\em Definition of the inverse}. This is standard: we let
$\psi\in D(\mathcal{A}_{\beta})$,  $\mathcal{A}_{\beta}\psi=\phi$ and define
$\mathfrak{R}_{\beta}\phi=\psi$. This is well defined since
$\mathcal{A}_{\beta}\psi_1=\mathcal{A}_{\beta}\psi_2$ entails, by Step 1,
$\psi_1=\psi_2$.  By Step 2, $\mathfrak{R}_{\beta}$ is defined on a
dense set. By Step 3, for any $p$ there is an $\epsilon>0$ such
that  $\|\mathfrak{R}_{\beta}\|<\epsilon^{-1}$. Thus
$\mathfrak{R}_{\beta}$ extends by density to $L^2(\RR^3)$ and by construction
$\mathcal{A}_{\beta}\mathfrak{R}_{\beta}\phi=\phi$ whenever
$\mathfrak{R}_{\beta}\phi\in D(\mathcal{A}_{\beta})$. Conversely, if $\phi\in
D(\mathcal{A}_{\beta})$, and $\mathcal{A}_{\beta}\phi=u$ then $\mathfrak{R}_{\beta}u=\phi$ entailing $\mathfrak{R}_{\beta}\mathcal{A}_{\beta}\phi=\phi$ on the dense
set D$(\mathcal{A}_{\beta})$.

For the regularity of $\mathfrak{R}_\beta$ in $x$,  
we first note that if we define
$\mathcal{Q}=(I-\Delta )^{-1}$, we have the following identity:
\begin{equation}
\label{ResolId}
\mathfrak{R}_\beta = \mathcal{Q} 
\left [ 1 - \left ( \frac{b}{r} + i \beta \bchi_B
+ i p +1 \right ) \mathcal{Q} \right ]^{-1}   
\end{equation}
It is clear that if 
$\phi \in L^2(\mathbb{R}^3)$, $\mathcal{Q} \phi \in H^2 (\mathbb{R}^3)$ and
so $\left ( b/r - i \beta \bchi_B + i p + 1 \right ) \mathcal{Q} 
\phi \in 
L^2$. Therefore, from (\ref{ResolId}), 
$\mathfrak{R}_\beta : L^2 (\mathbb{R}^3) \rightarrow H^2 (\mathbb{R}^3) $.

 \subsection{Proof of Proposition~\ref{P010}}\label{010}
The shift operator $S$, defined  by  $(S Y)_j=y_{j+1}$,  is quite
straightforwardly shown to be bounded in $\mathcal{H}$: the proof of
Lemma 27 in \cite{JSPLieb} goes through without changes. By the second
resolvent identity we have
$$\mathfrak R_{\beta}=\Big(1-i\beta\mathfrak R_0{\bchi_{\sf B}}\Big)^{-1}\mathfrak
R_0$$
Since $-\Delta-br^{-1}$ is self-adjoint, we have by the spectral
theorem, for some $C>0$ independent of $p$,
\begin{equation}
  \label{eq:eqcontR}
  \|(-\Delta-br^{-1}-ip)^{-1}\|_{L^2(\RR^3)}\leqslant  C(\Re\, p)^{-1}
\end{equation}
and thus $\|\mathfrak R_{\beta}\|_{L^2({\sf B})}\leqslant  C_1 (1+|\Re\,p|)^{-1}$. 
Since $\mathfrak{R}_{\beta}$ is diagonal (in $n$) and $\mathfrak{S}$
is bounded (cf. Lemma \ref{lemY0C}), we have  $\|\mathfrak C\|_{\mathcal{H}}\leqslant   C_2
  (1+|\Re\,p_1|)^{-1}$. 
Thus $\| \mathfrak{C} \|_{\mathcal{H}}$ is small
 for
 large $\Re\, p_1$ and therefore $(I-\mathfrak{C})Y = Y^{(0)}$ has
 unique solution $Y \in \mathcal{H}$ and the proof follows.

\subsection{
Proof of Proposition~\ref{P12}}
\begin{proof} \label{P12P1} 
 The estimate $\|{\bchi_{\sf B}} \mathfrak{R}_0 {\bchi_{\sf B}} \|=O(p^{-1/2})$  is shown right after the statement of Proposition~\ref{P12}.  

   We now consider the
  analyticity of $\mathfrak{R}_\beta$ in an open set on the imaginary $p$
  axis for $\Im ~p < -\epsilon$.  There, $\beta=0$ and
  $\bchi_B \mathfrak{R}_\beta \bchi_B =
\bchi_B \mathfrak{R}_0 \bchi_B $ is manifestly analytic from its
representation as an integral operator, 
  whose  kernel $G$  is given below (see \cite{Hostler} 
and Appendix \S \ref{CoulombG}
for details).

With $k=\sqrt{ip}$ (using the principal branch of the square root),
and $\nu = b/(2k)$, 
\begin{multline}
\label{eq:exprr}
G(x,x';k)\\=  \frac{i k(\eta-\xi) I (-ik \xi) J (-i k \eta) -k^2\xi\eta
[I (-i k \xi) \dot{J} (-i k \eta) -
 J (-i k \eta) \dot{I} (-i k \xi) ]}{\Gamma(1-i\nu)\Gamma(1+i\nu)} 
\frac{e^{\frac{ik}{2}(\xi+\eta)}}{4\pi |x-x'|},  
\end{multline}
where
\begin{equation}
\label{defxieta}
\xi = |x| + |x'| + |x-x'|  ~~~~~~~\\, ~~~~~ 
\eta = |x| + |x'| - |x-x'| ~~, 
\end{equation}
\begin{multline}\label{defIJ}
  I (z_1)=\int_0^{i \infty} e^{-z_1 t}
t^{-i\nu}(1+t)^{i\nu}dt,\;\ \dot{I} (z_1 ) =-\int_0^{i \infty} e^{-z_1 t}
t^{1-i\nu}(1+t)^{i\nu}dt\\
J (z_2)=\int_0^{1}
 e^{z_2 t}t^{-i\nu}(1-t)^{i\nu}dt;\ \dot{J} (z_2)=\int_0^{1}
 e^{z_2t}t^{1-i\nu}(1-t)^{i\nu}dt
\end{multline}
\end{proof}
\noindent  Further properties of function  $G$ 
are discussed
in
\S \ref{CoulombG}. 
\begin{Note}\label{NoteG}{\rm 
 Note that  (\ref{eq:exprr}) 
still holds for $p\in i\RR^+$, with $k=\sqrt{ip}$,  with the choice
$\arg k = \pi/2$  for $p \in i\RR^+$,  and with the upper limits
$i \infty$ in (\ref{defIJ}) replaced by $+\infty$.}
\end{Note}

\subsection{Proof of  Proposition~\ref{P131}}\label{pp131}

The function $f=\mathfrak{R}_{{\beta},l,m} \bchi_{\sf B} g$ 
is the solution
of the equation
\begin{eqnarray}
  \label{eq:eq323}
  \Big(-\frac{d^2}{dr^2}-\frac{2}{r}\frac{d}{dr}+\frac{l(l+1)}{r^2}-\frac{b}{r}+\lambda^2-ic\Big)f=g;\
  \ r\leqslant  a\nonumber\\
 \Big(-\frac{d^2}{dr^2}-\frac{2}{r}\frac{d}{dr}+\frac{l(l+1)}{r^2}-\frac{b}{r}+\lambda^2\Big)f=0;\
  \ r> a\nonumber\\
\end{eqnarray}
such that $f$ decays at infinity, is regular at the origin and $C^1$ at
$r=a$. 
We note $\lambda =\sqrt{-ip}$ is in the closure of
the fourth
quadrant for $\Re\, p \geqslant 0$.
We let  $\alpha=\sqrt{\lambda^2-i{{c}}}$,
$\kappa_1=b/(2\alpha)$, $\kappa=b/(2\lambda),\mu=2l+1$ and define
(in terms of the Whittaker functions $\mathfrak{M}$ and $\mathfrak{W}$  \footnote{See \cite{Buchholz}, 
pp. 60 eq. (1) and pp. 63 eq. (5).}) \index{$\mathfrak{M} $}\index{$\mathfrak{W} $} 
\begin{multline}\label{defm1}
  m_1(s):=s^{-1}\mathfrak{M}_{\kappa_1,\mu/2}(2\alpha s);\ \ w_1(s):=
s^{-1}\mathfrak{W}_{\kappa_1,\mu/2}(2\alpha s); \ w_2(s):=
s^{-1}\mathfrak{W}_{\kappa,\mu/2}(2\lambda s)
\end{multline}
For $r>a$ we have $f=\mathrm{B}w_2(r)$ since $r^{-1} 
\mathfrak{M}_{\kappa,
\mu/2} (2 \lambda r)$ grows with $r$ as $r \rightarrow \infty$.
For $r\leqslant a$ we must have 
\begin{equation}
  \label{eq:eqff}
  f=A
m_1+f_0
\end{equation}
\z where, using standard results about the Wronskian of $\mathfrak{M}$ and
$\mathfrak{W}$, see    
\cite{Buchholz}, pp. 25 and \cite{Abramowitz}, pp 505, 508, we have
\begin{multline}
  \label{eq:inner}
 \frac{2 \alpha \Gamma (1+\mu)}{ 
 \Gamma\left(\frac{1}{2}+\frac{1}{2}\mu-\kappa_1\right)} f_0\\=
 w_1(r)\int_0^r {\bchi}_{[0, a]}(s)s^2 m_1(s)g(s)ds+
m_1(r)\int_r^{a} s^2w_1(s) {\bchi}_{[0, a]}(s)g(s)ds
\end{multline}
The integral representations of the functions $\mathfrak{M}$ and
$\mathfrak{W}$ $^{\thefootnote}$  entail immediately
that the functions $f_0$, $f$ and $\mathfrak{M}_{\kappa_1,\mu/2}(2\alpha r)$
depend analytically on $\lambda$ for small $\lambda$. 
Continuity of $f$ and $f'$ at $a \geqslant 1$ imply that $A$ defined in
(\ref{eq:eqff}), is given by
\begin{equation}
  \label{eq:valA}
  A=\frac{f_0(a) w'_2(a)-f'_0(a) w_2(a)}{m'_1(a) w_2(a)-
m_1(a) w'_2(a)}
\end{equation}
In  \S \ref{Adep} it is shown that
that $A$ is analytic in 
$ \left ( \lambda , 
\exp \left [ i \pi b/(2 \lambda) \right ] \right ) $  in 
a domain corresponding to   
$\lambda$ small in the closure of the fourth quadrant, 
if $a$ and $c$ are chosen large enough.
It follows that resolvent $\mathfrak{R}_{\beta,m,n}$ is analytic in $X$
for 
$X = \left ( \sqrt{p},    
\exp \left [ i \pi b/\left ( 2 \sqrt{-ip} \right ) \right ] \right ) 
\in \overline{\mathbb{D}_\epsilon^+}  \times \overline{\mathbb{D}} $ 
for  small $\epsilon$. 

\subsection{
Proof of Proposition~\ref{Ccomp}}\label{020}

   By adding and subtracting $1$
  from $\mathcal{A}_\beta$ and  using the second
  resolvent formula, whenever everything is well
  defined, we have
\begin{multline}
  \label{eq:2res2}
 {\bchi_{\sf B}} \mathcal{A}^{-1}_{\beta}{\bchi_{\sf B}}=
:{\bchi_{\sf B}} \mathfrak{R}_{\beta}{\bchi_{\sf B}}= 
{\bchi_{\sf B}}(-\Delta+ 1)^{-1}{\bchi_{\sf B}}
\\ -{\bchi_{\sf B}}\mathfrak{R}_{\beta}(-br^{-1}-i\beta 
{\bchi_{\sf B}}-1 - i p)(-\Delta+ 1)^{-1}{\bchi_{\sf B}}
\end{multline}
The Green's function for $-\Delta+ 1$ is 
\begin{equation}
  \label{eq:gfdelta}
  G(x,y)=\frac{1}{4\pi|x-y|}e^{-|x-y|}
\end{equation}
Now if $\|\phi_j\|_{L^2({\sf B})}\leqslant  1$ then the functions
$f_j=(-\Delta+ 1)^{-1}{\bchi_{\sf B}} \phi_j$ are seen by straightforward
calculation to be equicontinuous on the one point 
compactification of $\RR^3$. A subsequence, without loss of
generality assumed to be the $f_j$'s themselves, 
converges in $L^2(\RR^3)$ as
well (to a function with exponential decay, 
since there is a  $\delta_1>0$ small enough and independent of $j$ so that
$e^{\delta_1|x|}(-\Delta+ 1 )^{-1}{\bchi_{\sf B}} \phi_j)$ 
is also equicontinuous on
the compactification of $\RR^3$).  
In particular, ${\bchi_{\sf B}}(-\Delta+ 1)^{-1}{\bchi_{\sf B}}$  
is compact. 

Now $f_j$ converge in the sup norm with weight $e^{\delta_1 |x|}$, and thus
 $(-br^{-1}-i\beta {\bchi_{\sf B}}- 1- ip) f_j$ converge in 
$L^2(\RR^3)$.
Since $\mathfrak{R_{\beta}}$ is  bounded, compactness
of ${\bchi_{\sf B}} \mathfrak{R}_\beta {\bchi_{\sf B}}$ follows.  

By Proposition~\ref{P12}, and the previous argument, $\mathfrak{C}$ is a norm
limit of compact operators (the truncations of $\mathfrak{C}$ to the subspaces
of $\mathcal{H}$ with vanishing components for $|n|>N$).  Therefore, 
$\mathfrak{C}_{l,m} = P_{l,m} \mathfrak{C}$ is also compact.

\subsection{Proof of Proposition \ref{P20} and final estimates
for Theorem \ref{T3}}\label{Fredh}

If (\ref{eq:chomo}) has no nontrivial solution for any $p_1 \in 
\overline{\mathbb{H}}$, then compactness of 
$\mathfrak{C}_{l,m}$ implies that $\left ( I - \mathfrak{C}_{l,m} 
\right )^{-1}$ exists. Lemma \ref{lemCanal}, Corollary
\ref{CorP131} and Proposition \ref{P12}
give the analytic and continuity properties of
$\mathfrak{C}_{l,m}$ and $\mathfrak{T}_{l,m} Y_0$.  
Analyticity of $\left (I-\mathfrak{C}_{l,m} \right )^{-1}$ in 
$X_1$ for $X_1 \in  
\overline{\mathbb{D}^+_\epsilon}\times \overline{\mathbb{D}}$, 
follows in a standard way from analyticity of $\mathfrak{C}_{l,m}$ and
the second resolvent formula, 
\begin{equation}
  \label{eq:secondr}
 A^{-1}-B^{-1}=B^{-1}(B-A)A^{-1}
\end{equation}
(see \S\ref{CX}).
The same resolvent identity 
can be applied to show analyticity of $(I-\mathfrak{C}_{l,m})^{-1} $
with respect to $p_1$ in a neighborhood of 
$$\overline{\mathbb{H}} 
  \setminus \left \{ (\ell_{p_c} + i \omega \ZZ) \cup  (\ell_{-\epsilon} + i
    \omega \ZZ) \cup ( \mathcal{I}_\epsilon + i \omega {\ZZ}) \right \}
$$ Hence, the  solution $Y = \left ( I - \mathfrak{C}_{l,m} \right
)^{-1} \mathfrak{T}_{l,m} Y_0 $ is analytic for $p_1 \in
\overline{\mathbb{H}} \setminus \{ i \omega \ZZ \}$, since $\epsilon$,
$\beta$ and $p_c$ are artificially introduced parameters the
value of which 
cannot affect $Y$, since $Y_0$ is independent of these choices (see
Remark \ref{remanal}.)

The function $\hat{y}(p;x) = y_n (p_1, x)$, with $p=in \omega + p_1$, 
is analytic in $p $ for
$p \in i\RR \setminus{i \omega \ZZ}$ and by
analyticity 
of $Y$ in $X_1$, boundedness at $p=i \omega \ZZ$ follows. In particular,
as $p \rightarrow i n \omega$ from the right half-plane,
${\hat y} (p, x)$ is analytic in the extended variable 
\begin{equation}
\label{extendedpar}
\left ( (p-i n \omega)^{1/2}, 
\exp \left [ \frac{i\pi b}{2\sqrt{-i (p-in\omega)} } \right ] \right )
\end{equation}

The regularity properties of $Y $ in $p$ and the decay properties in $|n|$ of its
components $y_n$ for large $|n|$, simply stemming from  $Y \in \mathcal{H}$,
imply that
$y(t,x)$ can be
expressed as an inverse Laplace transform of ${\hat y} (p, x)$
on $i \mathbb{R}$.
We now show that  
$$P (t, {\sf B} ) = 
\| \psi_0 (x) e^{-t} + y (x, t) \|_{L^2 ({\sf B})}^2 \le 
2 e^{-2t}\|\psi_0\|^2_{L^2 ({\sf B})}+ 
2 \|y (x, t) \|_{L^2 ({\sf B})}^2 \to 0 \ 
\text{as}\ t\to\infty$$

We note that
\begin{multline}
\int_{\sf B} dx |y (t, x) |^2 
= \int_{\sf B} dx \int_{-\infty}^\infty \int_{-\infty}^{\infty} 
e^{i t (s-s') } {\hat y} (is, x) \overline{{\hat y} (is', x)} ds ds'\\
= \int_{-\infty}^\infty e^{i {\tilde s} t} 
\left \{ \int_{-\infty}^\infty \left [ \int_{\sf B} 
{\hat y} (i{\tilde s}+i s', x) 
\overline{{\hat y} (i s' , x)} dx \right ] ds' \right \} d{\tilde s}  
\end{multline}
So, in order to show ionization, it suffices from Riemann-Lebesgue Lemma
to show that 
$$  
 \int_{-\infty}^\infty \left [ \int_{\sf B} 
\overline{{\hat y} (is', x)} 
{\hat y} (i s' + i {\tilde s}, x) dx \right ] ds' 
$$ 
is in $L^1 (d{\tilde s})$. This follows from Cauchy-Schwarz, since
\begin{multline}
  \int_{-\infty}^{\infty} 
\left \{ \int_{-\infty}^\infty \left [ \int_{\sf B} |{\hat y} (is', x)| 
|{\hat y} (i s' + i {\tilde s}, x)| dx \right ] ds' \right \} d{\tilde s}  
\le \left ( 
\int_{\mathbb{R}}  \| {\hat y} (i s', \cdot ) \|_{L^2 ({\sf B}) } ds' 
\right )^2 
\end{multline}

\z However, 
\begin{multline}
\int_{\mathbb{R}}  \| {\hat y} (i s', \cdot ) \|_{L^2 ({\sf B})}  ds' 
= \int_{0}^\omega \left [ \sum_{ n \in \mathbb{Z} } 
\| y_n (i q , \cdot ) \|_{L^2 ({\sf B})} \right ] dq \\ 
\le C \int_0^{\omega} \sum_{n \in \mathbb{Z}} 
(1+|n|)^{4/3} \| y_n (iq, \cdot ) \|^2_{L^2 ({\sf B})} dq
\le C \int_0^\omega \| Y (iq, \cdot) \|^2_{\mathcal{H}} dq < \infty
\end{multline} 
since $Y$ is bounded  in $p_1 = i q$ for $q \in [0, \omega]$. 

  Since ${\hat y} (p_1+in\omega, x)$, $n\in\ZZ$, is
 analytic in the variable (\ref{extendedpar}), standard
  stationary phase analysis (see Appendix \S \ref{Newsp} ) 
  shows that  $y (t, x)=O(t^{-5/6})$, and hence
  $P (t, {\sf B} )=O(t^{-5/3})$ as $t \rightarrow \infty$.

\subsection{Proof of Theorem~\ref{cdi(ii)} }\label{iproof}
Since (\ref{eqvv}) (restricted to ${\sf B}$) follows from the
homogeneous system $w=\mathfrak{C} w$ (see also Proposition ~\ref{P4}
for the necessary regularity), we look for a nontrivial solution of
(\ref{eqvv}) in $\mathcal{H}$.  We multiply (\ref{eqvv}) by
$\overline{w}_n$, integrate over the ball ${\sf B}_{\tilde a}$ (of
radius ${\tilde a} \in (1, a] $), sum over $n$ (this is legitimate
since $w\in\mathcal{H}$) and take the imaginary part of the resulting
expression.  Noting that
\begin{multline}
  \label{conjsym}
  \overline{\sum_{j,n\in\ZZ}\Omega_j(x)w_{n-j}\overline{w_n}}=\sum_{j,n\in\ZZ}
\Omega_{-j} \overline{w}_{n-j}w_{n}=\sum_{j,n\in\ZZ}
\Omega_{j} \overline{w}_{n+j}w_{n}\\=\sum_{j,m\in\ZZ}\Omega_j(x)\overline{w_m}w_{m-j}\end{multline}
so the sum (\ref{conjsym}) is real, we get from (\ref{eqvv}) 
\begin{multline}
  \label{eq:nonreal}
  0=\Im\left(+ip_1\sum_{n\in\ZZ}\int_{{\sf B}_{\tilde a}}|w_n(x)|^2dx +
\int_{{\sf B}_{\tilde a}}\sum_{n\in\ZZ} dx
    \overline{w}_n\Delta w_n\right)\\=
  +\Re\,p_1\sum_{n\in\ZZ}\int_{{\sf B}_{\tilde a}} |w_n(x)|^2dx 
+\frac{1}{2i} \int_{\partial
    {\sf B}_{\tilde a}}\left(\sum_{n\in\ZZ}\overline{w}_n\nabla w_n-w_n\nabla
    \overline{w}_n\right)\cdot \mathbf{n}\,dS\end{multline}

 It is
convenient to decompose $w_n$ using spherical harmonics. We write
\begin{equation}
  \label{sph}
  w_n=\sum_{l\geqslant 0, |m|\leqslant  l}R_{n,l,m}(r)\mathcal{Y}_l^m(\theta,\phi).
\end{equation}
The last integral in (\ref{eq:nonreal}), including the prefactor, then
equals 
\begin{multline}
  \label{sph2}
  -\frac{i}{2}\, a^2 \,\sum_{n\in\ZZ}\sum_{m,l}\Big[\overline{R}_{n,m,l}R'_{n,m,l}-
  \overline{R'}_{n,m,l}R_{n,m,l}\Big]\\=-\frac{i}{2}\,  a^2 \, \sum_{n\in\ZZ}\sum_{m,l}{\mathcal{W}}[\overline{R}_{n,m,l},R_{n,m,l}]
\end{multline}
where ${\mathcal{W}}[f,g]$ is the Wronskian of $f$
and $g$.  On the other hand,  we have outside of ${\sf B}_{\tilde a}$ 
\begin{equation}
  \label{eq:outside}
  \Delta w_n+br^{-1}w_n+(i p_1-n\omega)w_n=0
\end{equation}
and then by (\ref{sph}), the $R_{n,l,m}$ satisfy for $r>a$ the equation
\begin{equation}
  \label{Rnlm}
  R''+\frac{2}{r}R'+br^{-1}R-\frac{l(l+1)}{r^2}R=(-ip_1+n\omega)R
\end{equation}
where we have suppressed the subscripts. Let $g_{n,l,m}=rR_{n,l,m}$.
Then for the $g_{n,l,m}$ we get
\begin{equation}
  \label{gnlm}
  g''-\left[\frac{l(l+1)}{r^2}-ip_1+n\omega-br^{-1}\right]g=0
\end{equation}
Thus
\begin{equation}
  \label{r-g}
  \overline{R}R'=\frac{\overline{g}g'}{r^2}-\frac{|g|^2}{r^3}
\end{equation}
and
\begin{equation}
  \label{r-g1}
  r^2{\mathcal{W}}[\overline{R},R]={\mathcal{W}}[\overline{g},g]=:{\mathcal{W}}_n.
\end{equation}
Multiplying (\ref{gnlm}) by $\overline{g}$, and the conjugate of
(\ref{gnlm}) by $g$ and subtracting, we get for $r>a$,
\begin{equation}
  \label{difg}
  {\mathcal{W}}_n'=-i(p_1+\overline{p_1})|g|^2=-2i|g|^2 \Re\, p_1
\end{equation}

\begin{Remark}\label{condifty}
  Direct estimates using the Green's function representation (\ref{eq:exprr})
  imply that
\begin{equation}
  \label{cinf}
  w_n(x)=\frac{e^{-\kappa_n r}}{r^{1+\frac{b}{2\kappa_n}}}\Big(c_n(\theta,\phi)+O(r^{-1})\Big)\ \text{as}\ r\to\infty
\end{equation}
with $c_n(\theta,\phi)$ independent of $r$ and
with
\begin{equation}
  \label{defkappa}
  \kappa_n= \sqrt{-i p_1 + n \omega} \ \text{(when
$\Re\,p_1 >0$,  $\kappa_n$ is in the fourth quadrant when $n<0$})
\end{equation}
\end{Remark}
\z (i) We first take $\Re\, p_1>0$, to illustrate the argument.  
Using (\ref{cinf}) we get
\begin{equation}
  \label{ginfty}
  g\sim C e^{-\kappa_n r}r^{-\frac{b}{2\kappa_n}}(1+o(1))\ \ \text{as}\ \ r\to\infty
\end{equation}

There is a one-parameter family of solutions of (\ref{gnlm})
satisfying (\ref{ginfty}) and the asymptotic expansion can be
differentiated \cite{Wasow}. We assume, to get a contradiction, that
there exist $n < 0 $ for which $g=g_n\ne 0$.  For these $n$ we have, using
(\ref{ginfty}), differentiability of this asymptotic expansion and
the definition of $\kappa_n$ that
\begin{equation}
  \label{t2}
  \frac{1}{2i}\lim_{r\to\infty}|g_n|^{-2} {\mathcal{W}}_n=-\Im\,\kappa_n\, >0.
\end{equation}
It follows from (\ref{difg}) and (\ref{t2}) that $\mathcal{W}_n/(2i)$ 
is strictly
positive for all $r>a$ (by monotonicity and positivity at infinity) and all $n$ for which $g_n\ne 0$.  This implies that
the last term in (\ref{eq:nonreal}) is a sum of nonnegative terms which shows
that (\ref{eq:nonreal}) cannot be satisfied nontrivially.

(ii): $\Re\,p_1=0$. For $n<0$, we use Remark \ref{condifty} (and
differentiability of the asymptotic expansion as in Case (i)) to calculate
${\mathcal{W}}_n$ in the limit $r\to \infty$: 
${\mathcal{W}}_n=2i|c_n|^2|\kappa_n|
(1+o(1))$.  Since for $\Re\,p_1=0$, ${\mathcal{W}}_n$ is constant, cf.
(\ref{difg}), it follows that ${\mathcal{W}}_n=2i|c_n|^2|\kappa_n| |g_n|^2$ 
exactly. Thus, (\ref{eq:nonreal}) cannot be non-trivially satisfied,
implying that 
\begin{equation}
  \label{IonizCond}
  w_n(x)=0\ \ \text{for all}\ \ n<0 \ \text{and}\ |x| = r = {\tilde a}
\in (1, a]
\end{equation}
For ${\tilde a} > r>1$ (where $V(t,x)=0$) we have $\mathfrak{O} w_n=0$, where
$\mathfrak{O}$ is the elliptic operator $-\Delta-b/r-ip_1+n\omega$. The proof
that $w_n(x)=0$ for $r>1$ then follows immediately from (\ref{IonizCond}),
by standard unique continuation results \cite{Hormander,Treves} (in fact, $\mathfrak{O}$ is analytic hypo-elliptic).  See also
Note~\ref{N07}.

\subsection{Connection with the Floquet operator}\label{Fl3}
 It is easy to check that the discrete time-Fourier transform of the
eigenvalue equation for the Floquet operator, Eq. (\ref{eq:floquet1}),
$Kv=\phi v$, 
 with $p_1 = i \phi$, 
coincides with (\ref{eqvv}),
the differential version of
the homogeneous equation associated to (\ref{eq:intform}). Now,
(\ref{IonizCond}) 
shows that a solution of (\ref{eqvv}) is an eigenvector of $K$.

In the opposite direction the existence of a Floquet eigenfunction
entails failure of ionization since it implies the existence of a solution of
(\ref{sch-eq}) for which the absolute value is time-periodic. 

\subsection{Differential equation for $w$
}\label{Lts}

We seek to show that the
only solution to the homogeneous system 
\begin{equation}
\label{eqhom1}
w = \mathfrak{C}_{l,m} w
\end{equation}
in the space $\mathcal{H}$ is $w = 0$. 
 Since $w$ is piecewise ${C}^2$ (see Note \ref{N07}),
(\ref{eqhom1})
implies that the
components of $w = \left \{ \mathcal{Y}_{l,m} r^{-1} g_n (r)\right 
\}_{n\in\ZZ}$
satisfy
the differential-difference system (see Note \ref{N07}):
\begin{equation}
\label{eqhom2}
\frac{d^2}{dr^2} g_n - \left (-br^{-1} + n \omega -ip_1
+\frac{l (l+1)}{r^2}  \right ) g_n
= i \Omega \left (g_{n+1} - g_{n-1} \right )
\end{equation}
First, we notice that for $n < 0$,  Theorem \ref{cdi} implies that $g_n
(r) = 0$ for $r \geqslant 1$. Thus $g_n (1) =0$, $g_n^\prime (1) = 0$ for all $n <
0$. 
\subsection{Proof of Proposition~\ref{startpoint}}\label{Lemhom1} The gist 
of the proof is that contractive mapping arguments show that if the statement
was false then the 
solution would vanish. 
\begin{Lemma}
  If $Y \ne 0$, then there exists some $n_0 \geqslant 0$ so that
  either $g_{n_0} (1) \ne 0$ or $g_{n_0}^\prime (1) \ne 0$.  (As
  before, in the sequel, we shall define $n_0$ to be the smallest such
  integer.)
\end{Lemma}

\begin{proof}
To get a contradiction, assume the statement is false.
 Since the functions $w_n$ are in the domain of $\Delta$ (see Note
\ref{N07}), then, 
in particular,  for any
$n$, $g_n$ is continuous in $r$. Thus, the set $Z_n:=\{r:g_n(r)=0\}$ 
is closed and so
is the (possibly empty) left connected component of $1$ in $Z_n$, call it
$K_n$.  Let
$$K=\bigcap_{n\in\ZZ}K_n$$
Assume to get a contradiction that $K$ is nonempty:
let then $K=[a,1]$.  If $a=0$, then $Y\equiv 0$ since $g_n (1) =0 $,
$g_n^\prime (1) =0$ imply $g_n(r)=0$ for $r>1$.  Then $Y \ne 0$ implies 
$ a > 0$.   We first take $0 < a < 1$. We write the
differential equation for $g_{n} (r)$ in integral form and use the conditions
$g_n (a) =0 = g_n^\prime (a)$, since $g_n$ vanishes
on $[a,1]$:
\begin{multline}\label{intf}
  e^{\sqrt{n\omega} r } g_n (r) = 
\int_r^a \left ( \frac{[1 - e^{-2 \sqrt{n\omega} (s-r)}
      ]}{2 \sqrt{n\omega}} \right ) e^{\sqrt{n\omega} s} \\ 
\left \{ \left [\frac{l
        (l+1)}{s^2} -ip_1 + \tilde{V}(s) \right ] g_{n} (s) -i \Omega (s) \left (
      g_{n-1} (s) - g_{n+1} (s) \right ) \right \} ds
\end{multline}
Consider the Banach space of sequences
$$ \left \{ g_n (r) \right
\}_{n=-\infty}^\infty $$
in the norm
$$\sup_{n\in\ZZ,r\in[a-\epsilon,a]}\left|e^{\sqrt{n\omega} r} g_n (r)\right|$$
It is easy to see that the rhs of (\ref{intf}) is a contractive mapping if
$\epsilon$ is small enough and then $g_n (r)=0$ for $r\in[a-\epsilon,a]$
contradicting the definition of $a$. The same
is true if $a=1$, since 
$g_n (1) =0$ and $g_n^\prime (1)=0$ would imply, with the same  proof as
before, that  $g_n =0$ for $r \in [1-\epsilon, 1]$, for some $\epsilon>0$,
contradicting the definition of $a$.
\end{proof}
\subsection{Proof of Theorem~\ref{Rnasympt}:}\label{prf91}
For a heuristic discussion see
\S\ref{Hr}.

The proof is by rigorous WKB. 
The fact that there are two competing potentially large
variables, $k$ and $1/r$ makes it necessary to rigorously match two regimes.

\z First, note that
(\ref{DefRk}) implies
\begin{equation}
\label{n9}
g_{n_0-k} (r) = i^k m_k (r) h_k (r)
\end{equation}
We need a few more preliminary results.

\begin{Lemma}
\label{lgkk0}
For any $\epsilon_1 > 0$, there exists $C_3>0$
independent of $k$ and $\epsilon_1 $ so that for $k \geqslant k_0 = C_3
\epsilon_1^{-1}$, and for $r \in [\epsilon_1, 1]$,
\begin{equation}
\sup_{\epsilon_1\le r\le 1}| h_k^\prime |
\leqslant  C_4 k_0 \left ( \frac{k_0}{k} \right )^{1/2}
\end{equation}
where $C_4$ is independent of $\epsilon_1$ and $k$.
\end{Lemma}
The proof of  Lemma \ref{lgkk0} is given in \S \ref{plgkk0}.

\begin{Definition}\label{definehk}
   For  fixed $\epsilon $, we define \index{$ L_\epsilon$}
  $\displaystyle L_{\epsilon}=\alpha C_3\left( 2 C_4
  C_3/\epsilon \right)^{2}$, with $C_3$ and $C_4$ 
defined in Lemma \ref{lgkk0}, 
and $\zeta=\alpha k r$,
where $\alpha$ is given in (\ref{eq:defquants}). We will take
$\epsilon$  small enough so that $L_\epsilon \ge C_3 \alpha$.

Finally, in what follows, $c_*$ \index{$c_* $} is a positive ``generic''
constant, the value of which is immaterial.
\end{Definition}
\begin{Lemma}
\label{gkl}
 For $\epsilon > 0$ small enough and $k \alpha r = \zeta \in [L_{\epsilon}, k\alpha]$,
we have 
\begin{equation}
  \label{eq:eqt1}
  |h_k (r) -1 | \leqslant  \epsilon
\end{equation}
\end{Lemma}
The proof of  Lemma \ref{gkl} is given in \S \ref{pgkl}.

\begin{Definition}
Let ${\tilde h}_k (\zeta) = h_k (\zeta/(\alpha k))$.
\end{Definition}

\begin{Lemma}
\label{lemarzela}
For any small $\epsilon > 0$,
there exists a subsequence $S= \{{\tilde h}_{k_j}\}_{j\in\NN}$ that converges
to a continuous function ${\tilde h} $ for $ \zeta\in [0, L_{\epsilon}] $.
For the limiting function 
${\tilde h} (\zeta)$, we have $|{\tilde h} (\zeta) - 1 | 
\le 4 \epsilon $ for $\zeta \in [0, L_\epsilon ]$.
\end{Lemma}
The proof of this proposition is given in \S\ref{Sarz}.
\begin{Proposition}\label{P33.1}
 For any $r \in [0,1]$, $\lim_{j \rightarrow \infty} h_{k,j} (r) = 1$.
\end{Proposition}
\begin{proof}
  From Lemma \ref{lemarzela} and Lemma  \ref{gkl} it follows that  for any $r\in[0,1]$ and any $\epsilon>0$ we have $\lim_{j\to\infty}|h_{k_j}(r)-1|\le 4\epsilon$. 
\end{proof}
\z The proof of Theorem~\ref{Rnasympt} 
now follows from the definition of $h_k$  in (\ref{S7.3}),
Remark~\ref{R46}, Note \ref{N07}   
and Proposition \ref{P33.1}.

\subsection{Further results on $g_{n_0-k}$ and $h_k$
}\label{Addtl}

\begin{Lemma}
\label{l0.1}
For any $j, k\in\NN\cup\{0\}$ we
have, at $r=1$, {\it i.e.} at $\mathfrak{s}=0$,
$$\frac{\partial^{j+\tau} g_{n_0-k} }{\partial \mathfrak{s}^{j+\tau}}  |_{\mathfrak{s} =0} = 
\delta_{j, 2k} ~i^{k}
\ \ \ \ \ \ {\rm for}~~0\leqslant  j\leqslant  2k 
$$
\end{Lemma}

\begin{proof}
   In case {\bf (i)} (corresponding to $\tau=0$),
 note that (\ref{eqhom2}) may be
  rewritten, cf. (\ref{defmk}), as  \index{$Q_k $}
\begin{equation}
\label{Rxi}
 (g_{n_0 -k})_{\mathfrak{s}\mathfrak{s}} - 
\frac{\Omega^\prime}{2 \Omega^{3/2}} (g_{n_0 -k})_{\mathfrak{s}}
+\frac{Q_{k}}{\Omega} g_{n_0-k}= i \left ( g_{n_0-k+1} - g_{n_0-k-1} 
\right ), 
\end{equation}
 where
\begin{equation}
\label{eqQk}
Q_k = \frac{b}{r} + (k - n_0 )\omega + i p_1 - 
\frac{l (l+1)}{r^2} 
\end{equation} 
Since $g_{n_0-k} (1) =0 = g_{n_0-k}^\prime (1) $ for all $ k \geqslant 1$, while
$g_{n_0} (1) = 1$, the statement follows from (\ref{Rxi}) for any $0 \leqslant 
j \leqslant  2$, if $2 k \geqslant j$.  Assuming the statement holds for some $j \geqslant
2$ for $2 k \geqslant j$, we prove it for  
$(j+1)$ for $2 k \geqslant (j+1)$.

Taking $(j-1)$ derivatives in  $\mathfrak{s}$ 
of (\ref{Rxi}) at $\mathfrak{s} =0$, we obtain
$$\frac{\partial^{j+1} g_{n_0-k}}{\partial \mathfrak{s}^{j+1}}=i
\frac{\partial^{j-1}}{\partial \mathfrak{s}^{j-1}} g_{n_0-(k-1)} - i
\frac{\partial^{j-1}}{\partial \mathfrak{s}^{j-1}} g_{n_0-(k+1)}+L$$
where $L$ is a
linear combination of derivatives of $g_{n_0-k}$ up to order $j$, which are
all zero since $2k \geqslant (j+1) > j$.  The first two terms on the rhs give a
contribution of $i i^{k} \delta_{(j-1), 2(k-1)} + 0$ since $2k \geqslant (j+1)$
implies $2 (k-1) \geqslant (j-1)$ and $ 2(k+1) > (j-1)$ completing the inductive
step.

\bigskip

 In case {\bf (ii)} (corresponding to $\tau=1$):
since $g_{n_0} (1) = 0$ and $g_{n_0-k}
(1) = 0 = g_{n_0-k}^\prime (1) $ for all $ k \geqslant 1$, 
it follows from
(\ref{Rxi}) that $g_{n_0-k}^{\prime \prime} = 0$ for all $k \geqslant 1$ 
implying
the conclusion for $j=0$ and $j=1$.  By taking an additional derivative of
(\ref{Rxi}) with respect to $\mathfrak{s}$ and evaluating at $\mathfrak{s}=0$, we obtain
  $$\frac{\partial^3 g_{n_0-k}}{ \partial \mathfrak{s}^3} = i \delta_{2,2k}
  \frac{\partial}{\partial \mathfrak{s}} g_{n_0} |_{\mathfrak{s} =0} 
  = i \delta_{2,2k}
  \frac{g_{n_0}^\prime (1)}{-\sqrt{\Omega (1)}}  = i \delta_{2,2k} 
$$
  so the statement holds
  for $j=2$ and any $k $ with $2 k \geqslant j$.
  The rest of the proof is very similar to that for $\tau=0$.
\end{proof}
\z Let $\psi_{1,k} $, $\psi_{2,k}$ be two independent solutions of
\begin{equation}
\label{n5}
\mathcal{L}_k \psi = 0 ~~; {\rm and}~W_k =
\psi_{1,k} (r) \psi_{2,k}^\prime (r) - \psi_{2,k} (r) \psi_{1,k}^\prime (r)
\end{equation}
where
\begin{equation}
\label{mathLk}
\mathcal{L}_k \psi = \psi^{\prime \prime} + 
Q_k
\psi  
\end{equation}
{From} the form of the equation we see that $W_k$ is independent of $r$.
\begin{Lemma}
\label{l1}
For $n = n_0-k$, $k \geqslant 1$, the system (\ref{eqhom2}) is equivalent to
\begin{equation}
\label{n3}
g_{n_0-k} (r) = i \int_r^1 \Omega (s) \left ( g_{n_0-k+1} (s)
- g_{n_0-k-1} (s) \right ) G_k (r, s) ds\ \ \ k\geqslant 1
\end{equation}
where \index{$G_k $}
\begin{equation}
\label{n4}
G_k (r, s) = W_k^{-1}[{\psi_{1,k} (r) \psi_{2,k} (s) -
\psi_{2,k} (r) \psi_{1,k} (s) }]
\end{equation}
\end{Lemma}
\begin{proof} The proof simply follows from variation of parameters,
  the two boundary conditions at $r=1$ and
  $g_{n_0-k}(1)=g'_{n_0-k}(1)=0$.
\end{proof}

\begin{Definition}
\label{defjk}
Define \index{$ j_k$}
\begin{equation}
j_k = \frac{\mathfrak{s}}{m_k} \left [\mathcal{L}_k m_k
  -\Omega m_{k-1} \right ]
\end{equation}
\end{Definition}
\begin{Lemma}
\label{l2}
 For $ k \ge 1$, 
there exist constants $C_1$, $C_2$ and $c_*$, independent of
$k$ so that for any $r \in (0, 1]$ we have  $|j_k| \le c_*$ . For $ r \ge
\frac{1}{k} $, we have $| j_k^\prime (r) | \le C_1/\left (kr^2 \right ) + C_2 $   
\end{Lemma}
\begin{proof}
In the Appendix, (\ref{An7}), we obtain an explicit expression
for $j_k$. Routine asymptotics for large
$k$ in different regimes of $r \in (0, 1] $, discussed in the
Appendix \S \ref{jkcal}, show that 
$k^2 j_{k}^{(2)} + k j_k^{(1)} = O(1)$ in all cases and hence $j_k = O(1)$.
In fact,  as $r \rightarrow 0$ and $k\to\infty$ 
 with $\zeta = k \alpha r = O(1)$ fixed, we have
$j_k \rightarrow g(\zeta)$,  where $g(\zeta)$ is bounded.
Also taking the $r$- derivative of $j_k$ for
$r =O(1)$ not small, we get $j_k^\prime (r) = O(1)$. 
When $r \ll 1$, the asymptotics
in the regime $ \frac{1}{k} \ll r \ll 1$  gives
$j_k = O\left (\zeta^{-1} \right ) = O \left ( 1/(kr) \right )$. 
Since the asymptotics is differentiable, we have  
$j_k^\prime (r) = O \left ( 1/(k r^2) \right ) $. Finally, we look at
 $\zeta = O(1)$, $\zeta \ge 1$. 
Since $\frac{d}{dr} j_k 
= k \frac{d}{d\zeta} j_k  \sim k g^\prime (\zeta)$ where
$\zeta^2 g^\prime (\zeta) $ is bounded for all $\zeta$, it follows
that  $| j_k^\prime (r) | \le C_1/\left ( kr^2 \right ) + C_2 $  for
 $ r \ge
1/k $.  
\end{proof}
\begin{Lemma}
\label{l3}
For $k \geqslant 1$,
$h_k (r)$ defined in (\ref{n9})
satisfies the system of differential equations:
\begin{multline}
\label{S7.7}
h_k^{\prime \prime} + 2 h_k^\prime \frac{m_k^\prime}{m_k}
+ \left ( \frac{\Omega m_{k-1} }{ m_k} + \frac{j_k}{\mathfrak{s}} \right ) h_k=
\Omega \left ( \frac{m_{k-1}}{m_k} h_{k-1} (r) + \frac{m_{k+1}}{m_k} h_{k+1} (r) \right ),
\end{multline}
and the system of integral 
equations (\ref{n3}) is equivalent, for $k\geqslant 1$,
to
\begin{multline}
\label{S7.9}
h_k (r) =
\int_r^1 \frac{\Omega (s) m_{k-1} (s)}{m_k (r)} G_k (r,s)
h_{k-1} (s) ds
\\+ \int_r^1
\frac{\Omega (s) m_{k+1} (s)}{m_k (r)} G_k (r,s)
h_{k+1} (s) ds := \mathcal{A}_k h_{k-1} + \mathcal{H}_k h_{k+1}
\end{multline}  \index{$  \mathcal{A}_k $}\index{$ \mathcal{H}_k $} 
\end{Lemma}

\begin{proof}
This simply follows by substituting $g_{n_0-k} (r) =i^k m_k (r) h_k (r)$
into (\ref{eqhom2}) and (\ref{n3}), and using
$$ \frac{m_k^{\prime \prime}}{m_k} + Q_k = 
\frac{\mathcal{L}_k m_k}{m_k} = \frac{\Omega m_{k-1} }{ m_k} +
\frac{j_k}{\mathfrak{s}},
$$
in turn a consequence of Lemma \ref{l2}.
\end{proof}

\begin{Remark}\label{1stremark}{\rm 
 Let now  $ r \in [{\hat \epsilon}, 1]$, 
where ${\hat \epsilon} \geqslant {C_2}k^{-1}$ 
for sufficiently large $C_2$
independent of $k$. It is convenient to rewrite
$\mathcal{A}_k $ and
$\mathcal{H}_k $ in
(\ref{S7.9}) in terms of $\mathfrak{s}$ (see (\ref{eqxi})).
Furthermore, changing the variable of integration from $s$ to
$ t = {\mathfrak{s} (s)}/{\mathfrak{s} (r)}$, we obtain \index{$ T_k $}
\begin{equation}
\label{S7.20}
[\mathcal{A}_k h_{k-1}] (\mathfrak{s})
= (2k+\tau) (2k +\tau -1 )
\int_0^1 t^{2 k - 2+\tau} 
T_k (\mathfrak{s}, t) h_{k-1} (\mathfrak{s} t) dt 
\end{equation}
where,  using (\ref{S7.3}), we get
\begin{equation}
\label{eqTkdef}
T_k (\mathfrak{s}, t) = \frac{\sqrt{\Omega ( r (\mathfrak{s} t) )}
F_{k-1} (r(\mathfrak{s} t)}{
\mathfrak{s} F_k (r( \mathfrak{s}))}
G_k (r(\mathfrak{s}),r(\mathfrak{s} t))
\end{equation}
and 
\begin{multline}
\label{S7.20.5}
[\mathcal{H}_k h_{k+1}] (\mathfrak{s})= 
\frac{\mathfrak{s}^3}{(2 k+2+\tau)(2k+1+\tau)}\\ \times 
\int_0^1 \sqrt{\Omega (r(\mathfrak{s} t)} t^{2 k+2+\tau}
\frac{F_{k+1} (r(\mathfrak{s} t)}{F_k (r(\mathfrak{s}))}
G_k (r(\mathfrak{s}),r(\mathfrak{s} t))
h_{k+1} (\mathfrak{s} t) d t 
\end{multline}
In evaluating $\mathcal{A}_k$ for large $k$, it is useful to calculate the
Taylor expansion of $T_k (\mathfrak{s}, t)$ 
and its $\mathfrak{s}$ derivative at $t=1$. 
 To do so, we first note that
\begin{equation}
\label{new.tk1}
\frac{\partial T_k}{\partial t} = 
\left ( -\frac{\Omega^\prime (r') F_{k-1} (r') }{2 \Omega (r') F_k (r)}  
- \frac{F_{k-1}^\prime (r')}{F_k (r)} \right ) G_k (r, r') 
- \frac{F_{k-1} (r')}{F_k (r)} \frac{\partial G_k}{\partial r'} (r, r')  
\end{equation}
where, to simplify notation, we wrote  $r(\mathfrak{s} ) = r$  and
$r( t \mathfrak{s} ) = r'$ and used 
$ \partial_{r'} \mathfrak{s} (r') = - 
\sqrt{\Omega (r')} $.    
{From} (\ref{n5}) and (\ref{n4}) we get 
$G_k (r, r) =0$ and $\partial_{r'} G_k (r, r') = 1$ at $r'=r$;
(\ref{new.tk1}) implies
\begin{equation}
\label{new.tk2}
\frac{\partial T_k}{\partial t} \Big\rvert_{t=1}
= -\frac{F_{k-1} (r)}{F_k (r)} 
\end{equation}
Using (\ref{new.tk1}), taking an additional derivative with respect
to $t$, using also (\ref{eqQk}) and (\ref{mathLk}) to see that
$\partial_{r' r'} G_k = -Q_k G_k$, we obtain
\begin{equation}
\label{new.tk3}
\frac{\partial^2 T_k}{\partial t^2} \Big\rvert_{t=1}
= \frac{F_{k-1} (r)}{F_k (r) } 
\left ( \frac{ \xi \Omega^\prime (r) }{2 \Omega^{3/2} (r) }  
+ \frac{2 \xi F_{k-1}^\prime (r)}{\sqrt{\Omega (r)} F_{k-1} (r)} \right ) 
\end{equation}
A similar calculation can be carried out for the third derivative. We
only write down the  potentially largest term in the regime 
$kr \ge C_2$ (for large $k$ and small $r$)
\begin{equation}
\label{new.tk4}
\frac{\partial^3 T_k}{\partial t^3} \Big\rvert_{t=1}
= 
\frac{\xi^2 F_{k-1} (r) Q_k (r)}{\Omega (r) F_k (r) } 
+ O \left ( 1, \frac{1}{k r^3} \right )   
= \frac{\xi^2 F_{k-1} (r) }{\Omega (r) F_k (r) } 
\left ( k \omega - \frac{l (l+1}{r^2} \right ) + 
O \left ( 1, \frac{1}{k r^3} \right ) .    
\end{equation}
Note that if $k r$ is sufficiently large, (\ref{S7.6}) gives
\begin{equation}
\label{new.tk41}
\frac{F_{k-1} (r)}{F_k (r) } = \frac{H (\alpha (k-1) r) H(\alpha k)}{
H(\alpha k r) H(\alpha (k-1)} = 1 + O (k^{-2} r^{-1} ) 
\end{equation}
and 
\begin{equation}
\label{new.tk5}
\frac{F_{k-1}^\prime (r)}{F_{k-1} (r) } 
\sim  - \frac{l (l+1)}{2 \alpha k r^2} + O \left ( \frac{1}{k^2 r^3} \right )  
\end{equation}
Note also that
(\ref{eq:defquants}) implies
$\alpha - 2 \sqrt{\Omega (r)}/\left ( \mathfrak{s}(r) \right ) = O(r)$ for
small $r$. Including all terms that become important when $r $ is small,
we note that
in the regime when $k r$ is sufficiently large, we have
\begin{multline}
\label{S7.25}
T_k = (1-t) + \left ( - \frac{k}{4} f_1
+ \frac{f_2}{r^2} \right )
\left ( \frac{2}{3} (1-t)^3 - \frac{(1-t)^2}{k} \right )
\\+ O \left (\frac{(1-t)^4}{r^3}, \frac{(1-t)^3}{k r^3},
\frac{(1-t)^3}{r}, \frac{(1-t)^2}{k r} , 
\frac{(1-t)^2}{k^2 r^3} \frac{(1-t)}{k^2 r} \right )
\end{multline}
\begin{multline}
\label{S7.25.2}
\frac{\partial T_k}{\partial \mathfrak{s}}  =
\left ( - \frac{k}{4} f_1^\prime
+ \frac{f_3}{r^3} \right )
\left ( \frac{2}{3} (1-t)^3 - \frac{(1-t)^2}{k} \right )
\\+ O \left (\frac{(1-t)^4}{r^4}, \frac{(1-t)^3}{k r^4},
\frac{(1-t)^3}{r^2}, \frac{(1-t)^2}{k r^2}, \frac{(1-t)^2}{k^2 r^4},   
\frac{(1-t)}{k^2 r^2} \right )
\end{multline}
where \index{$f_1 $}
\begin{equation}
\label{S7.27}
f_1 (\mathfrak{s}) = \frac{\omega \mathfrak{s}^2}{\Omega (r(\mathfrak{s}))}
\end{equation}\index{$f_2 $}
\begin{equation}
\label{7.27.1.5}
f_2 (\mathfrak{s}) =
\frac{l(l+1) \mathfrak{s}^2 }{4 \Omega}
\end{equation}\index{$f_3 $}
\begin{equation}
\label{7.27.1}
f_3 (\mathfrak{s}) =
\frac{l(l+1) \mathfrak{s}^2 }{2 \Omega^{3/2}}
\end{equation}
When $r \in [0, {\hat \epsilon}]$, for ${\hat \epsilon} = C_2/k$, 
it is sometimes more convenient to express
$\mathcal{A}_k $ in terms of $\zeta = k \alpha r$. For that purpose,
we define
\begin{equation}
\label{n10}
Q(\zeta) = -2k \log \left [ 1 - \frac{\mathfrak{s} (0)\!\! - \!\!\mathfrak{s}}{\mathfrak{s} (0)} \right ]
-\log \left [ \left (\frac{\Omega (0)}{\Omega (r)} 
\right )^{\frac{1}{4}}\!\!\!
\exp \left ( \frac{1}{4}
\int_0^r \!\!\!\! dr' \frac{\omega \mathfrak{s} (r')}{\sqrt{\Omega(r')}} 
\right ) \right ],
\end{equation} \index{$ Q(\zeta)$}

\noindent where we recall the relation (\ref{defmk}) 
between $\mathfrak{s}$ and $r =\zeta/(k\alpha)$,
$ \zeta \in [0, k \alpha \epsilon]$. A series expansion in $k^{-1}$ 
leads to
\begin{equation}
\label{n11}
Q(\zeta) = \zeta - \frac{\zeta}{k} \left (
\frac{\omega}{2 \alpha^2} - \frac{\Omega^\prime (0)}{4 \Omega (0) \alpha}
\right ) + \frac{\zeta^2}{4 k} \left ( 1 +
\frac{\Omega^\prime (0)}{\alpha \Omega (0)} \right )
+ O\left (\frac{\zeta^3}{k^2}, \frac{\zeta^2}{k^2} \right )
\end{equation}
We choose ${\hat \epsilon}_1 = {\tilde C}_2 k^{-1}\log k$,
for some  $k$-independent ${\tilde C}_2$ (chosen more precisely later).
We define  
${\hat \delta}_1$, dependent or $r$, so that
\begin{equation}
\label{eqdelta1}
(1-{\hat \delta}_1) \mathfrak{s} (r)  
= \mathfrak{s} (\epsilon_1) 
\end{equation} 
 From (\ref{defmk}), it follows 
that for sufficiently large ${\tilde C}_2$   we have
\begin{equation}
\label{eqdelta11}
{\hat \delta}_1 \geqslant 1- 
\frac{\mathfrak{s} (\epsilon_1)}{\mathfrak{s} (\epsilon)} \geqslant
\frac{(5+l) \log k}{(4k+2\tau) k},
\end{equation}
It follows from the definition of $\mathcal{A}_k$ in (\ref{S7.9}) that
for $r \in [0, {\hat \epsilon}]$, {\it i.e.} $\zeta \in [0, k 
\alpha {\hat \epsilon}]$,
\begin{multline}
\label{n12}
\left [ \mathcal{A}_k h_{k-1} \right ] (\zeta)=
\int_\zeta^{k \alpha {\hat \epsilon}_1} \!\!\!\!\!\! e^{-Q(\eta) + Q(\zeta)}
\left ( 1 + \frac{a_1}{k} \right )
\frac{H(\eta (1-k^{-1}))}{H(\zeta)}
\mathcal{G} (\zeta, \eta) h_{k-1} (\eta (1-k^{-1}) ) d\eta \\
+ (2k+\tau) (2k+\tau-1)
\\\times \int_{{\hat \epsilon}_1}^{1} 
\left [ \frac{\mathfrak{s} (r')}{\mathfrak{s} (r)} \right ]^{2k-2+\tau}
\Omega (r')  G_k \left (r,
r' \right ) \frac{F_{k-1} (r')}{\mathfrak{s}^2 F_{k} (r) }
h_{k-1} (r') dr'\\ =: [\mathcal{A}_k^{0} h_{k-1}](\zeta) +
[\mathcal{A}_k^1 h_{k-1} ] (r)
\end{multline}
\index{$ \mathcal{A}_k^{0} $}\index{$\mathcal{A}_k^1 $}
where $\mathcal{G}(\zeta, \eta)$ is defined by \index{$ \mathcal{G}(\zeta, \eta)$}
\begin{equation}
\label{n13a}
\mathcal{G}(\zeta, \eta) = k \alpha
G_k (r(\zeta), r(\eta)) 
\end{equation}
while
\begin{multline}
\label{n13}
\frac{a_1 (\eta,\zeta)}{k} = 
\frac{H (\alpha k)}{H (\alpha (k-1))} 
\left ( 1 + \frac{\tau}{2k} \right ) 
\left ( 1 + \frac{\tau-1}{2k} \right ) 
\frac{\mathfrak{s}^2 (0) \Omega (\eta/(k\alpha)) }{ 
      \mathfrak{s}^2 (\eta/(k\alpha)) \Omega (0) } \\
\times  
\left ( \frac{\mathfrak{s} (\eta/(k\alpha)) }{
\mathfrak{s} (\zeta/(k\alpha)) }  
\right )^{\tau}
- 1 
\end{multline}
\index{$ a_1$}
while  for large $k$ and $0<\zeta\le \eta\le {\hat \epsilon}_1\alpha$ we have
\begin{equation}
a_1 (\eta,\zeta) = \tau - \frac{1}{2} + 
\left ( 1 + \frac{\Omega^\prime (0)}{\alpha \Omega (0)} \right )
\eta + 
\frac{\tau}{2} (\zeta -\eta) +
O \left ( \frac{\eta^2}{k}, \frac{\eta}{k} \right )
\end{equation}
Similarly, for $k \alpha r = \zeta
\in (0, k {\hat \epsilon}_1 \alpha)$, we define \index{$ b_1 $}
\begin{equation}
\frac{b_1 (\eta,\zeta)}{k} = 
\frac{H (\alpha k)}{H (\alpha (k+1))} 
\frac{\mathfrak{s}^2 (\eta/(k\alpha) ) \Omega (\eta/(k\alpha)) }{ 
      \mathfrak{s}^2 (0) \Omega (0) } 
\left ( \frac{\mathfrak{s} (\eta/(k\alpha)) }{
\mathfrak{s} (\zeta/(k\alpha)) }  
\right )^{\tau}
- 1 
\end{equation}
We then have
\begin{multline}
\label{n12.5}
\left [ \mathcal{H}_k h_{k+1} \right ] (\zeta) = 
\frac{\Omega (0) \mathfrak{s}^2 (0)}{\alpha^2 k^2 (2k+2+\tau) (2 k+1 +\tau)}
\\\times
\int_\zeta^{k \alpha {\hat \epsilon}_1} \!\!\!\!\!\!\!\!e^{-Q(\eta) + Q(\zeta)}
\left ( 1 + \frac{b_1}{k} \right )
\frac{H(\eta (1+k^{-1}))}{H(\zeta)}
\mathcal{G} (\zeta, \eta) h_{k+1} (\eta (1+k^{-1}) ) d\eta \\
+ \frac{\mathfrak{s}^2}{(2k+2) (2k +1+2\tau)}\\\times 
\int_{0}^{1-{\hat \delta}_1}
\sqrt{\Omega (r(\mathfrak{s} t))} G_k \left (r (\mathfrak{s}),
r(\mathfrak{s} t) \right ) t^{2 k+2+\tau} \frac{F_{k+1} (r(\mathfrak{s} t))}{F_{k} (r (\mathfrak{s})) }
h_{k+1} (\mathfrak{s} t) dt\\=: [\mathcal{H}_k^{0} h_{k+1}] +
[\mathcal{H}_k^1 h_{k+1} ]
\end{multline}}
\end{Remark}

\begin{Lemma}
\label{l2.5}
For $k \geqslant 2$ and $k_1\in\{k-1, k,k+1\}$ we have
\begin{enumerate}

\item{} If $r \in (0,1)$ and $s \in (r, r+\delta)$, where
$\delta \le \min \left \{ {C_2} k^{-1}\log k, 1-r \right \}$, then
$$
\Big\lvert G_k (r, s) \frac{F_{k_1} (s)}{F_k (r)} \Big\rvert \leqslant 
\frac{c_*}{k^{1/2}},
\ \ \Big\lvert \frac{\partial}{\partial r} \left (G_k (r, s)
  \frac{F_{k_1} (s)}{F_k (r)} \right ) \Big\rvert < c_* k^{1/2}$$

\item{} If $r \in (0, 1)$, $\delta \le C_2 k^{-1}\log k$ with $r+\delta < 1$,
  then for $s \in (r+\delta, 1)$, 
  $$
  \Big\lvert G_k (r, s) \frac{F_{k_1} (s)}{F_k (r)}
  \Big\rvert < c_* k^{l/2-1/2}, \Big\lvert \frac{\partial}{\partial r} \left (G_k (r, s)
    \frac{F_{k_1} (s)}{F_k (r)} \right ) \Big\rvert < c_* k^{l/2+1/2}$$

\end{enumerate}

\end{Lemma}

\begin{proof}
It suffices to find
bounds for $G_k (r, s) H(\alpha k_1 s)/H(\alpha k r)$  
since the other functions  involved are regular everywhere for
$r , s \in [0,1]$, see (\ref{S7.3}).
We first consider  $k \rightarrow +\infty$.

It is easily verified
that $\mathcal{G}(\zeta, \eta)$, defined  in (\ref{n13a}), 
is the Green's function (see (\ref{eqQk}), (\ref{mathLk}) ) for
\begin{equation}
\label{n11.1}
\mathcal{L}:=\Psi\mapsto\Psi^{\prime \prime} - \frac{l (l+1)}{\zeta^2} \Psi
+ \frac{\Psi}{k} \left [ \frac{\omega}{\alpha^2} + \frac{b}{\alpha \zeta}
\right ]
+\frac{\Psi}{k^2 \alpha^2} \left [
ip_1 - n_0 \omega \right ]
\end{equation}
and is given by
\begin{equation}
\label{n10new}
\mathcal{G} (\zeta, \eta)
:= k \alpha G_k (r (\zeta), r(\eta)) = \frac{\Psi_1 (\zeta)
\Psi_2 (\eta) - \Psi_2 (\zeta) \Psi_1 (\eta) }{W},
\end{equation}
where $\Psi_1$, $\Psi_2$ are two independent solution of $\mathcal{L}\Psi=0$
and $W=\Psi_1 (\zeta) \Psi_2^\prime (\zeta) -
\Psi_2 (\zeta) \Psi_1^\prime (\zeta)$ is their constant Wronskian.

Standard asymptotic results show there exist two independent solutions
$\Psi_1$, $\Psi_2$ such that for large $k$,
we have uniformly in $z \in [0,\sqrt{\omega k}]$
\begin{equation}
\label{eqPsi1}
\Psi_1 \sim
-\frac{2^l l!}{(2 l)!} \sqrt{\frac{\pi z}{2} } Y_{l+1/2} (z) ~;
~{\rm where}
~z=\sqrt{\frac{\omega}{\alpha^2 k} } \zeta = \sqrt{\omega k} r
\end{equation}
\begin{equation}
\label{eqPsi2}
\Psi_2 \sim
\frac{2^{-l-1} (2 l+ 2)!}{(l+1)!} \sqrt{\frac{\pi z}{2} } J_{l+1/2} (z)
\end{equation}
The Wronskian $W$ is asymptotic, for large $k$, to $(2 l
+1)\sqrt{\omega}/\sqrt{\alpha^2 k} $.  The expressions (\ref{eqPsi1}) and
(\ref{eqPsi2}) may also be used to determine the asymptotics of
$\Psi_1^\prime$ and $\Psi_2^\prime$.  Using (\ref{n10new}), (\ref{eqPsi1}),
(\ref{eqPsi2}) and (\ref{S7.3}) and the bounds on $W$, with $l_1=l+\frac{1}{2}$ it follows that
\begin{multline}
\label{eqGkFkbound}
\left\lvert \frac{F_{k_1} (s)}{F_k (r)}  G_k (r, s) \right\rvert 
\leqslant  \frac{c_*|z z'|^{1/2}}{k^{1/2}}\Bigg \lvert 
\frac{H \left (\alpha \sqrt{\frac{k_1}{\omega}} z' 
\right)}{  
H\left (\alpha \sqrt{\frac{k}{\omega}} z \right )} 
[Y_{l_1} (z) J_{l_1} (z') 
- J_{l_1} (z) Y_{l_1} (z')]
\Bigg\rvert
\end{multline}
where $z'=\eta\sqrt{\omega}/\sqrt{\alpha^2 k} = \sqrt{\omega k} s $.  A
similar bound holds for
$$ \left\lvert \frac{\partial}{\partial r} 
\left \{ \frac{F_{k_1} (s)}{F_k (r)}  G_k (r, s) \right \} \right\rvert $$

We now prove part (1). We break this case up into two subcases: (a) $r \in
[k^{-2/3}, 1]$ and (b) $r \in [0, k^{-2/3}]$. In case (a), we note that $s \in
[r, r+\delta]$ implies $s/r $ and therefore $1 \le z'/z = O(1)$. The function
$H$ in (\ref{eqGkFkbound}) is close to $1$ because its argument is large.
Furthermore, note that $\sqrt{z} Y_{l+{1/2}} (z)$ and $\sqrt{z} J_{l+1/2} (z)$
are bounded for large $z$, while they are asymptotic to constant multiples of
$z^{-l}$ and $z^{l+1}$ for small $z$.  Using (\ref{eqGkFkbound}), part 1 of
the Lemma follows by inspection in  case (a).  For the case (b),
(\ref{eqGkFkbound}) further simplifies since $z$,$z'$ are small and
\begin{multline}
\label{eq12.37}
  \frac{H(k_1 \eta/k)}{H(\zeta)} G_k (r(\zeta), r(\eta)) =
\frac{H(k_1 \eta/k)}{k\alpha H(\zeta)}
\mathcal{G}(\zeta, \eta)\sim
\frac{H(k_1 \eta/k) \left ( \eta^{l+1} \zeta^{-l} -
\zeta^{l+1} \eta^{-l} \right )}{k \alpha H(\zeta)
(2 l+1)}
\end{multline}
 When $\zeta \in [\log k, \alpha k^{1/3}]$ and $\eta \in 
[\zeta, \zeta + \alpha k \delta ]$, we have
$1 \le [\eta/\zeta]^l \le c_* $ and therefore
$$
\left \lvert \frac{H(k_1 \eta/k)}{H(\zeta)} G_k (r(\zeta), r(\eta))
\right\rvert = \left\lvert \frac{H(k_1 \eta/k)}{k\alpha H(\zeta)}
  \mathcal{G}(\zeta, \eta) \right\rvert \leqslant  \frac{c_*}{k^{1/2}}$$
The same
inequality holds if $\zeta \in [0, \log k]$, since $\eta \in [\zeta, (C_2+1)
\log k] $ since in this regime $\zeta^{-l}/ H(\zeta)$ is bounded and the
logarithmic growth in $k$ of terms involving $\eta$ can be bounded by, say,
$k^{1/2}$,  while for small $\eta$,
$\eta^l H \left ( k_1 \eta/k \right )$
is bounded.  The bounds on derivatives follow in a similar manner using
$\frac{d}{dr} = k \alpha \frac{d}{d\zeta}$.

Part 2 (which is only relevant for $r +\delta \leqslant  1$) follows similarly on
careful inspection of (\ref{eqGkFkbound}), from the asymptotic behavior in
different regimes of $z$ and $z'$.
\end{proof}

\begin{Lemma}
\label{14.0.0}
Let 
$r \in (0, {\hat \epsilon} ]$, with ${\hat \epsilon}_1
= 
\frac{C_2}{k} \log k$. We choose $C_2$ large enough so that
$\frac{\mathfrak{s} (\epsilon_1) }{\mathfrak{s} (r) }  
= (1- \delta_1) \le \frac{(5+l) \log k}{4k + 2 \tau} $.  Then
$ |[\mathcal{A}_k^1 f] (r) | \leqslant  
c_* k^{l/2-1/2} (1-\delta_1)^{2k-2+\tau} \| f \|_\infty
\leqslant  c_* k^{-3}\|f \|_\infty$ and
$ |\frac{d}{dr} [\mathcal{A}_k^1 f] (r) | \leqslant 
c_* k^{l/2+1/2} (1-\delta_1)^{2k-2+\tau} \| f \|_\infty
\leqslant  c_*k^{-2} \|f \|_\infty$.

\end{Lemma}

\begin{proof} 
  Consider $\mathcal{A}_k^1$ given by (\ref{n12}).  We note that
  $\mathfrak{s}^{-2}\Omega(s)$ and its $r-$derivative are bounded, while $ G_k (s, r )
  {F_k (s)}/{F_k (r)}$ and its $r$-derivative are bounded by $ c_* k^{l/2-1/2}$
  and $ c_* k^{l/2+1/2}$ respectively for any $\tau$ (cf. Lemma \ref{l2.5}).
  Further $\lvert {\mathfrak{s} (s)}/{\mathfrak{s} (r)} \rvert \leqslant  (1-\delta_1)$ and from
  (\ref{eqdelta11}), we have
$$ (1-\delta_1)^{2k-2+\tau}  \leqslant  \frac{c_*}{k^{l/2+5/2}}, $$
and  the lemma follows.   
\end{proof}

\begin{Remark}
  Since for $r \in (0, {\hat \epsilon}]$, the bound in Lemma \ref{14.0.0} on
  $\mathcal{A}_k^1$ is $O(k^{-2})$, we will see later that $\mathcal{A}_k$ is
  dominated by $\mathcal{A}_k^0$ (defined in
(\ref{n12})) as $k\to\infty$.
\end{Remark}

\begin{Lemma}\index{$\mathcal{G}_0 (\zeta, \eta) $}
\label{l4.0.0}
Define $\mathcal{G}_0 (\zeta, \eta) =\lim_{k \rightarrow \infty} \mathcal{G}
(\zeta, \eta)$ and 
$ H_0 (\zeta)$ $=$ $\lim_{k \rightarrow \infty}$ 
$ H(\zeta)$, where
$\zeta, \eta \ll k^{1/2}$ as $k \rightarrow \infty$.  Then,
\begin{equation}
\int_{\zeta}^\infty e^{-\eta + \zeta} \mathcal{G}_0 (\zeta, \eta)
\frac{H_0(\eta)}{H_0(\zeta)} d\eta = 1
\end{equation}
\begin{equation}
\int_{\zeta}^\infty e^{-\eta + \zeta} {\mathcal{G}_0}_\zeta (\zeta, \eta)
\frac{H_0(\eta)}{H_0(\zeta)} d\eta=-1 + \frac{H_0^\prime (\zeta)}{H_0(\zeta)}
\end{equation}

\end{Lemma}

\begin{proof}
Using (\ref{eqPsi1}) and (\ref{eqPsi2}) and the
behavior of Bessel functions for small argument, \cite{Abramowitz}, it follows 
that for
$\zeta , \eta \ll k^{1/2} $ we have
\begin{equation}
  \label{eq:nonum00}
 \mathcal{G}_0 (\zeta, \eta) = \lim_{k\rightarrow \infty}
\mathcal{G} (\zeta, \eta) =
\frac{ {\eta}^{l+1} \zeta^{-l}
- \zeta^{l+1} {\eta}^{-l} }{2 l + 1}
 \end{equation}
and $H_0 (\zeta) = \lim_{k \rightarrow \infty} H(\zeta) = 
\sqrt{\frac{2}{\pi}} \zeta^{1/2}
e^{\zeta} K_{l+1/2} (\zeta)$.  Now, using the modified Bessel function
equation, it is easily verified that $ f (\zeta) = e^{-\zeta} H_0 (\zeta)$
satisfies
$$
f^{\prime \prime} - \frac{l (l+1)}{\zeta^2} f = f $$
with $f(\zeta) \sim
e^{-\zeta}$ as $\zeta \rightarrow \infty$.  Using variation of parameters to
invert the left hand side of the above equation, and using the boundary
conditions at $\infty$ we obtain
$$
f (\zeta) = \int_\zeta^\infty \mathcal{G}_0 (\zeta, \eta) f(\eta) d\eta $$
Dividing through by $f(\zeta)$, the first identity in the Lemma follows. By
differentiating the first identity with respect to $\zeta$, and using the
first identity in the resulting expression, we obtain the second identity.
\end{proof}

\begin{Lemma}
\label{l4.0.0.0}
For any $r \in (0, 1)$,
\begin{equation}
\label{claim}
\Big\lvert
\mathcal{A}_k [1] (r) - 1\Big | = \Big\lvert
\int_{r}^1 \Omega (s) \frac{m_{k-1} (s)}{m_k (r)} G_k (r, s) ds - 1 \Big\rvert
\leqslant  \frac{c_*}{k^2}.
\end{equation}
 For $ \frac{1}{k} \le r \le \frac{1}{2}$ we get
\begin{equation}
\label{claim1}
\Big\lvert \frac{d}{dr} \mathcal{A}_k [1] (r) \Big\rvert =
\Big\lvert \frac{d}{dr}
\int_r^1 \Omega (s) \frac{m_{k-1} (s)}{m_k (r)} G_k (r, s) ds \Big\rvert
\leqslant  \frac{c_*}{k^2} + \frac{c_*}{k^3 r^2} , 
\end{equation}
while for any $r \in [0, \frac{1}{2} ] $,
\begin{equation}
\label{lplem44}
\int_{r}^1\!\!\!\! \Omega (s) 
\Big\lvert \frac{\partial}{\partial r}
\left ( G_k (r, s) \frac{m_{k-1} (s)}{m_k (r)} \right ) \Big\rvert ds \leqslant  c_* k \\
\end{equation}
\end{Lemma}
\begin{proof} Recalling the definition (\ref{S7.9}),   it follows 
from (\ref{defjk}) and Lemma \ref{l2}  that 
\begin{equation}
\label{S7.14}
\mathcal{L}_k m_k - \Omega m_{k-1} = \frac{j_k(r)}{\mathfrak{s}} m_k
\end{equation}
where $j_k(r) = O(1)$ as $k \rightarrow +\infty$ for any $r \in [0,1]$.  We
can  check from (\ref{S7.3})
that $m_k (1) = 0$, $m_k^\prime (1) =0$ for $k \geqslant 1$. 
{From} (\ref{S7.14}), inversion of $\mathcal{L}_k$ yields
\begin{equation}
\label{S7.15}
m_k (r) = \int_{r}^1 G_k (r, s) \left \{ \Omega (s) m_{k-1} (s) +
\frac{j_k (s)}{\mathfrak{s} (s)}
m_k (s)
\right \} ds
\end{equation}
Therefore,
\begin{equation}
\label{S7.16}
\int_{r}^1 G_k (r, s) \frac{\Omega (s) m_{k-1} (s)}{m_k (r) } ds = 1- \int_{r}^1 G_k (r,s) \frac{j_k (s) m_k (s)}{\mathfrak{s} (s) m_k (r)} ds
\end{equation}
First, we choose ${\hat \delta}_1$ so that
$\displaystyle  1- {\hat \delta}_1 = (5+l) \log
  k / \left ( 4 k+2\tau \right ) $.  
We then define $\hat \delta$ so that 
$(1- {\hat \delta}_1 ) \mathfrak{s} (r) = \mathfrak{s}
(r+{\hat \delta})$. 
It is clear that for large $k$ we have 
$\displaystyle {\hat \delta}  \sim
\left ( 5/2 +
  l/2 \right ) \mathfrak{s} (r) \log k / \left (
(2 k +\tau) \sqrt{\Omega (r) }\right )
$.
Lemma \ref{l2.5}, and the fact that $k^{l/2-1/2} (1-{\hat
    \delta}_1)^{2k+1+\tau} /( 2k+1+\tau) 
\leqslant  \frac{1}{k^3} $
 give
\begin{multline}
\label{eq12.45}
   \Big\lvert
\int_{r+{\hat \delta}}^1 G_k (r, s)
\frac{j_k (s) m_k (s) }{\mathfrak{s} (s) m_k (r)} ds \Big\rvert\leqslant  
\Big\lvert \int_{0}^{1-{\hat \delta}_1} t^{2k+\tau} \\ \times\frac{F_k (r(\mathfrak{s} t)}{
\sqrt{\Omega (r(\mathfrak{s} t) ) } F_k (r(\mathfrak{s}) }
G_k (r (\mathfrak{s}), r(\mathfrak{s} t) ) j_k (r (\mathfrak{s} t))  dt \Big\rvert
\leqslant  \frac{c_*}{k^3} \|j_k \|_\infty \leqslant  \frac{c_*}{k^3}
\end{multline}
Now, consider the contribution from $\int_r^{r +{\hat \delta}}$.
There are again two cases: (i) $1 \geqslant r \geqslant k^{-2/3} $ and (ii)
$ 0 < r \leqslant  k^{-2/3} $.

In the first case, Taylor expanding 
$G_k (r, s)$ near $s=r$ we get $G_k = (s-r) + O((s-r)^3 Q_k)
= \sqrt{\Omega (r)} \mathfrak{s} (1-t) + O(k^{4/3} (1-t)^3, (1-t)^2)$.
Hence,
\begin{equation}
\label{eq12.46}
\Big\lvert
\int_{r}^{r+{\hat \delta}} G_k (r, s) \frac{j_k (s) m_k (s) }{\mathfrak{s}(s) m_k (r)} ds \Big\rvert
\leqslant  c_* \| j_k \|_\infty \int_{1-{\hat \delta}_1}^1 t^{2 k+\tau-1} (1-t) dt
\leqslant  \frac{c_*}{k^2}
\end{equation}
For the case (ii), we rewrite the integral in terms of $\zeta = k \alpha r$, to
obtain
\begin{multline}
\label{eq37}
\Big\lvert
\int_{r}^{r+{\hat \delta}} \!\!\!\!G_k (r, s) \frac{j_k (s) m_k (s) }{\mathfrak{s} (s) m_k (r)} 
ds \Big\rvert
\leqslant  \frac{c_*}{k^2} \| j_k \|_\infty \int_\zeta^{\zeta + k \alpha {\hat \delta}}\!\!\!\!
e^{-Q(\eta) + Q(\zeta) }
\mathcal{G} (\zeta, \eta) \frac{H(\eta)}{H(\zeta) } d\eta \\
\leqslant  \frac{c_*}{k^2} \int_\zeta^{\zeta + k \alpha {\hat \delta}} d\eta
e^{-\eta + \zeta} \mathcal{G}_0 (\zeta, \eta) \frac{H_0 (\eta)}{
H_0 (\zeta)}\leqslant  \frac{c_*}{k^2} \int_\zeta^{\infty} d\eta 
e^{-\eta + \zeta} \mathcal{G}_0 (\zeta, \eta) \frac{H_0 (\eta)}{
H_0 (\zeta)} \leqslant  \frac{c_*}{k^2}
\end{multline}
by Lemma 
\ref{l4.0.0}.
Using (\ref{S7.16}) and
(\ref{eq37}), the first part 
follows.

To prove (\ref{claim1}), we note that if $C_3$ is large and 
$ r k > C_3 $,
Taylor expansion gives
\begin{multline}
  U_1 (\mathfrak{s}, t) := \frac{F_k (r(\mathfrak{s} t)}{
\sqrt{ \Omega ( r(\mathfrak{s} t) ) } F_k (r(\mathfrak{s}) }
G_k (r (\mathfrak{s}), r(\mathfrak{s} t) )\\ = f_4 (\mathfrak{s}) (1-t) + 
O\left ( (1-t)^2, k (1-t)^3,
\frac{(1-t)^3}{r^2},
\frac{(1-t)^2}{k r^2}  \right )
\end{multline}
 for $f_4 = -\mathfrak{s}/\Omega (r(\mathfrak{s}))$,
while
$$\frac{\partial}{\partial \mathfrak{s}}  U_1 (\mathfrak{s}, t)
= f_4^\prime (\mathfrak{s}) (1-t) +
O\left ( (1-t)^2, k (1-t)^3,
\frac{(1-t)^3}{r^3},
\frac{(1-t)^2}{k r^3} \right )
$$
{From} (\ref{S7.16}) we note that 
\begin{multline}
\label{Akrdiff}
\frac{d}{dr} \mathcal{A}_k [1] (r) =
-\sqrt{\Omega (r(\mathfrak{s}) ) }  
\left \{ \int_{1-{\hat \delta}_1}^1 
t^{2 k + \tau+1} \frac{j_k^\prime (r (\mathfrak{s} t) )}{
\sqrt{\Omega (r (\mathfrak{s} t) ) } } {U}_1 (\mathfrak{s}, t) dt \right . \\
\left .  
- \int_{1-{\hat \delta}_1}^1 
t^{2 k + \tau} j_k (r (\mathfrak{s} t) ) {U}_{1,\mathfrak{s}} 
(\mathfrak{s}, t) dt \right \}
- \frac{d}{dr} \int_{r+{\hat \delta}}^1 
\left ( \frac{\xi (s)}{\xi (r) } \right )^{2k+\tau} 
\frac{j_k (s) F_{k} (s)}{ \xi (s) F_k (r) } G_k (r, s) ds   
\end{multline}
We note further that
\begin{multline}
\label{Akrdiff2}
- \frac{d}{dr} 
\int_{r+{\hat \delta}}^1 
\left ( \frac{\xi (s)}{\xi (r) } \right )^{2k+\tau} 
\frac{j_k (s) F_{k} (s)}{ \xi (s) F_k (r) } G_k (r, s) ds   
\\= 
-\int_{r+{\hat \delta}}^1 
\left ( \frac{\xi (s)}{\xi (r) } \right )^{2k+\tau} 
\frac{j_k (s)}{\xi (s)} \frac{\partial}{\partial r} \left [
\frac{F_{k} (s)}{F_k (r) } G_k (r, s) \right ] ds \\
+(2k+\tau) \frac{\xi^\prime (r)}{\xi (r)}  
\int_{r+{\hat \delta}}^1 
\left ( \frac{\xi (s)}{\xi (r) } \right )^{2k+\tau} 
\frac{j_k (s) F_{k} (s)}{ \xi (s) F_k (r) } G_k (r, s) ds 
\\
+\left ( \frac{\xi (r+{\hat \delta})}{\xi (r) } \right )^{2k+\tau} 
\frac{j_k (r+{\hat \delta}) F_k (r+{\hat \delta})}{
\xi (r+{\hat \delta}) F_k (r) } 
 G_k (r, r+{\hat \delta} ) \left ( 1 + {\hat \delta}^\prime (r) \right )  
\end{multline}
 From the bounds in Lemmas \ref{l2} and \ref{l2.5} and the fact that 
 $\xi (s)/\xi (r) \le (1-{\hat \delta}_1 )$, 
we easily conclude
that the contribution of $\int_{r+{\hat \delta}}^1$ in 
(\ref{Akrdiff2}) to $\frac{d}{dr} \mathcal{A}_k [1] (r)$ 
is  $O(1/k^2)$. 
   
Since Lemma \ref{l2} implies $|j_k (r)|< c_*$ 
and $|j_k^\prime (r)| < c_* + c_*/(k r^2)$ for $\frac{1}{2} \ge r \ge
\frac{1}{k} $,
it follows from the
local expansion of $U_1 (\mathfrak{s}, t)$ and its $\mathfrak{s}$-derivative
 in a neighborhood of $t =1$
in the first integral in  (\ref{Akrdiff}) that 
$$
\Big\lvert \frac{d}{dr} \int_{r}^{r+{\hat \delta}} G_k (r, s) \frac{j_k (s)
  m_k (s) }{\mathfrak{s} (s) m_k (r)} ds \Big\rvert \leqslant  \frac{c_*}{k^3 r^2} +
\frac{c_*}{k^2} $$
and (\ref{claim1}) 
follows.

We now prove (\ref{lplem44}).  We first note that
for $r \ge k^{-2/3}$, $s \in (r, r+{\hat \delta})$,
from (\ref{n4}), $\partial_r G (r, s) = -1$ at $s=r$ and therefore,
from (\ref{eqPsi1}), (\ref{eqPsi2}), 
it follows that for $s-r =
O(k^{-1} \log k) \ll k^{-1/2}$, $\partial_r G (r, s) \sim -1 < 0$ for
$ s \in (r, r+{\hat \delta})$. 
The same is true for $r \in [0, k^{-2/3}]$ since in this
regime, $\partial_r G_k (r, s) \sim 
\partial_\zeta \mathcal{G}_0 (\zeta, \eta)$ (see  
(\ref{eq:nonum00})), 
with $\zeta= r/(\alpha k)$, $\eta = s/(\alpha k )$.
Therefore, from (\ref{defmk}) and
(\ref{S7.3}),  we get  
\begin{multline}
\label{pGkr}
-\frac{\partial}{\partial r} 
\left ( \frac{m_{k-1} (r) G_k (r,s)}{m_k (r)} \right ) 
=  - \left [ (2k+\tau-2) \frac{\sqrt{\Omega (r)}}{\xi (r)}
 - \frac{F_{k}^\prime (r)}{F_k (r) } \right ]  
\frac{m_{k-1} (r) G_k (r, s) }{m_k (r)} \\ 
- \partial_r G_k (r, s) \frac{m_{k-1} (s)}{m_k (r)} 
\end{multline}
Since the contributions to the integrals 
from $\int_{r+{\hat \delta}_1}^1$
is $O(\frac{1}{k^2} )$, and the first term on the right on (\ref{pGkr})
is negative for large $k$, while the second is positive, 
it follows that 
\begin{multline}
\label{inequalabs}
\int_{r}^1 \Omega (s) \Big\lvert \frac{\partial}{\partial r} 
\left ( \frac{m_{k-1} (r)}{m_k (r) G_k (r, s) } \right )
\Big\rvert ds  \\
\le \big\lvert \frac{d}{dr} \mathcal{A}_k [1] (r) \big\rvert  
+ 2 \left [ (2k+\tau-2) \frac{\sqrt{\Omega (r)}}{\xi (r)}
 - \frac{F_{k}^\prime (r)}{F_k (r) } \right ] \lvert
\mathcal{A}_k [1] (r) \rvert + O \left ( \frac{1}{k^2} \right ) \le c_*k 
\end{multline}  
for $r \in \left [ C_2 k^{-1} , 1 \right ]$. 
For $r \in \left [ 0, C_2 k^{-1} \right ]$, we note that
since the contribution from $\int_{r+{\hat \delta}}^1 $ for
$\frac{d}{dr} \mathcal{A}_k [1] (r)$ is negligible,  we have 
\begin{multline}
\frac{d}{dr} \mathcal{A}_k [1] (r) \sim  
k \alpha \frac{d}{d\zeta} \int_{\zeta}^{\zeta+ k\alpha {\hat \delta}_1 }   
e^{-Q (\eta)+Q(\zeta) } 
\left ( 1 + \frac{a_1}{k} \right ) \frac{H(\eta (1-1/k))}{H(\zeta)} 
\mathcal{G} (\zeta, \eta) \\ 
\sim  
k \alpha \frac{d}{d\zeta} \int_{\zeta}^{\zeta+ k\alpha {\hat \delta}_1 }   
e^{-\eta+\zeta} 
\frac{H_0 (\eta )}{H_0 (\zeta)} 
\mathcal{G}_0 (\zeta, \eta) \sim 
k \alpha \frac{d}{d\zeta} \int_{\zeta}^{\infty}
e^{-\eta+\zeta} 
\frac{H_0 (\eta )}{H_0 (\zeta)} 
\mathcal{G}_0 (\zeta, \eta) 
\end{multline}
it follows immediately from Lemma \ref{l4.0.0} that  
in this case , $ |\frac{d}{dr} \mathcal{A}_k [1] (r) | \le c_* k$.
Hence the inequality in (\ref{inequalabs}) is valid for
 all $r \in [0, 1/2 ]$.
\end{proof}
\begin{Lemma}
\label{l4}
For any $f \in L^\infty [0, 1]$,
\begin{eqnarray}
\label{8180} (a)\ \ \text{~~For}~~r \in [0, 1], ~~~ 
\| \mathcal{A}_k f \|_{\infty} &\leqslant&  \left (1 + \frac{c_*}{k^2} \right )
\| f \|_\infty \\
\ \ (b)
\text {~For} ~~r \in \left [0, \frac{1}{2} \right ], 
~~~~\left\| \frac{d}{dr}
[\mathcal{A}_k f] (r) \right\|_{\infty} &\leqslant&  c_* k
\| f \|_\infty
\label{l4.0-0}
\end{eqnarray}
\end{Lemma}

\begin{proof}
Consider the expression for $\mathcal{A}_k f$ from (\ref{S7.9}).
We break up the integral into 
$ \int_{r}^{r+\delta}$ and $\int_{r+\delta}^1$, where $\delta ={C_2} k^{-1}
\log k $, with $C_2$ large enough so that 
$$(1-\delta_1)^{2 k-2+\tau} \leqslant  \frac{1}{k^{l/2+7/2} };\ \ \ 
(1-\delta_1) :=\frac{ \mathfrak{s} (r+\delta)}{\mathfrak{s} (r) }$$
{From}  (\ref{S7.3}) and  Lemma 
\ref{l2.5},
part (2), transforming the integration variable to $t$, it follows that
\begin{equation}
  \label{eq:451}
  \Big\lvert \int_{r+\delta}^1 \Omega (s) \frac{m_{k}(s)}{m_k (r)} G_k (r, s)
f(s) ds \Big\rvert
\leqslant  \frac{c_*}{k^2} \| f \|_\infty 
\end{equation}
In $\int_r^{r+\delta}$ (we replace upper limit $r+\delta$ by 1 if
$r+\delta > 1$).  Since $\delta_1 = O\left (k^{-1} \log k \right )$
and $t \in (1-\delta_1, 1)$ then $T_k (\mathfrak{s}, t)\ge 0$ and $G_k
(r, s) \geqslant 0$ for $r\in[k^{-2/3},1]$ Therefore,
\begin{equation}
\label{S7.12}
 \| \mathcal{A}_k f \|_\infty
\leqslant  \|f \|_\infty \left \{
\left [ \int_{r}^{r+\delta} \frac{\Omega(s) m_{k-1} (s)}{ m_k (r)}
G_k (r, s) \right ]  + \frac{c_*}{k^2} \right \}
\end{equation}
{From} (\ref{eq:451}) we get
we have $$\Big\lvert \int_{r+\delta}^1
\Omega(s) \frac{m_{k-1} (s)}{m_k (r)} G_k (r, s) ds \Big\rvert \leqslant 
\frac{c_*}{k^2} $$ Hence
\begin{equation}
\label{S7.13}
\int_{r}^{r+\delta}
\Omega(s) \frac{m_{k-1} (s)}{m_k (r)} G_k (r, s)  ds
= \int_{r}^1 \Omega(s) \frac{m_{k-1} (s)}{m_k (r)} G_k (r, s)  ds +
O \Big(\frac{1}{k^2}\Big)  
\end{equation}
Using Lemma 
\ref{l4.0.0.0},
(\ref{8180}) (a) follows.
For (b) we write
\begin{multline}
\label{15l.1}
\frac{d}{dr} \int_r^1 \Omega (s) \frac{m_{k-1} (s)}{m_k (r)}
G_k (r, s) f(s) ds
\\= 
\int_r^1 \Omega (s) 
\frac{\partial}{\partial r} \left ( G_k (r, s)
\frac{m_{k-1} (s)}{m_k (r)} \right ) f(s) ds \\
\end{multline}
 By Lemma \ref{l4.0.0.0}, the quantity above
is bounded by $c_* k \|f
\|_\infty $.
\end{proof}

\begin{Lemma}
\label{boundsBk}
For any $f \in \mathcal{L}_{\infty} [0, 1]$,
$$\| \mathcal{H}_k f \|_\infty \leqslant  \frac{c_*}{k^2} \| f \|_\infty
$$
$$ \Big\| \frac{d}{dr} [\mathcal{H}_k f] (r) \|_\infty
\leqslant  \frac{c_*}{k^2} \| f \|_\infty
$$
\end{Lemma}

\begin{proof}
  As before, we choose $\delta = {C_2}k^{-1} \log k$ large $C_2$ independent
  of $k$.  Using Lemma \ref{l2.5}, it follows that
  \begin{multline}
    \Big\lvert \int_{r+\delta}^1
\frac{\Omega (s) m_{k+1} (s)}{ m_{k (r)} } G_{k} (r, s) f(s) ds \Big\rvert
\leqslant  c_*  (1-\delta_1)^{2k + 2} k^{l/2 -5/2} \| f \|_\infty 
\leqslant  \frac{c_*}{k^4} \| f \|_\infty
  \end{multline}
\begin{multline}
  \Big\lvert \int_{r+\delta}^1 \frac{\partial}{\partial r} \left \{
\frac{\Omega (s) m_{k+1} (s)}{ m_{k (r)} } G_{k} (r, s) \right \}
f(s) ds \Big\rvert
\\ \leqslant  c_*  (1-\delta_1)^{2k + 2} k^{l/2 -3/2} \| f \|_\infty 
\leqslant  \frac{c_*}{k^3} \|f \|_\infty
\end{multline}
Now, Lemma \ref{l2.5}
implies
\begin{multline}
  \Big\lvert \int_{r}^{r+\delta}
\frac{\Omega (s) m_{k+1} (s)}{ m_{k (r)} } G_{k} (r, s) f(s) ds \Big\rvert
\leqslant  
 c_* \frac{\| f \|_\infty}{k^2} \int_0^1 t^{2k+2+\tau} dt
\leqslant  \frac{c_* \| f \|_\infty}{k^3}
\end{multline}
\begin{multline}
  \Big\lvert \int_{r}^{r+\delta}\!\!\! \frac{\partial}{\partial r} \left \{
\frac{\Omega (s) m_{k+1} (s)}{ m_{k (r)} } G_{k} (r, s) \right \} f(s) ds
\Big\rvert
\leqslant  \frac{c_* \| f \|_\infty}{k} \int_0^1 t^{2k+2+\tau} dt
\leqslant  \frac{c_* \| f \|_\infty}{k^2}
\end{multline}
\end{proof}

\begin{Lemma}
\label{lbounded}
There exist $k_0$ and 
$c_*$, independent of $k$, so that for $ k > k_0$, over the $r$-interval
$(0, 1)$,
\begin{equation}
\| h_k \|_{\infty} < c_*
\end{equation}
\end{Lemma}

\begin{proof}
First we note that for $k_0$ sufficiently large,
$\| h_{k_0} \|_{\infty}$ exists
since $g_{k_0} $ is continuous for $r \in [0, 1]$ and 
the expression for $m_k$ in (\ref{S7.3}) shows that 
$1/m_{k_0}$ is bounded as well for
sufficiently large $k_0$ since  $K_{l+1/2}$ has no zeros 
in the region of interest. 
Define $r_k = \mathcal{H}_k h_{k+1}$.
Note that
\begin{equation}
\label{S7.17}
h_{k} = \mathcal{A}_k \left ( \mathcal{A}_{k-1} h_{k-2} + r_{k-1} \right )
+ r_k
\end{equation}
In  $k-k_0$ inductive steps we get
\begin{multline}
\label{S7.18}
h_{k} = \mathcal{A}_k \mathcal{A}_{k-1} ..\mathcal{A}_{k_0+1}
h_{k_0} + \mathcal{H}_k h_{k+1} +\!\!\!\!
\sum_{m=1}^{k-k_0-1}
\left ( \prod_{j=1}^m \mathcal{A}_{k-j+1} \right )
\mathcal{H}_{k-m} h_{k-m+1}
\end{multline}
We write this abstractly as
\begin{equation}\label{eqhh}
\mathfrak{h} = \mathfrak{h}^0 + \mathfrak{N} \mathfrak{h}
\end{equation}
where
\begin{multline}
  \mathfrak{h}^0_k = \mathcal{A}_k \mathcal{A}_{k-1} ..\mathcal{A}_{k_0+1}
h_{k_0};\\ \left [ \mathfrak{N} \mathfrak{h} \right ]_{k} = \mathcal{H}_k
h_{k+1} +\sum_{m=1}^{k-k_0-1} \left ( \prod_{j=1}^m \mathcal{A}_{k-j+1} \right
) \mathcal{H}_{k-m} h_{k-m+1} 
\end{multline}
and $\mathfrak{N}$ is defined on the space $\mathcal{S}$ of sequences
$ \mathfrak{h} = \left \{ h_k \right \}_{k=k_0+1}^\infty$ in the norm
\begin{equation}
\label{S7.19}
\| \mathfrak{h} \| = \sup_{k\geqslant k_0+1} \|h_k \|_\infty ,
\end{equation}
Lemmas \ref{l4} and \ref{boundsBk} imply 
\begin{multline}
  | [\mathfrak{N} \mathfrak{h}]_k | 
\leqslant \|\mathfrak{h} \|_\infty \left ( \frac{c_*}{k_0^2}  + c_* \!\!\!\!\sum_{m=1}^{k-k_0-1}
\left \{ \prod_{j=1}^m \left [ 1 + \frac{c_*}{(k-j+1)^2} \right ] \right \}
\frac{1}{(k-m)^2} \right )  
< \nu \| \mathfrak{h} \|_\infty
\end{multline}
where, if $k_0$ is large $\nu < 1$ is independent of
$k$. Thus, $\mathfrak{N}$ is contractive and there is a unique
solution of (\ref{eqhh}) in $\mathcal{S}$.
\end{proof}

\begin{Lemma}
\label{lgkprime}  For any $r \in \left [ 0, \frac{1}{2} \right ]$ 
and for large enough $k$ we have
$\|  \frac{d}{dr} h_k  \|_\infty
\leqslant  c_* k $.
\end{Lemma}
\begin{proof}
Since by Lemma \ref{lbounded} $h_k$ is bounded, Lemmas \ref{boundsBk} and \ref{l4} imply
$$ | h_k^\prime (r) | \leqslant  |
\frac{d}{dr} [\mathcal{A}_k h_{k-1} ](r) |
+
|\frac{d}{dr} [\mathcal{H}_k h_{k+1} ](r) |
\leqslant  c_* k $$

\end{proof}

\begin{Lemma}
\label{l3.1}
For all $k \geqslant 1$,
$h_k (1) = 1$.
\end{Lemma}

\begin{proof}
In case ${\bf (i)}$, a 
simple computation  shows that
$$
\frac{\partial^{2k}g_{n_0-k}}{\partial \mathfrak{s}^{2k}} \rvert_{\mathfrak{s} =0}  = i^k h_k
(1) ;\ \ \ (g_{n_0-k}: = i^k m_k h_k )$$ (By the differential equation for $h_k$, all derivatives exist.) 
Lemma \ref{l0.1} with $j=2k$ gives
$$i^k = \frac{\partial^{2k}}{\partial \mathfrak{s}^{2k}} \rvert_{\mathfrak{s} =0} g_{n_0-k} =
i^k h_k (1) $$
implying the result in case ${\bf (i)}$. \index{bfi@${\bf i}$} In
case ${\bf (ii)}$, \index{bfii@${\bf ii}$} using Lemma~\ref{l1}, a similar computation shows that
$$ i^k =
\frac{\partial^{2k+1}}{\partial \mathfrak{s}^{2k+1}} \rvert_{\mathfrak{s} =0} g_{n_0-k} =
i^k h_k (1) \ \ (g_{n_0-k}: = i^k m_k h_k )$$
\end{proof}

\begin{Definition}
Let
\begin{equation}
\label{7.27.2}
{\hat T}_k (\mathfrak{s}, s) = s^{-2k+1-\tau}
\int_0^s t^{2k-2+\tau} \mathfrak{s} \frac{\partial}{\partial \mathfrak{s}} T_k (\mathfrak{s}, t) dt,
\end{equation}
where $T_k (\mathfrak{s}, t)$ is defined in (\ref{eqTkdef}).
\end{Definition}

\begin{Lemma}
\label{lhatTk} Let $\delta = k^{-1} \log k$ and $S_k (\mathfrak{s}): = \frac{\partial}{\partial \mathfrak{s}} \int_0^1 t^{2k-2} T_k (\mathfrak{s}, t) dt$. If $C_2$ is large enough,  $s \in (0, \delta)$ and $r (\mathfrak{s}) \geqslant k^{-1}{C_2}$ we have
\begin{multline}
{\hat T}_k (\mathfrak{s}, s) = \mathfrak{s} S_k (\mathfrak{s}) -
\frac{\mathfrak{s} f_1^\prime (\mathfrak{s})}{12} (1-s)^3
+ \frac{\mathfrak{s} f_3 (\mathfrak{s})}{3 k r^3} (1-s)^3  \\
+
O \left (\frac{(1-s)^4}{k r^4}, \frac{(1-s)^3}{k^2 r^4},
\frac{(1-s)^2}{k^3 r^4},
\frac{(1-s)}{k^4 r^3} , \frac{(1-s)^3}{k r^2},
\frac{(1-s)^2}{k^2 r^2} \right )
\end{multline}
\end{Lemma}

\begin{proof} This simply follows by integrating (\ref{S7.25}) from  $t=1$ to $s$    of $T_k$  and the fact that ${\hat T}_k
  (\mathfrak{s}, 1) = \mathfrak{s} S_k (\mathfrak{s})$.
\end{proof}
\subsection{Proof of Lemma \ref{lgkk0}}\label{plgkk0}
First choose $\epsilon_1 > 0$.
{From} Lemma \ref{boundsBk}, it follows that
$$\left\|\frac{d}{d \mathfrak{s}} [\mathcal{H}_k h_{k+1} ] \right\|_\infty \leqslant  \frac{c_*}{k^2} \|
h_{k+1} \|_\infty \leqslant  \frac{c_*}{k^2} $$
where we applied Lemma ~\ref{lbounded}.
Further, we note that
\begin{multline}
  \frac{1}{(2k+\tau) (2k+\tau -1)} \frac{d}{d\mathfrak{s}} \mathcal{A}_k h_{k-1} (\mathfrak{s})
\\= \int_0^1 t^{2k+\tau-2}
\frac{\partial T_k}{\partial \mathfrak{s}} (\mathfrak{s}, t) h_{k-1} (\mathfrak{s} t) dt
+\int_0^1 t^{2 k+\tau-1} T_k (\mathfrak{s}, t) h_{k-1}^\prime (\mathfrak{s} t) dt
\end{multline}
We have
\begin{multline}
  \int_0^1 t^{2k+\tau-2} \frac{\partial T_k}{\partial \mathfrak{s}} (\mathfrak{s}, t) h_{k-1}
  (\mathfrak{s} t) dt = h_{k-1} (\mathfrak{s}) S_k (\mathfrak{s})\\ -\int_0^1
  dt \,\,t^{2k+\tau-2} \mathfrak{s}
  \frac{\partial T_k}{\partial \mathfrak{s}} (\mathfrak{s}, t) \int_t^1 h_{k-1}^\prime
  (\mathfrak{s} s) ds = h_{k-1} (\mathfrak{s}) S_k (\mathfrak{s}) \\- \int_0^1 h_{k-1}^\prime (\mathfrak{s} s) \left
    [\int_0^s t^{2 k+\tau -2} \mathfrak{s} \frac{\partial T_k}{\partial \mathfrak{s}} (\mathfrak{s}, t)
    dt \right ] ds
  \\
  = h_{k-1} (\mathfrak{s}) S_k (\mathfrak{s}) -\!\!\! \int_0^1\!\!\! h_{k-1}^\prime (\mathfrak{s} s)
  s^{2k-1+\tau} {\hat T}_k (\mathfrak{s}, s) ds  = (2k+\tau-1) S_k (\mathfrak{s}) \\
  \times \int_0^1 s^{2 k -2+\tau} h_{k-1} (\mathfrak{s} s) ds
  -\int_0^1 s^{2k -1+\tau} [{\hat T}_k (\mathfrak{s}, s)-{\hat T}_k (\mathfrak{s}, 1)]
  h_{k-1}^\prime (\mathfrak{s} s) ds
\end{multline}
Therefore,
\begin{multline}
\frac{\frac{d}{d\mathfrak{s}} \mathcal{A}_k [h_{k-1} ] (\mathfrak{s} )}{(2k+\tau) (2k+\tau-1)} 
= \int_0^1 [T_k (\mathfrak{s}, s) - {\hat T}_k (\mathfrak{s}, s)+\mathfrak{s} S_k (\mathfrak{s})]
s^{2k +\tau-1}
\\\times h_{k-1}^\prime (\mathfrak{s}  s) ds 
+ (2k+\tau -1) S_k (\mathfrak{s}) 
\int_0^1 s^{2 k +\tau-2} h_{k-1} (\mathfrak{s} s) ds 
\end{multline}
We note that
$$(2k+\tau) (2k+\tau-1) S_k (\mathfrak{s}) 
= \frac{\partial}{\partial \mathfrak{s}} \left [
\mathcal{A}_k [1] (\mathfrak{s}) \right ] = 
O \left(\frac{1}{k^3 \epsilon_1^2}, \frac{1}{k^2}\right)
$$
and that $ (2k+\tau -1) \int_0^1 s^{2k+\tau-2} h_{k-1} (\mathfrak{s}
s) ds$ has a bound independent of $k$.  Combining (\ref{S7.25}) with
Lemma \ref{lhatTk}, if $k$ is large so that $k \epsilon_1 $ is large, then
\begin{multline}
\label{TkhatTkexp}
T_k (\mathfrak{s}, s) - [{\hat T}_k (\mathfrak{s}, s) -\mathfrak{s} S_k (\mathfrak{s})]
= (1-s) +
\left (-\frac{k f_1}{4} + \frac{f_2}{r^2} \right )
 \\\times \Big[-\frac{(1-s)^2}{k} + \frac{2}{3} (1-s)^3\Big]
- \mathfrak{s}
\left ( -\frac{f_1^\prime}{12} + \frac{f_3}{3 kr^3} \right  ) (1-s)^3 \\
+
O \Big (\frac{(1-s)^4}{k r^4}, \frac{(1-s)^3}{k^2 r^4},
\frac{(1-s)^2}{k^3 r^4},
\frac{(1-s)}{k^4 r^3} , \frac{(1-s)^3}{k r^2},
\frac{(1-s)^2}{k^2 r^2},\\
\frac{(1-s)^4}{r^3}, \frac{(1-s)^3}{k r^3},
\frac{(1-s)^3}{r}, \frac{(1-s)^2}{k r} \Big )
\end{multline}
{From} (\ref{TkhatTkexp}), it is clear that $T_k (\mathfrak{s}, s) -
{\hat T}_k (\mathfrak{s}, s) + \mathfrak{s} S_k (\mathfrak{s}) > 0$ if
$s \in (1-\delta, 1)$ and $k \epsilon_1 $ is sufficiently large. Now,
$\mathfrak{s} f_3/ \left ( 3k r^3 \right ) (1-s)^3 >0$ 
exceeds any term  following it in
(\ref{TkhatTkexp}),  except possibly when $1-r$, {\it i.e.} $\mathfrak{s}$ 
is small. Thus, if we define
\begin{equation}
M_k = \sup_{r(\mathfrak{s}) \in [\epsilon_1, 1]} |h_{k}^\prime (\mathfrak{s}) |
\end{equation} we get
\begin{multline}
|h_k^\prime (\mathfrak{s}) | \leqslant 
(2k+\tau) (2k+\tau -1) M_{k-1} \left \{ \int_{1-\delta_1}^1 s^{2k+\tau -1}
\left [ (1-s)
+ \left (-\frac{k f_1}{4} + \frac{f_2}{r^2} \right ) \right . 
\right .  \\\times
\left . \left . [-\frac{1}{k} (1-s)^2 + \frac{2}{3} (1-s)^3] 
+ \mathfrak{s}
\frac{f_1^\prime}{12} (1-s)^3 \right ] ds \right \}  + \frac{c_*}{k^2}
+ \frac{c_*}{k^3 \epsilon_1^2}
\end{multline}
When $(1-r)$ (and thus $\mathfrak{s}$)  is small, we can replace
the term $\mathfrak{s} f_1^\prime/(12) (1-s)^3 $ on
the right side of the above equation simply by  $(1-s)^3$, which is
clearly bigger.
{From} the fact that 
$\int_{1-\delta}^1 s^{2k-1}
[-k^{-1} (1-s)^2 + (2/3) (1-s)^3] ds = O(k^{-5} )$, 
it follows that
\begin{equation}
M_k \leqslant  M_{k-1} \left ( \frac{2k-1+\tau}{2k+1+\tau} + \frac{c_*}{k^2} + 
\frac{c_*}{k^3 \epsilon_1^2}
\right ) + \frac{c_*}{k^2} + \frac{c_*}{k^3 \epsilon_1^2}
\end{equation}
Let $C_3$ be large enough and define $k_0 (\epsilon_1) =
C_3/\epsilon_1,$ so that for $k \geqslant k_0$ we have 
$$  \left ( \frac{2 k +\tau-1}{2 k +\tau+1} + \frac{c_*}{\epsilon_1^2 k^3} + \frac{c_*}{k^2}
\right )
\leqslant  \left ( \frac{k-1}{k} \right )^{1/2} $$
Then for $k \geqslant k_0$, 
\begin{equation}
\label{S7.35}
M_{k} \leqslant  \left ( \frac{k-1}{k} \right )^{1/2} M_{k-1} +
\frac{c_*}{k^2}  + \frac{c_*}{k^3 \epsilon_1^2}
\end{equation}
implying
\begin{multline}
\label{S7.36}
M_k \leqslant  \left ( \frac{k_0}{k} \right )^{1/2} M_{k_0}
+ \frac{c_*}{k^{1/2}} \sum_{j=k_0}^k
\frac{1}{j^{3/2}} + \frac{c_*}{k^1/2} \sum_{j=k_0}
\frac{1}{j^{5/2} \epsilon_1^2}
\\\leqslant  c_* \frac{k_0^{3/2}}{k^{1/2}} + 
\frac{c_*}{k^{1/2} k_0^{1/2}}
+ \frac{c_*}{k^{1/2} k_0^{3/2} \epsilon_1^2}
\end{multline}
The result follows from the definition of $M_k$ and noting that last two 
terms in (\ref{S7.36}) are $O( c_* {k_0^{3/2}}{k^{-1/2}} )$

\subsection{Proof of Lemma \ref{gkl}}\label{pgkl}
{From} Lemma \ref{lgkk0} and the definition of $k_0$, it follows that
$$ |h_k^\prime (\epsilon_1) | \leqslant   
\frac{C_4 C_3^{3/2}}{k^{1/2} \epsilon_1^{3/2}}
$$
for $k \geqslant {C_3}{\epsilon_1}^{-1} = k_0$.
Using $h_k (1) =1$, it follows that for $ k \ge C_3/r $, 
$$
| h_k (r) - 1 | \leqslant  \int_{r}^1 |h_k^\prime (r')| dr' \leqslant  \frac{C_4
  C_3^{3/2}}{\frac{1}{2} (kr)^{1/2}} $$
Additionally,  if \footnote{
It is to be noted that for small enough $\epsilon$ the inequality $\alpha k r \ge L_\epsilon$ always implies
 $k \ge {C_3}/{r}$.} 
$\displaystyle \alpha k r \geqslant
\left ( \frac{C_4 C_3^{3/2}\alpha^{1/2}}{\frac{1}{2} \epsilon } \right
)^{2}=L_\epsilon \ \ \ \text{then\ \ \ } | h_k (r) - 1 | \leqslant  \epsilon $.

\subsection{Proof of Lemma \ref{lemarzela}}\label{Sarz}
For $\zeta \in [0, L_\epsilon] $, using the {\it a priori} boundedness of
$h_k$ in $k$ and Lemma \ref{lgkprime}, we note
that both ${\tilde h}_k(\zeta) := h_k (r(\zeta))$ and 
$(\tilde{h}_k)_\zeta$ are bounded independently of $k$.  Hence the sequence $
\{ {\tilde h}_k\}_{k\geqslant 2}$ is bounded and equicontinuous.  
By Ascoli-Arzel\`a's theorem, there exists a subsequence 
${\tilde h}_{k_j} (\zeta)$
converging to a continuous function ${\tilde h} $. The first part of
the result is proved.
 We first prove that $|{\tilde h} (\zeta) - 1| \le 4 \epsilon$.
Now, from 
Lemma \ref{gkl}
\begin{equation}
\label{tildehk1}
| {\tilde h}_k (\zeta) - 1 | \leqslant   \epsilon ~~ ~{\rm for}~ \zeta \in 
[L_\epsilon, 
\alpha k] 
~~{\rm for ~sufficiently~large}~k
\end{equation} 
Let ${\tilde h}_{k, j}$ be a subsequence that converges to ${\tilde h} $
for $\zeta \in [0, L_\epsilon]$.  Let $\zeta_m$, $\zeta_M$ be a minimum,
and a maximum point of $\tilde{h}$ on 
$ [0, L_\epsilon]$ and the corresponding minimum and maximum values are
denoted by $m$ and $M$ respectively.  Continuity at the endpoint $\zeta=L_\epsilon$ implies that 
$M \ge 1-\epsilon$, $m \le 1+\epsilon$.  If both $M-1-\epsilon < 0$
and $m-1+\epsilon > 0$, there is nothing to prove because in that case
it is clear that $|{\tilde h} (\zeta) - 1| \le 2 \epsilon $. Now, consider
the possibility that
({\bf i}): $M > 1 +\epsilon$.
In a similar manner, we will also consider the possibility
({\bf ii}): $m < 1-\epsilon $. 
Consider ({\bf i}) first.
Since at the end point of the interval,
${\tilde h} (L_\epsilon) < 1+\epsilon$, from continuity
there exists an interval $[a,b] \subset [\zeta_M, L_\epsilon] $ 
of nonzero length for which
\begin{equation}
  \label{eq:nonum0}
  {\tilde h} (\eta) \le \frac{1}{2} (M+1+\epsilon) < M ~{\rm for}~
\eta \in [a, b]
\end{equation}
For some ${\hat L} > L_\epsilon$, independent of $k$ (to be 
determined shortly), we write
\begin{multline}
\label{nn12}
\left [ \mathcal{A}_k^0 f \right ] (\zeta)=
\left ( \int_\zeta^{\hat L} + \int_{\hat L}^{k \alpha \epsilon_1 } 
\right )K(\zeta,\eta)f(\eta (1-k^{-1}) )  d\eta\\ {\text{with\ \ }} K(\zeta,\eta):= e^{-Q(\eta) + Q(\zeta)}
\left ( 1 + \frac{a_1}{k} \right )
\frac{H(\eta (1-k^{-1})}{H(\zeta)}
\mathcal{G} (\zeta, \eta) d\eta \\
=: [\mathcal{A}_k^{00} f](\zeta) +
[\mathcal{A}_k^{01} f ] (\zeta)
\end{multline}
 For fixed $\zeta$ and $\eta$ we have
\begin{equation}
\lim_{k \rightarrow \infty} K(\zeta, \eta) = K_0 (\zeta, \eta)= e^{-\eta+\zeta} \frac{H_0 (\eta)}{H_0 (\zeta)} \mathcal{G}_0 
(\zeta, \eta) 
\end{equation} 
On our interval we have $\eta \geqslant \zeta$. Thus $\mathcal{G}_0 \geqslant 0$
(see
(\ref{eq:nonum00})); 
$\mathcal{G}_0$ can vanish only if $\eta=\zeta$.
Furthermore, by (\ref{eq:nonum0}) we have
$\zeta_M\not\in[a,b]$. We can then define
$$
J = \frac{3 \sup_{[0, L_\epsilon]} |{\tilde h} |}{(b-a) K_m}, 
~{\rm where}~ K_m = \min_{\eta \in [a,b]} K_0
(\zeta_M, \eta) > 0$$
 Note that $Q(\eta)\sim\eta$ for large
  $k$ and, aside from the exponential term, $K$ is
  algebraically bounded.  
We can thus choose ${\hat L} > L_\epsilon$ large
enough independently of $k$, so that

\begin{equation}
\label{tildeLcont}
| [\mathcal{A}_k^{01} f] (\zeta) | 
\leqslant  {\epsilon}{J}^{-1} \| f \|_{\infty, [L_\epsilon, k\alpha
\epsilon_1 ]} 
\end{equation}

There is a subsequence of ${\tilde h}_{k_j}$ that converges uniformly
on $ \in [0, {\hat L}]$; for simplicity, we will use the same notation
${\tilde h}_{k,j}$ for the subsequence.  
It is clear that the limit is ${\tilde h} (\zeta)$ 
if  $\zeta \in [0, L_\epsilon]$.
We keep the notation ${\tilde h}$ for the limit
on $[0, {\hat L}]$.  We note that
(\ref{tildehk1}) implies
\begin{equation}
\label{tildeh1}
| {\tilde h} (\zeta) - 1 | \leqslant  \epsilon ~~ ~{\rm for}~ \zeta \in 
[L_\epsilon, {\hat L}]
\end{equation} 
Now choose a small $\epsilon_2 > 0$.
It is clear that in the interval $[L_\epsilon, {\hat L}]$, 
${\tilde h} (\zeta) \le 1+\epsilon
< M $. 
For sufficiently large $k_j$, using continuity of ${\tilde h} (\zeta)$,
we have
\begin{multline*}
[\mathcal{A}_{k,j}^{00} {\tilde h} (\zeta_M) ]  
\le \int_{\eta \in [\zeta_M, {\hat L}]\setminus [a,b]} K(\zeta_M, \eta)
{\tilde h} (\eta) d\eta + 
\int_{a}^b K(\zeta, \eta) {\tilde h} (\eta) d
\eta + M \epsilon_2 \\
\leqslant    
 M \int_{\eta \in [\zeta_M, {\hat L}] \setminus [a,b]} K(\zeta_M, \eta) d\eta
+ \frac{1}{2} (M+1+\epsilon) \int_{a}^b K(\zeta_M, \eta) d
\eta + \epsilon_2 M \\
= M \int_{0}^{\hat L} 
K(\zeta_M, \eta) d\eta
- \frac{1}{2} (M-1-\epsilon) \int_{a}^b K(\zeta_M, \eta) d\eta 
+ M \epsilon_2
\\ \leqslant  M \mathcal{A}_{k_j}^{00} [1](\zeta_M) 
-\frac{(b-a)}{3} (M-1-\epsilon) K_m + M \epsilon_2 
\end{multline*}   
Since $\mathcal{A}_{k_j} [1] = \mathcal{A}_{k_j}^{00}[1] +
\mathcal{A}_{k_j}^{01}[1] + \mathcal{A}_{k_j}^1 [1]$ (see (\ref{n12})
and (\ref{nn12})) Lemmas \ref{l4.0.0.0}, \ref{14.0.0} and
(\ref{tildeLcont}) imply that for large $k_j$ we have
$$ [\mathcal{A}_{k_j}^{00} [1] ](\zeta_M) \leqslant  1 + \frac{\epsilon}{J}
+\epsilon_2 $$     
Hence, for large $k_j$ we have
\begin{equation}
\label{Akj00ub}
[\mathcal{A}_{k_j}^{00} {\tilde h}] (\zeta_M)  
\leqslant  M \left (1+\frac{\epsilon}{J} 
+ 2 \epsilon_2 \right ) - \frac{K_m}{3} (M-1-\epsilon) (b-a) 
\end{equation}
Now, there exists $N$ so that if $j \geqslant N$,
$ \| {\tilde h}_{k_j} - {\tilde h} \|_{\infty, [0,{\hat L}]} < 
\epsilon_2 $
and
$ \Lambda_j = \mathcal{A}_{k_{j+1}} ... \mathcal{A}_{k_j + 1}$
satisfies
$$ \| \Lambda_j - I \|_\infty \leqslant
\epsilon_2 $$
while
$$ r_{j+1} :=
B_{k_{j+1}} +
\sum_{m=1}^{k_{j+1}-k_j-1} \prod_{l=1}^m \mathcal{A}_{k_{j+1}-l+1}
B_{k_{j+1}-m}, $$
where $B_l = \mathcal{H}_l h_{l+1}$, 
satisfies the estimate
$$ | r_{j+1}| < \epsilon_2$$
Therefore, from 
$$ {\tilde h}_{k_{j+1}} = \Lambda_j \mathcal{A}_{k_j} {\tilde h}_{k_{j}}
+ r_{j+1} $$
it follows that
$$ {\tilde h}_{k_{j+1}} (\zeta_M) \geqslant {\tilde h} (\zeta_M) -\epsilon_2
= M -\epsilon_2$$
On the other hand, at $\zeta = \zeta_M $ we have
\begin{multline}
   \Lambda_j \mathcal{A}_{k_j} {\tilde h}_{k_j} + r_{j+1} \leqslant 
(1 + \epsilon_2) \left [
 M (1+ \frac{\epsilon}{J} + 2 \epsilon_2) + \epsilon_2
- \frac{K_m}{3} (M-1-\epsilon) (b-a) \right ]
+\epsilon_2 
\end{multline}
Thus,
$$ M - \epsilon_2 \leqslant
(1 + \epsilon_2) \left [
 M (1+ \frac{\epsilon}{J}+ 2 \epsilon_2) + \epsilon_2
- \frac{K_m}{3} (M-1-\epsilon) (b-a) \right ]
+\epsilon_2 $$
This is true for any $\epsilon_2$, hence as
$\epsilon_2 \downarrow 0$. 
Thus, 
$$ M \leqslant
\left [
 M \left (1+ \frac{\epsilon}{J} \right )
- \frac{K_m}{3} (M-1-\epsilon) (b-a) \right ]
$$

However, from the definition of $J$, 
this implies $M-1-\epsilon \leqslant  \epsilon$.
We note that for ({\bf ii}), we repeat the above argument for  
$-{\tilde h}$,
which has a maximum at $\zeta_m$, to conclude that either
$(-m) - (-1+\epsilon) \le 0 $ or 
$ (-m) - (-1+\epsilon) = 1-\epsilon-m \leqslant \epsilon$. 
Therefore,
$$ 1-2 \epsilon \leqslant m \leqslant M \leqslant 1+ 2 \epsilon $$ 
implying that $|{\tilde h} - 1| \le 4 \epsilon $.

\section{Appendix}
\subsection{Short proof of the regularity of the unitary propagator}\label{lapts1}
\begin{Theorem}
  Assume that $H_1=H+V(x,t)$, where $H$ is time independent and
  self-adjoint, and $V(\cdot,t)$ is in $L^\infty(\RR^n)$ for every $t$
  and is differentiable in time, with integrable derivative. Consider
  the Schr\"odinger problem
\begin{equation}
  \label{eq:eqin}
  i\psi_t=H_1\psi;\ \ \psi(x,0)\in D(H)
\end{equation}
Then there exists a strongly differentiable unitary propagator on $L^2(\RR^n)$
$U(t)$ so that $\psi(x,t)=U(t)\psi_0\in D(H)$ for all $t$ and $\psi(x,t)$
solves(\ref{eq:eqin}).

\end{Theorem}
\begin{proof}
We note that it enough to 
prove this property on a finite interval $[0,\epsilon]$, since the problem
can be restarted at $t=\epsilon$. 
Let $y=\psi-e^{-t}\psi_0$. Then $y$ satisfies the inhomogeneous Schr\"odinger
equation
\begin{equation}
  \label{eq:eqy1}
  iy_t=y_0e^{-t}+Hy+Vy;\  y_0:=i\psi_0+H\psi_0+V\psi_0,\ \ y(0)=0
\end{equation}

 We transform this equation into an integral equation, {\em formally for now}. Straightforward calculations show that
\begin{equation}
  \label{eq:eqtr1}
  i(e^{iHt}y)_t=e^{iHt}e^{-t}y_0+e^{iHt}Vy
\end{equation}
or (still formally)
\begin{multline}
  \label{eq:eq4}
 i e^{iHt} y 
=\left(\int_0^t e^{(iH-1)s}ds\right)y_0+\int_0^t e^{iHs}V(s)y(s)ds\\=(iH-1)^{-1}(e^{iHt-t}-1)y_0+\int_0^t e^{iHs}V(s)y(s)ds
\end{multline}
or, equivalently,
\begin{equation}
  \label{eq:eqtr2}
iy 
=(iH-1)^{-1}(e^{-t}-e^{-iHt})y_0+e^{-iHt}\int_0^t e^{iHs}V(s)y(s)ds
\end{equation}
It is clear that (\ref{eq:eqtr2}) is contractive in the norm
$\sup_{t\in [0,\epsilon]}\|\cdot\|_{L^2(\RR^3)}$ for small $\epsilon$, 
and
has a unique solution. Clearly, the first term on the right
side of (\ref{eq:eqtr2}) is differentiable in time and the derivative
is continuous since $e^{-iHt}$ is; let $u_0$
denote this derivative. 

We now write a {\em formal}  equation for $u=y_t$. We have
\begin{equation}
  \label{eq:equ5}
  iu=u_0+\int_0^te^{-iHs}V'(t-s)\left(\int_0^{t-s}u(s')ds'\right) ds+\int_0^te^{-iHs}V(t-s)u(t-s)ds
\end{equation}
This equation is also contractive, and has a unique solution, in the same 
space. Thus both sides of (\ref{eq:equ5}) are integrable in time. By integration
and appropriate changes of variables and order of integration, we see that $\int_0^t u(s)ds$ satisfies
the same equation as $y$, which has a unique solution. Thus $y=\int_0^t u(s)ds$ is strongly differentiable. 
Since both $y$ and $e^{iHt}y$ are strongly differentiable (the latter
by inspection from (\ref{eq:eq4})), $y\in D(H)$ for all $t$ and is strongly
differentiable. It is clear  that $\psi\in D(H)$ and
easy to check that it is differentiable and satisfies (\ref{eq:eqin}).
\end{proof}
\subsection{Laplace transform of the Schr\"odinger equation}\label{lapts}
We look more generally at equations of the form
\begin{equation}
  \label{eq:eqs1}
  i\psi_t=H\psi+V(t,x)\psi
\end{equation}
where $H$ is self-adjoint and time independent, and $V(x,t)$ 
is bounded 
on $\RR^3$ and differentiable  and bounded in $t$, 
and $\psi(x,0)\in D(H)$. The conditions
on $V$ can be relaxed.  (For the purpose of this paper, $H$ would be taken to be $H_C$.)
\begin{Proposition}\label{LT1}
  Under the assumptions above, 
the  Laplace transform ${\hat \psi} (p, \cdot)$ of $\psi (t, \cdot)$ exists
for $\Re\, p > 0$; 
it is in $D(H)$ and satisfies
\begin{equation}
  \label{eq:eqpsi3}
  (p+iH)\hat{\psi}=\psi_0-i\widehat{V\psi}
\end{equation}
\end{Proposition}
\begin{proof}
  We take the unitary propagator of the time-independent problem,
  $U=e^{-iHt}$ and apply $U^*(t)=U^{-1} (t)$ to both sides of
  (\ref{eq:eqs1}).  Since (cf. \S\ref{setting}) $U^{-1}$ is strongly
  differentiable, with derivative $iU^{-1}H$, and $\psi$ is
  $t-$differentiable in $L^2$, $U^{-1}\psi$ is differentiable and we
  get
\begin{equation}
  \label{eq:eqg}
  (U^{-1}\psi)_t=iU^{-1}H\psi+U^{-1}\psi_t=-iU^{-1}V\psi
\end{equation}
 Since $U^{-1} V \psi$ is 
continuous in $t$, we can integrate both sides
and get, after multiplication by $U$ and using the fact that 
$U^{-1}(t)=U(-t)$, 
\begin{multline}
  \label{eq:eq3}
   \psi=U\psi_0-iU\int_0^tU^{-1}V\psi(s)ds=U\psi_0-i\int_0^tU(t-s)(V\psi)(s)ds\\
=U\psi_0 -iU*(V\psi)
\end{multline}
where $*$ is the usual  Laplace convolution.
Taking the Laplace transform (which clearly exists) in (\ref{eq:eq3}) 
and using standard functional
calculus we get
\begin{equation}
  \label{eq:eqpsi2}
  \hat{\psi}=(p+iH)^{-1}\psi_0-i(p+iH)^{-1}\widehat{V\psi}
\end{equation}
and thus $\hat\psi$ is a $D(H)$ solution of (\ref{eq:eqpsi3}). 
\end{proof}

Now, from eq. (\ref{psidecomp}), it follows that ${\hat y}$ satisfies
(\ref{eq:S-lt-y}). Furthermore, using (\ref{eq:eqpsi2}) and the 
fact that  $y^0$ and $\Omega_j$ are compactly supported, we see that $\hat{y}$ 
also satisfies
\begin{equation}\label{eq:lps}
{\hat y} (p,\cdot ) = \mathfrak{R}_0 \bchi_{\sf B}
{\hat y}^0 (p, \cdot) - \mathfrak{R}_0 
\bchi_{\sf B} 
\left [ \sum_{j \in \mathbb{Z}} \Omega_j {\hat y} ({\hat p} - i j \omega,  \cdot)
\right ] 
\end{equation}
where $\mathfrak{R}_0 = (H_C-ip)^{-1} $.

\subsection{Analyticity of $(I-\mathfrak{C_{l,m}})^{-1}$ in $X$}\label{CX}
This is standard, and can be seen directly from analytic functional
calculus. We provide a self-contained argument, for completeness.  We
write $\mathfrak{C}_{X}$ to emphasize the $X-$ dependence of
$\mathfrak{C}$, and for simplicity of notation we drop the $(l,m)$
subscript.  We have
\begin{multline}  \label{eq:eqinv2}
    (I-\mathfrak{C}_{X_1})^{-1}-(I-\mathfrak{C}_{X'})^{-1}=
(I-\mathfrak{C}_{X'})^{-1}(\mathfrak{C}_{X_1}-\mathfrak{C}_{X'})
(I-\mathfrak{C}_{X_1})^{-1}\ \ \text{and} \\
  (I-\mathfrak{C}_{X'})^{-1}\Big[I+(\mathfrak{C}_{X_1}-
\mathfrak{C}_{X'}) (I-\mathfrak{C}_{X_1})^{-1}\Big]=
(I-\mathfrak{C}_{X_1})^{-1}
\end{multline}
We fix $X_1$ and let $X'\to X_1$. Since $(I-\mathfrak{C}_{X_1})^{-1}$ 
is bounded, then
$\|(\mathfrak{C}_{X_1}-\mathfrak{C}_{X'})(I-
\mathfrak{C}_{X_1})^{-1}\|\to 0$ as $X'\to X_1$ and
\begin{equation}
  \label{eq:inv3}
  I+(\mathfrak{C}_{X_1}-\mathfrak{C}_{X'})
(I-\mathfrak{C}_{X_1})^{-1}
\end{equation}
is invertible when $X_1$ and $X'$ are close enough and
$[I+(\mathfrak{C}_{X_1}-\mathfrak{C}_{X'})(
I-\mathfrak{C}_{X_1})^{-1}]^{-1}\to I$ in operator norm as $X'\to X_1$. Thus
\begin{equation}
  \label{eq:inv5}
(I-\mathfrak{C}_{X'})^{-1}\to   (I-\mathfrak{C}_{X_1})^{-1}
\end{equation}
in operator norm, as $X_1 \to X'$. Now diferentiability in $X$ follows from 
(\ref{eq:eqinv2}).

\subsection{Coulomb Green's function representation}\label{CoulombG}
\label{Green}
The retarded Green's functions $G=G_+$ is defined as the solution of the
equation,
\begin{equation}
  \label{eq:CoulombG}
 \mathcal{A}_0G(x,x';k)=\delta(x-x')
\end{equation}
in distributions, satisfying the radiation condition
\begin{equation}
  \label{eq:condinf}
 G(x,x';k)\sim  F(\theta,\phi)e^{ikr}r^{-1-i\nu};
\text{ as }\ \ r\to\infty
\end{equation}
where
\begin{equation}
  \label{eq:defnu}
  k=\sqrt{ip} \ \ \text{($\Im\,k>0$ if $\Re\,p>0$)},\ \ \nu=\frac{b}{2k}
\end{equation}
Equivalently, $G$ is the $\RR^3\setminus\{0\}$ solution of
(\ref{eq:CoulombG}) with zero right hand side, satisfying (\ref{eq:condinf})
and $ |x-x'|G(x,x';k)\to (4\pi)^{-1}$ as $x-x'\to 0$.

\begin{Proposition}\label{valr0}
  \begin{equation}
  \label{eq:formula}
\mathfrak{R}_0 {\bchi_{\sf B}} g=\int_{\sf B} G(x,x';k)g(x')dx'
\end{equation}
\end{Proposition}
\begin{proof} The function
\begin{equation}
  \label{eq:defGG}
  f:= \int_{\sf B}G(x,x';k)g(x')dx'
\end{equation}
solves, as can be checked, the equation
\begin{equation}
  \label{eq:loceq}
 \mathcal{A}_0 f={\bchi_{\sf B}} g
\end{equation}
with the radiation condition (\ref{eq:condinf}). Such a solution is
unique since the difference of two solutions satisfies the equation $
\mathcal{A}_0 f=0$ (with the radiation condition (\ref{eq:condinf})).
Multiplying by $G(x,x';k)$, integrating over a volume and passing to
the limit where the volume approaches $\RR^3$ we see that $f\equiv 0$.
\end{proof}
Symmetries of the Coulomb potential $-b/r$ allow for a closed
form of $G$ (cf. \cite{Hostler}-- where the sign is chosen differently) in terms of Whittaker functions
$\mathfrak{W} $ and $\mathfrak{M} $,
\begin{equation}
  \label{eq:explsol}
  G(x;x';k)=
\frac{\Gamma(1-i\nu)}{4\pi ik|x-x'|}\left(\frac{\partial}{\partial
  \xi}-\frac{\partial}{\partial \eta }\right)
\mathfrak{W}_{i\nu,\frac{1}{2}}(-ik\xi)\mathfrak{M}_{i\nu,\frac{1}{2}}(-ik\eta )
\end{equation}
where $\Im\,k>0$ , $2k\nu=b$ and
\begin{equation}
\label{defxinu1}
\xi = |x| + |x'| + |x-x'| ~~~~\mbox , ~~~
\eta = |x|+|x'| - |x-x'| 
\end{equation}
The Whittaker functions are defined in terms the Kummer functions $M$ and
$U$ by the relations, see \cite{Abramowitz} Chapter 13,
\begin{multline}
  \label{eq:hyperg}
  \mathfrak{M}_{\kappa,\mu}(z)=e^{-\frac{z}{2}}z^{\frac{1}{2}+\mu}M\left(\frac{1}{2}+\mu-\kappa,1+2\mu,z\right),\ 
 -\pi<\arg z\leqslant  \pi\\
 \mathfrak{W}_{\kappa,\mu}(z)=
e^{-\frac{z}{2}}z^{\frac{1}{2}+\mu}U\left(\frac{1}{2}
+\mu-\kappa,1+2\mu,z\right),\ 
 -\pi<\arg z\leqslant  \pi
\end{multline}
The following integral representation follows from \cite{Abramowitz} Chapter
13, for the
values we are interested in, $z_1=-ik\xi$, $z_2=-ik\eta$, $a=1-i\nu$, $b=2$
(a different ``$b$''
than the one in our Coulomb potential)
\begin{equation}\label{defM}
 \mathfrak{M}_{i\nu;\frac{1}{2}}(z)
 =\frac{e^{-\frac{1}{2}z  }zJ (z)}{\Gamma(1-i\nu)\Gamma(1+i\nu)};
W_{i\nu;\frac{1}{2}}(z)
 =\frac{e^{-\frac{1}{2}z}z I (z)}{\Gamma(1-i\nu)}
 \end{equation}
where $I$ and $J$ are as defined in (\ref{defIJ}) and the expression is valid
in the regions where the integrals converge (in particular, $|\Im\,\nu|<1$).
 For other values of $\nu$ of interest, the integrals can be replaced by
 appropriate contour integrals. For instance $J$ would be replaced by
 $$\Big(1-e^{-2\pi\nu}\Big)^{-1}\oint_C e^{z t}t^{-i\nu}(1-t)^{i\nu}dt$$
 where $C$ is smooth simple curve encircling $[0,1]$, as it can be
 checked by calculating the jump across the cut of the integrand.  It
 follows from these integral representations that the Green's function
 is analytic at any (small) $p$, $\Re\,p\ne 0$.  
Substituting (\ref{defM}) into (\ref{eq:explsol}), 
we obtain (\ref{eq:exprr}).

\smallskip

\subsection{Dependence of $A$ in equation (\ref{eq:valA}) on $Z$, 
$p$}\label{Adep}

We now seek to determine the asymptotics of $A$ in (\ref{eq:valA})
in the
resolvent 
$\bchi_{\sf B} 
\mathfrak{R}_{\beta} \bchi_{\sf B}$ 
in terms of 
$\lambda = \sqrt{-i p}$ and
$Z = \exp \left [ i \pi b/ (2 \lambda )  
\right ]$ for  $X = (\sqrt{p}, Z) \in  
\overline{\mathbb{D}^+_\epsilon}\times \overline{\mathbb{D}}$ 
for
sufficiently small $\epsilon$.

Recall the expression $A$ in (\ref{eq:valA}). Note that since 
\begin{equation}
\label{a.A.1}
\alpha
= \sqrt{\lambda^2-ic} \sim 
e^{-i \pi/4} c^{1/2} \left (1 + O(\lambda^2) \right ) \equiv \alpha_0 +
\lambda^2 \alpha_1 ~+~.. 
\end{equation}
\begin{equation}
\label{a.A.2} 
\kappa_1 = \frac{b}{2 \alpha} = \frac{b e^{i \pi/4}}{\sqrt{c}} 
\left [ 1 + O(\lambda^2) \right ] \equiv \kappa_{1,0} + \lambda^2 
\kappa_{1,2} ~+~, 
\end{equation} 
each of $m_1 (a)$ and $w_1 (a)$ is analytic in $\lambda$ for small
$\lambda$, with the expansion 
\begin{equation}
\label{a.A.3}
m_1 (a) = \frac{1}{r} \mathfrak{M}_{\kappa_1, l+1/2} (2 \alpha a) 
\sim m_{1,0} (a) + \lambda m_{1,1} (a) ~+~... 
\end{equation}
\begin{equation}
\label{a.A.3.1}
w_1 (a) = \frac{1}{a} \mathfrak{W}_{\kappa_1, l+1/2} (2 \alpha a) 
\sim w_{1,0} (a) + \lambda w_{1,1} (a) ~+~... 
\end{equation}
The asymptotics in this case is also differentiable with respect to
$a$ and we get similar expressions as above for $m_1^\prime (a)$ and
$w_1^\prime (a)$. 
It follows from the expression for $f_0$ in (\ref{eq:inner}) also possesses
a regular series expansion in $\lambda$:
\begin{equation}
\label{a.A.4}
f_0 (a) = f_{0,0} (a) +\lambda f_{0,1} (a) + ...
\end{equation}
To simplify $A$ as in (\ref{eq:valA}) for small $\lambda$, we now
now consider the asymptotics of $w_2 (a)$ and $w_2^\prime (a)$ for
small $\lambda$. 

\subsection{Asymptotics of $w_2 (a)$, $w_2^\prime (a)$ for 
small $\lambda$}\label{Acalc}

Since $w_2 (a) = \frac{1}{a} \mathfrak{W}_{\kappa, l+1/2} (2 \lambda a)$,
with $\kappa = b/(2 \lambda)$, 
it follows from 
formula (13.1.33) and  analytic continuation 
to larger values of $\kappa$  of (13.2.5)  of \cite{Abramowitz}, p. 505 
and the identity 
$\Gamma (x) \Gamma (1-x) 
=\pi/\sin [\pi x] $ that 
\begin{multline}
\label{A.2} 
w_2 (a) = -\frac{e^{-i \pi (l-\kappa)}   
e^{-\lambda a} (2 \lambda a)^{(l+1)} \Gamma (\kappa-l)}{
2 \pi i a } H(2 \lambda a; \kappa, l) \\
{\rm where} ~~   
H(z; \kappa, l) = \int_C e^{-z t} t^{l-\kappa} (1+t)^{l+\kappa} dt
\end{multline}
where the contour $C$ starts at $\infty e^{i 0}$, circles around the
origin once counter-clockwise to the right of $t=-1$ and goes to
$\infty e^{i 2 \pi}$. In defining the integrand, we choose
$\arg t \in [0, 2\pi]$, $\arg (1+t) \in (-\pi, \pi]$ so that there
is no branch cut on the real axis between $-1$ and $0$.

It follows from (\ref{A.2}) that
\begin{multline}
\label{A.2.0} 
w_2^\prime (a) = \left (-\lambda + \frac{l}{a} \right ) w_2 (a)
+\frac{e^{-i \pi (l-\kappa)}   
e^{-\lambda a} (2 \lambda)^{(l+2)} a^l \Gamma (\kappa-l)}{
2 \pi i } H_1 (2 \lambda a; \kappa, l) \\
{\rm where} ~~   
H_1 (z; \kappa, l) = \int_C e^{-z t} t^{l-\kappa+1} (1+t)^{l+\kappa} dt
\end{multline}

We now seek to determine $H(2\lambda a; b/(2 \lambda), l)$
and $H_1 \left (2 \lambda a, b/(2 \lambda), l \right ) $  
asymptotically for small $\lambda$. 
For that purpose it
is convenient to define 
\begin{equation}
\label{A.3}
\epsilon_2 = 2 \lambda \left ( \frac{a}{b} \right )^{1/2} ~~\\,
\tau = \epsilon_2 t ~~\\,~~~P(\tau; \epsilon_2 ) = -\frac{1}{\epsilon_2} 
\log \left ( 1+ \frac{\epsilon_2}{\tau} \right ) + \tau , 
\end{equation} 
where we use  the
principal branch of $\log$ in defining $P(\tau; \epsilon_2)$ 
above. 
Then, noting that in the definition of
$\log ~\tau$ and $\log ~(\tau+\epsilon_2)$, $\arg \tau \in [0, 2\pi )$ and
$\arg ~(\tau + \alpha) \in (-\pi, \pi ]$, we have
$$ \log \tau - \log (\tau+\epsilon_2)
= -\log \left ( 1+ \frac{\epsilon_2}{\tau} \right )$$
for $\tau$ in the upper-half plane,
while
for $\tau$ in the lower-half plane, we have
$$ \log \tau - \log (\tau+\epsilon_2)
= i 2 \pi -\log \left ( 1+ \frac{\epsilon_2}{\tau} \right )$$
it is readily checked that
\begin{multline}
\label{A.4}
H(2 \lambda a; \frac{b}{2 \lambda} , l) 
= \frac{b^{l+1/2}}{2^{2l+1} \lambda^{2l+1} a^{l+1/2} } \left \{
\int_{C_1}
\tau^{l} (\tau+\epsilon_2)^l 
\exp \left[  -\sqrt{ab} P(\tau; \epsilon_2) \right ] d\tau \right .\\
\left. + \exp \left (-\frac{i \pi b}{\lambda} \right )
\int_{C_2}
\tau^{l} (\tau+\epsilon_2)^l 
\exp \left [ -\sqrt{ab} P(\tau; \epsilon_2) \right ] d\tau \right \}
\end{multline}
\begin{multline}
\label{A.4.1}
H_1 (2 \lambda a; \frac{b}{2 \lambda} , l) 
= \frac{b^{l+1}}{2^{2l+2} \lambda^{2l+2} a^{l+1} } \left \{
\int_{C_1}
\tau^{l+1} (\tau+\epsilon_2)^l 
\exp \left [ -\sqrt{ab} P(\tau; \epsilon_2 ) \right ] d\tau \right .\\
\left. + \exp \left (-\frac{i \pi b}{\lambda} \right ) 
\int_{C_2}
\tau^{l+1} (\tau+\epsilon_2 )^l 
\exp \left [ -\sqrt{ab} P(\tau; \epsilon_2 ) \right ] d\tau \right \}
\end{multline}
Here $C_1$ is a contour in the upper-half complex $\tau$-plane 
from $+\infty$ to $-\epsilon_2$ along a steepest descent line, passing through the saddle  point 
$\tau_{s,1} = i(1+o(1)) $, where 
$P^\prime (\tau_{s,1}; \epsilon_2 ) = 0$. The contour 
$C_2$ is the steepest descent line in the lower-half $\tau$-plane
from
$\tau=-\epsilon_2$ to $+\infty$ through the  saddle point 
$\tau_{s,2}, 
= -i(1+o(1))$ where  $P^\prime (\tau_{s,2}; \epsilon_2 ) =0$.  
We rewrite $w_2$ and $w'_2$ as 
\begin{equation}
\label{A.5.0}
w_2 (a) = \frac{(-1)^{l+1} e^{-\lambda a} b^{l+1/2} \Gamma (\kappa-l) }{2^{l+1} 
\sqrt{a} \lambda^{l} } \left [ Z M_1 (\sqrt{ab}, \epsilon_2 ) + Z^{-1} 
M_2 (2 \sqrt{ab}, \epsilon_2 ) \right ] 
\end{equation}
where 
\begin{equation}
\label{A.9}
M_1 (\zeta, \epsilon_2 ) 
= \frac{1}{\pi i} 
\int_{C_1} e^{-\zeta P (\tau; \epsilon_2 ) } \tau^{l} (\tau+\epsilon_2)^l 
d\tau \\,  ~~~ 
M_2 (\zeta, \epsilon_2 ) 
= \frac{1}{\pi i} 
\int_{C_2} e^{-\zeta P (\tau; \epsilon_2 ) } \tau^{l} (\tau+\epsilon_2 )^l 
d\tau 
\end{equation}
\begin{multline}
\label{A.5.0.0}
w^\prime_2 (a) = \left ( -\lambda+ \frac{l}{a} \right ) w_2 (a) + 
\frac{(-1)^{l} e^{-\lambda a} b^{l+1} \Gamma (\kappa-l) }{2^{l+1} 
a \lambda^{l} } \left [  Z M_3 (\sqrt{ab} )  + Z^{-1} M_4 (\sqrt{ab} ) 
\right ] 
\end{multline} 
where
\begin{multline}
\label{A.9.1}
M_3 (\zeta, \epsilon_2 ) 
= \frac{1}{\pi i}
\int_{C_1} e^{-\zeta P (\tau; \epsilon_2) } \tau^{l+1} (\tau+\epsilon_2)^l 
d\tau \\ \text{and}\ \ \ \  
M_4 (\zeta, \epsilon_2) 
= \frac{1}{\pi i}
\int_{C_2} e^{-\zeta P (\tau; \epsilon_2 ) } \tau^{l+1} (\tau+\epsilon_2)^l 
d\tau 
\end{multline}
It follows that, with $\epsilon_2 =2\lambda \sqrt{a/b } $, we have
\begin{equation}
\label{A.9.1.1}
\frac{w_2^\prime (a)}{w_2 (a)} 
= -\lambda + \frac{l}{a} -\left ( \frac{b^{1/2}}{a^{1/2}} \right )   
\left ( \frac{Z^2 M_3 (\sqrt{ab}, \epsilon_2 ) + M_4 (\sqrt{ab}, \epsilon_2) 
}{  
Z^2 M_1 (\sqrt{ab}, \epsilon_2 ) + M_2 (\sqrt{ab}, \epsilon_2 ) } \right )
\end{equation}

\subsubsection{Analyticity in $\epsilon_2$} 
\begin{Proposition}
  The functions $M_i(\sqrt{ab},\cdot), \,\, i=1,...,4,$ are analytic near zero.
\end{Proposition}
\begin{proof}
 We look at $M_1$, the others being similar.
We can make a change of variable
\begin{equation}
\label{A.10}
q=P(\tau; \epsilon_2 ) - P(\tau_{s,1}; \epsilon_2 ),   
\end{equation} 
 where the function $q$ is real 
on the steepest descent contour and changes
monotonically from $\infty$ to 0, as we move from $+\infty$ to 
$\tau=\tau_{s,1} $, and then increases monotonically
again from 0 to $\infty$ as
we move along the steepest descent path from $\tau=\tau_{s,1}$ to 
$\tau=-\epsilon_2 $. 
We denote the two branches of the inverse function $\tau(q)$
in  (\ref{A.10}) by $\tau_1 (q)$ and $\tau_2 (q)$.  
Noting that
$$\dfrac{dP(\tau; \epsilon_2 )}{d\tau} = 
\dfrac{1}{\tau (\tau+\epsilon_2 )} + 1, $$ 
we have    
\begin{multline}
\label{A.11}
M_1 (\zeta, \epsilon_2 ) =
e^{-\zeta \tau_{s,1} } 
\left ( \int_0^\infty e^{-\zeta q} 
\left ( \frac{\tau_2^{l+1} (\tau_2+\epsilon_2)^{l+1}}{\tau_2^2 +1 + 
\epsilon_2 
\tau_2 } \right ) dq    \right . \\   
\left . - \int_0^\infty e^{-\zeta q} 
\left ( \frac{\tau_1^{l+1} (\tau_1+\epsilon_2 )^{l+1}}{\tau_1^2 +1 
+ \epsilon_2 
\tau_1 } \right ) dq  \right )   
\end{multline}
It is easy to check that  $(\tau_i-\tau_{s_1})^2$ is  analytic for small 
$\epsilon_2 ,$, regular in $q$ and nonzero at $\epsilon_2=0$ for all $q$.

Furthermore, the integrands in (\ref{A.11}) are clearly bounded by 
an $L^1$ function uniformly in $\epsilon_2$ 
(see (\ref{A.3}) and (\ref{A.10})), 
ensuring $\epsilon_2$-analyticity of the
integrals.
\end{proof}
Returning to the original variable $\tau$ we get
\begin{multline}
\label{A.M1.0}
\frac{1}{\pi i} M_1 (\sqrt{ab}, 0) = \frac{1}{\pi i}
\int_{C_{1,0}} 
e^{-\sqrt{ab} \left (-\frac{1}{\tau} + \tau \right )  } \tau^{2l} 
d\tau  
= \int_{0}^\pi \exp \left [ i (2l+1) \theta - 2 \sqrt{ab}  
\sin \theta \right ] d\theta \\ 
= J_{2l+1} \left (2 \sqrt{ab} \right ) - i \left [ Y_{2l+1} 
\left (2 \sqrt{ab} \right ) + \mathcal{G}_{2l+1} \left ( 2 \sqrt{ab} 
\right ) \right ]  
\end{multline}
and 
\begin{multline}
\label{A.M2.0}
\frac{1}{\pi i} M_2 (\sqrt{ab}, 0) = \frac{1}{\pi i}
\int_{C_{2,0}} 
e^{-\sqrt{ab} \left (-\frac{1}{\tau} + \tau \right )  } \tau^{2l} 
d\tau  
= \int_{-\pi}^0 \exp \left [ i (2l+1) \theta - 2 \sqrt{ab}  
\sin \theta \right ] d\theta \\ 
= J_{2l+1} \left (2 \sqrt{ab} \right ) + i \left [ Y_{2l+1} 
\left (2 \sqrt{ab} \right ) + \mathcal{G}_{2l+1} \left ( 2 \sqrt{ab} 
\right ) \right ]  
\end{multline}
where $J_{2l+1}$ and $Y_{2l+1}$ are the usual Bessel functions of 
order $2l+1$ and 
\begin{multline}
\label{A.6.1}
\mathcal{G}_{2l+1} (\nu) \equiv
\frac{1}{\pi} \int_0^\infty \left \{ \exp [ (2l+1) t ] 
+ (-1)^{2l+1} \exp [-(2l+1) t ]  \right \} e^{-\nu \sinh t} dt \\
=\frac{2}{\pi} \int_0^\infty \sinh \left ( (2l+1) t \right )
e^{-\nu \sinh t} dt 
\end{multline}

\begin{equation}
\label{A.12}
\tau_{s,1} = i \sqrt{1-\frac{\epsilon_2^2}{4} } - \frac{\epsilon_2}{2} 
= i + {\rm Series ~in~} ~\epsilon_2 
\end{equation}

Thus, asymptotically, to the leading order in $\lambda$, we have
with $\nu=2\sqrt{ab}$,
\begin{multline} 
-\sqrt{\frac{a}{b}} \frac{w_2^\prime (a)}{w_2 (a)} = \\
\frac{ \left [ Z^2 \left ( J_{2l+2} (\nu) - i Y_{2l+2} 
(\nu) - i \mathcal{G}_{2l+2} (\nu) \right ) + 
\left ( J_{2l+2} (\nu) + i Y_{2l+2} 
(\nu) + i \mathcal{G}_{2l+2} (\nu) \right ) \right ]}{
\left [ Z^2 \left ( J_{2l+1} (\nu) - i Y_{2l+1} 
(\nu) - i \mathcal{G}_{2l+1} (\nu ) \right ) + 
\left ( J_{2l+1} (\nu ) + i Y_{2l+1} 
(\nu ) + i \mathcal{G}_{2l+1} (\nu ) \right ) \right ]} \\
\times 
\left ( 1+ O(\lambda) \right ) 
\end{multline}

The discussion on $w_2^\prime (a)/w(a)$ shows that 
\begin{equation}
  \label{eq:valA1}
  A=\frac{f_0(a) w'_2(a)-f'_0(a) w_2(a)}{m'_1(a) w_2(a)-
m_1(a) w'_2(a)}
\end{equation}
is an analytic function of the extended parameter set $X$ for 
$X = \left (\sqrt{p}, Z \right ) \in
\overline{\mathbb{D}^+_\epsilon}\times \overline{\mathbb{D}} $
as long as the denominator for $A$ is nonvanishing as $\lambda \rightarrow 0$. 

We can prove it is nonvanishing by simplifying
the leading order expression in $\lambda$
for $w_2^\prime (a)/w_2 (a)$, defined as $w_{2,0}^\prime (a)/w_{2,0} (a)$
under further assumption that $a$ and $c$ (as in the definition of $\beta$) 
are sufficiently large.

\subsubsection{Further simplification for large $a$}

For large $a$, there is additional simplification 
since
\begin{equation}
\label{A.7}
J_{2l+1} \left ( 2 \sqrt{ab} \right ) 
\pm{i} Y_{2l+1} \left ( 2 \sqrt{ab} \right ) 
\sim \left ( \frac{(-1)^{l+1}}{\pi^{1/2} a^{1/4} b^{1/4} } \right ) 
~\exp \left [ \pm{i} \left ( 2 \sqrt{ab} + \frac{\pi}{4} \right ) \right ]
\end{equation}
\begin{equation}
\label{A.7.0}
J_{2l+2} \left ( 2 \sqrt{ab} \right ) 
\pm{i} Y_{2l+2} \left ( 2 \sqrt{ab} \right ) 
\sim \left ( \frac{(-1)^{l+1}}{\pi^{1/2} a^{1/4} b^{1/4} } \right ) 
~\exp \left [ \pm{i} \left ( 2 \sqrt{ab} - \frac{\pi}{4} \right ) \right ]
\end{equation}
and from Watson's Lemma, we get
\begin{equation}
\label{A.7.1}
\mathcal{G}_{2l+1} \left (2 \sqrt{ab} \right ) = O(1/a) \\, ~~~
\mathcal{G}_{2l+2} \left (2 \sqrt{ab} \right ) = O(1/\sqrt{a}) 
\end{equation}
It follows that for large $a$, 
\begin{equation}
\label{A.8} 
\frac{w_{2,0}^\prime (a)}{w_{2,0} (a)}  \sim 
 \frac{r_2}{a} 
\left ( \frac{n_1 - Z^2}{Z^2+n_1 } \right ) \left ( 1 + O(a^{-1/2} ) 
\right ),  
\end{equation}
where
\begin{equation}
\label{add.2.0}
n_1 = i e^{4 i \sqrt{b a}} ~~,~~Z = \exp \left [ \frac{i \pi b}{2 \lambda}
\right ] ~~,~~
r_2 = i \sqrt{ab} 
\end{equation}

\subsubsection{Nonvanishing of the denominator of $A$ in  (\ref{eq:valA1})}

Now, defining 
\begin{equation}
\label{add.6}
m = m_1 (a)~~,~~m'=m_1^\prime (a) ~~, ~~f=f_0 (a) ~~,~~f'=f_0^\prime (a) 
\end{equation}
We have to the leading order in $\lambda$, for large $a$, 
\begin{multline}
\label{add.7}
-A = \frac{r_2 \left [ 
\dfrac{n_1-Z^2+O(\lambda)}{n_1+Z^2+O(\lambda, a^{-1/2})} 
\right ] 
- a f^\prime }{ r_2 m \left [  
\dfrac{n_1-Z^2+O(\lambda)}{n_1+Z^2+O(\lambda, a^{-1/2})} 
- \dfrac{3}{4 r_2} \right ] - a m^\prime } \\
= \dfrac{ f \left [ 4 r_2 (n_1 - Z^2) + O(\lambda)
\right ] - 4 a f^\prime \left [ n_1 + Z^2 +O(\lambda)\right ]}{ 
m \left [ 4 r_2 (n_1 - Z^2) + O(\lambda)
\right ] - 4 a m^\prime \left [ n_1 + Z^2 + O(\lambda) \right ]}
\end{multline}
The denominator of $A$ is 
\begin{equation}
\label{add.8}
D = -m \left \{ Z^2 \left ( 4 a \frac{m^{\prime}}{m} + 4 r_2 +O(\lambda)\right )
+ n_1 \left ( 4 a \frac{m^\prime}{m} - 4 r_2 +O(\lambda)\right ) \right \}
\end{equation}
We note that  
\begin{multline}
\label{add.9}
m \equiv m_1 (a) = \frac{1}{a} \mathfrak{M}_{\frac{b}{2\alpha}, l+1/2} 
\left ( 2 \alpha a \right )  = e^{-\alpha a} (2 \alpha)^{l+1} a^{l} 
M \left ( l+1-\frac{b}{2 \alpha} , 2 l + 2 , 2 \alpha a \right )  \\
\sim \left ( 2 \lambda \right )^{l+1} a^l \frac{e^{\alpha a} }{ 
\Gamma \left (l+1-\frac{b}{2 \alpha} \right ) } \left [1
+O\left ( \alpha a)^{-1} \right ) \right ]  ~~{\rm for} ~\alpha ~~
{\rm large}~
\end{multline}
and for large $\alpha$ in the fourth quadrant
\begin{equation}
\label{add.10}
m^\prime \equiv m_1^\prime (a) 
\sim 
\left ( 2 \alpha \right )^{l+1} a^{l} \alpha
\frac{e^{\alpha a}}{ 
\Gamma \left (l+1-\frac{b}{2 \alpha} \right ) }  \left [ 1 + O (\alpha a)^{-1}
\right ]
\end{equation}
\begin{equation}
\label{add.11}
\alpha = \sqrt{\lambda^2-ic} \rightarrow 
c^{1/2} \exp \left [ -i \frac{\pi}{4} \right ] ~~{\rm as} 
~\lambda \rightarrow 0 
\end{equation}
Therefore, $D$ can be zero for large enough $c$ ({\it i.e.}
large $\beta$) only if
\begin{multline*}
Z^2 \left ( 2 \sqrt{2} a c^{1/2} (1-i) \left [ 1
+O \left ( (ca)^{-1} \right )\right ] + 4 i \sqrt{ba} \right ) \\
=-n_1 \left ( 2 \sqrt{2} a c^{1/2} (1-i) \left [1 + O \left ( (ca)^{-1} 
\right ) \right ] - 4 i \sqrt{ba} \right )+O(\lambda)
\end{multline*}
Taking the absolute square of both sides, we obtain, 
\begin{multline*}
|Z|^2 \left ( [ 2 a \sqrt{2 c} ]^2 \left (1+O(c^{-1} a^{-1} \right )  + 
[ 4 \sqrt{a b} - 2 a \sqrt{2 c} \left ( 1 + O(c^{-1} a^{-1} \right ) ]^2 
\right )  \\
=     
\left ( [ 2 a \sqrt{2 c} ]^2 \left ( 1 + O(c^{-1} a^{-1} \right )  + 
[ 4 \sqrt{a b} + 2 a \sqrt{2 c} (1 + O(c^{-1} a^{-1} ) ]^2 \right ) +
O(\lambda)
\end{multline*}
This is 
impossible, since $|Z| \le 1$. This means that for large
enough $a$ and $c$ (that is, $\beta$ large),  
$D$ cannot be zero.  It means that the
resolvent is well-defined as $p =0$ is approached from the closure of
$\mathbb{H}$.  \smallskip

\begin{Note}\label{NoteAcal}
Note that the denominator of 
$A$ in (\ref{add.7}) vanishes at points in the region $|Z| > 1$, where, as 
a result, the resolvent $\mathfrak{R}_\beta$
has poles. From the relation between $Z$ and  $p$, it follows that 
$p =0 $ is an accumulation point of a sequence of
poles in the left half plane approaching zero tangentially to $i\RR$.   
\end{Note}

\subsection{Stationary phase analysis needed to calculate the
  ionization rate}\label{Newsp}

We know that the solution ${\hat y} (is, x)$ is 
analytic in the extended parameter  $\left ( \sqrt{i s}, Z)\right ) $, 
where
\begin{equation}
\label{a.A.13}
Z= \exp \left [ i \pi b / \left ( 2 \sqrt{s} \right ) \right ]
\end{equation}     
So, for $X = \left ( \sqrt{is} ,Z \right ) 
\in \overline{\mathbb{D}^+_\epsilon}\times \overline{\mathbb{D}}$,
\begin{equation}
\label{a.A.0}
{\hat y} (i s, x) = \sum_{l=0}^\infty s^{l/2} F_l (Z) 
\end{equation}  
Consider 
\begin{equation}
\label{a.A.0.0}
G(s) \equiv
\sum_{l=4}^\infty s^{l/2} F_l 
\left (\exp \left [ \frac{i \pi b}{2 \sqrt{s}} 
\right ] \right ) 
\end{equation}
It is clear that $G(s)$ is a $C^1$ function of $s$ in $[-a,
a]$. Integration by parts gives
\begin{equation}
\label{a.A.0.0.0}
\int_{-a}^a G(s) e^{i st} ds = O(t^{-1} )
\end{equation}
Now note that 
\begin{equation}
\label{a.A.14}
F_l \left (\exp \left [ \frac{i \pi b }{2 \sqrt{s}} \right ] \right )
= \sum_{j \ge 0 } {D}_{j,l} 
\exp \left [ i \frac{\pi b j}{2\sqrt{s}} 
\right ]      
\end{equation}
with $D_{j,l}$ decreasing exponentially with $j$, because of
analyticity of $F_l (Z)$ for $|Z| \le 1$.    
For $0 \le l \le 3$, it follows there exists constants $c$ and $C$
independent of $j$
so that
\begin{equation}
\label{a.A.14.1}
\sum_{l=0}^3 |D_{j,l}| \le C e^{-c j} 
\end{equation}
It follows that for large $t$,  we have
\begin{equation}
\label{a.A.14.1111}
\lvert \sum_{l=1}^3 \sum_{j=[\sqrt{t} ]+1}^\infty  D_{j,l} \int_{-a}^a 
\exp \left [ \frac{i b j}{2 \sqrt{s}} \right ] e^{i s t} s^{l/2} ds \rvert 
\le C_1 e^{-c \sqrt{t}}   
\end{equation}
Further, for large $t$, 
\begin{equation}
\label{a.A.14.2}
\left\lvert \sum_{l=0}^3 D_{0,l} \int_{-a}^a e^{i s t} s^{l/2} ds \right\rvert 
\le \frac{C_2}{t}      
\end{equation}
Therefore, 
\begin{multline}
\label{a.A.12} 
\int_{-a}^a \sum_{l=0}^3 \left ( s^{l/2} 
F_l (Z) \right )  
e^{i s t} ds 
= \sum_{0 \le l \le 3} \int_{-a}^a s^{l/2} 
F_l (Z) e^{i s t)} ds
\\
\sim  
\sum_{0 \le l \le 3} \sum_{j=1}^{[\sqrt{t}]}
D_{j,l} 
\int_{-a}^a s^{l/2} \exp \left [ i \left \{ s t 
+ j \frac{\pi b}{2 \sqrt{s}} \right \} \right ]  ds  
~+~O \left ( \frac{1}{t} \right )
\end{multline}
   
We first evaluate the terms of the form
\begin{equation}
\label{a.A.16}
\int_{-a}^a s^{l/2} e^{i t s + i d_j/\sqrt{s} } ds 
\end{equation}
for large $t$,
where
$$ d_j = \frac{j \pi b}{2}. $$
The contribution from $\int_{-a}^0$ is obviously
small, at most $O(1/t)$, uniformly for all $t$, 
since the integrand vanishes
exponentially as $s \rightarrow 0^{-}$. So we only consider, 
for $ 1 \le j \le \left [ \sqrt{t} \right ] $, 
\begin{equation}
\label{a.A.16.1}
\int_{0}^a s^{l/2} e^{i t s + i d_j s^{-1/2} } ds 
\end{equation}
We have a point of stationary phase at $s=s_{0,j}$, where
\begin{equation}
s_{0,j} = \left ( \frac{d_j}{2 t} \right )^{2/3}    
\end{equation}
Note that $s_{0,j} \ll 1$ for $t$ large since $j$ is restricted to
$j \le \sqrt{t}$.
It is then convenient to rescale  $s = s_{0,j} q$, to obtain
\begin{equation}
\label{a.A.17.0}
s_{0,j}^{1+l/2} \int_{0}^{\frac{a}{s_{0,j}} } \exp \left [ i \nu_j 
\left ( q + 2 q^{-1/2} \right ) \right ] q^{l/2} dq, ~~~~{\rm where} ~~
\nu_j = \frac{2^{-2/3}}{d_j^{2/3}} t^{1/3}  
\end{equation}
Using standard  stationary phase arguments we obtain that, 
for large $t$, and hence large $\nu_j$,
\begin{equation}
\label{a.A.17}
\lvert s_{0,j}^{1+l/2} \int_{0}^{a/{s_{0,j}}} \exp \left [ i \nu_j 
\left ( q + 2 q^{-1/2} \right ) \right ] q^{l/2} dq -  
\frac{\sqrt{2 \pi} s_{0,j}^{l+1/2} e^{i \nu_j} }{ 
\sqrt{\nu_j} }  e^{-i\pi/4} \rvert \le C \frac{s_{0,j}^{1+l/2}}{\nu_j}        
\end{equation}
For large $t$, 
the dominant contribution comes from the term with $l=0$ and so
\begin{equation}
\label{a.A.17.1}
\left \lvert \int_{-a}^a {\hat y} (is, x) e^{i s t} ds  
- \sum_{j=0}^{\left [ \sqrt{t} \right ]} D_{j,0}    
\frac{\sqrt{2 \pi} s_{0,j} e^{i \nu_j} }{ 
\sqrt{\nu_j} }  e^{-i\pi/4}    \right\rvert  
\le C t^{-1} \sum_{j=0}^{\left [\sqrt{t} \right ]}  e^{-c j} 
\le C_1 t^{-1}
\end{equation}
The sum over $j$ is clearly convergent because of
the exponential decay of $D_{j,0}$; hence $\left [ \sqrt{t} \right ]$
in the upper limit can be replaced by $\infty$.
{From} the definition of $s_{0,j}$ and 
$\nu_j$, it follows that
\begin{equation}
\int_{-a}^a {\hat y} (is, x) e^{i s t} ds = O\left ( t^{-5/6} \right )   
\end{equation}

At all other singular points, $ p = i n \omega $, $n \in \mathbb{Z}$,
the behavior is similar, and a similar calculation gives a $e^{i n
  \omega t} t^{-5/6}$ contribution. Since $Y \in \mathcal{H}$, there
is sufficient decay in $n$ to ensure that the sum over all such
contribution is convergent.  
  
\subsection{Calculation of $j_k$}\label{jkcal}
Substituting the explicit expressions for
$m_k (r)$ and $m_{(k-1)} (r)$, it may be checked that in both cases,
$\tau=0$ and $\tau=1$, corresponding to ({\bf i}) and ({\bf ii}) 
respectively 
\begin{multline}
\label{An7}
j_k = k^2 \alpha^2 \mathfrak{s} (r) j^{(2)}_k +  k \alpha^2 \mathfrak{s} (r)
j^{(1)}_k
+ j_k^{(0)} ~~, ~~~{\rm where} \\
j^{(2)}_k =  
\frac{4 \Omega}{\alpha^2 \mathfrak{s}^2}
\Big ( 1 - \frac{H(\alpha k) H(\zeta-\zeta/k)}{H(\alpha (k-1)) H(\zeta)}
\Big ) +
\frac{H^{\prime \prime} (\zeta)}{H(\zeta)} -
\frac{l (l+1)}{\zeta^2} - 
\frac{4 \sqrt{\Omega (r)}}{\alpha \mathfrak{s} (r)}
\frac{H^\prime (\zeta)}{H(\zeta)} \\
j^{(1)}_k = 
-\frac{2 (1-2 \tau) \Omega}{\alpha^2 \mathfrak{s}^2}
\Big ( 1 - \frac{H(\alpha k) H(\zeta-\zeta/k))}{H(\alpha (k-1)) H(\zeta)}
\Big ) + \frac{b}{\alpha \zeta} +
\Bigg [ - \frac{2 \tau \sqrt{\Omega}}{\alpha \mathfrak{s} } +  
\frac{\omega \mathfrak{s}}{2 \sqrt{\Omega} \alpha} -
\frac{\Omega^\prime}{2 \alpha \Omega}
\Bigg]
\frac{H^\prime (\zeta)}{H(\zeta)} \\
j^{(0)}_k = \frac{5 \mathfrak{s} {\Omega^\prime}^2}{16 \Omega^2}
-\frac{\omega \mathfrak{s}^2 \Omega^\prime}{4 \Omega^{3/2}}
-{\frac {\mathfrak{s}
 \Omega^{\prime \prime} }{4 \Omega}}
-\frac{(1+2 \tau)}{4}  \omega \mathfrak{s} 
+ {\frac {{\omega}
^{2}{\mathfrak{s}}^{3} }{16 \Omega}} 
-(\omega n_0 -ip_1) \mathfrak{s}
\end{multline}
where $\mathfrak{s} (r) = \int_r^1 \sqrt{\Omega (s)} ds $,
$\zeta=k\alpha r$ and 
$j_k := \mathfrak{s} [\mathcal{L}_k m_k - \Omega m_{k-1}]/m_k $. 
Recall   that $H(\zeta)$ satisfies
\begin{equation}
\label{An7a}
H^{\prime\prime} = 2 \left ( 1 - \frac{\omega}{2 k \alpha^2} +
\frac{\Omega^\prime (0) (1 + 2 \zeta)}{4 k \alpha \Omega (0)} + \frac{\tau}{2k}
\right ) H^\prime + \left ( \frac{l(l+1)}{\zeta^2} - \frac{b}{\alpha \zeta k}
\right ) H, 
\end{equation}  
where
\begin{equation}
\alpha = 2 \frac{\sqrt{\Omega (0)}}{\mathfrak{s} (0)}  
\end{equation}
and that $H(\zeta)$ has the following asymptotic behavior 
\begin{equation}
\label{Hasymp}
H(\zeta) \sim 1 + \frac{l (l+1)}{2 \zeta} + \frac{b}{2 k \alpha} 
\log \zeta + O \left ( \frac{\log \zeta}{ k \zeta} , \frac{1}{\zeta^2} 
\right ) ~~~,~~~\left ( \zeta, k \rightarrow \infty, ~~,~~\zeta \le k \alpha
\right ) 
\end{equation}

Now, we claim that for any $r \in (0, 1)$,  
$| j_k^{(2)} + k^{-1} j_k^{(1)} | \le C k^{-2} $.  
In the regime $r \ll 1$, we use Taylor expansion:
\begin{equation}
\Omega = \Omega (0) + \Omega^\prime (0) \frac{\zeta}{k \alpha} ~+
~O \left (\frac{1}{k^2} \right ) ~~,~~
\mathfrak{s} = \mathfrak{s} (0) - \sqrt{\Omega (0)} \frac{\zeta}{k \alpha}
~+~O \left ( \frac{1}{k^2} \right )
\end{equation}
and substitute $r = \zeta/(k\alpha)$ in 
$j_k^{(2)} + k^{-1} j_k^{(1)} $; we then  use 
$\alpha =2 \sqrt{\Omega (0)}/\mathfrak{s}(0)$,
(\ref{An7a}) and the asymptotic behavior (\ref{Hasymp}) to evaluate 
$H(\alpha k)$ and $H(\alpha (k-1))$ to find 
$j_k^{(2)} + k^{-1}j_k^{(1)} \sim k^{-2} g(\zeta) $
for some bounded differentiable function $g(\zeta)$, with asymptotic
behavior $g(\zeta) \sim {{\rm const.}}/{\zeta}$ for large $\zeta$. 
When $r $ is not small, we use the asymptotic behavior
(\ref{Hasymp}) to evaluate all terms involving the function $H$ and
to find the same inequality $|j_k^{(2)} + k^{-1} j_k^{(1)} | \le C k^{-2} $. 

Therefore, $j_k (r) = O(1)$ in all regimes. 
Further, it is easily checked that in the regime
$ k \gg \zeta \gg 1$,
$j_k (r ) = O\left ( 1, \zeta^{-1} \right ) = O(1/(kr), 1)$. 
Since the asymptotics is differentiable (since $H$ satisfies a second
order differential equation), it follows
$j_k^\prime (r) = O (k^{-1} r^{-2}, 1)$.
When $r $ is not small, using
(\ref{Hasymp}), it is readily checked that
$j_k^\prime = O(1)$.

\subsection{Generalizations} 
\label{Hr}
In fact, the same asymptotic arguments hold more generally if 
$$V(t,x)=\sum_{j=-M}^M e^{i j \omega t} \Omega_j (r)$$
with $\Omega_j (r)$ satisfying the conditions we used for $\Omega$. We
substitute for $r = O(1)$
$$ g_{n_0-k} (r) = 
\frac{c_*}{\Gamma (2k/M+1)}
\exp \left [ 
k \log f_0 (r) + 
\sum_{j=1}^{M} k^{1-j/M} f_j (r) \right ] 
$$
and calculate the error term $R_k$ as before. By requiring that the
$O(k^{2-2j/M})$ terms vanish for $j=0,..,M$, we obtain
$(M+1)$ first order differential equations for $f_j$. 
To leading order
$$ f_0 (r) =  
\left [ \int_{r}^1 \sqrt{\Omega_{-M} (s)} ds \right ]^{2/M} $$
The expressions for $f_j (r)$ for $j \geqslant 1$ are more complicated and
involve  arbitrary constants  to be determined from
the information for small $k$ at $r=1$.
Again because of the presence of $r^{-2}l (l+1)$ in 
$\mathcal{L}_k$, the remainder is $O(r^{-2})$, which is $O(k^2)$
when $r = O(k^{-1})$. We write
\begin{equation}
\label{A.0}
g_{n_0-k} (r) \sim  
 c_*
\exp \left [ 
k \log f_0 (r) + \sum_{j=1}^{M} k^{1-j/M} f_j (r) \right ] 
\frac{H(\alpha k r)}{\Gamma (2k/M+1)}
\end{equation}
Then, if $\zeta = O(1)$, we find  to leading order 
$H(\zeta) \sim H_0 (\zeta)$ where
$$ H_0^{\prime \prime} - 2 H_0^\prime - \frac{l (l+1)}{\zeta^2} H_0 = 0$$
where now $\alpha = 2 \sqrt{\Omega_{-M} (0)}/{\mathfrak{s} (0)} $
and $ \mathfrak{s} (r) = \int_r^1 \sqrt{\Omega_{-M} (s)} $.
As for $M=1$, we 
have to require $H_0 (\zeta) \sim 1$ as $\zeta \rightarrow \infty$.
This leads to 
$$H_0 (\zeta) = \sqrt{\frac{2}{\pi}} e^{\zeta} \zeta^{1/2} 
K_{l+\frac{1}{2}} (\zeta)
$$
For nonzero $g_{n_0-k}$, 
the constant multiple in (\ref{A.0}) is expected to
be nonzero. 
On the other hand, the asymptotic behavior as $\zeta \downarrow 0$,
$H_0(\zeta) \sim c_* \zeta^{-l}$ implies that the behavior
at $r=0$ of 
${g_{n_0-k}}/{r}$ is not acceptable unless  every $g_n $ vanishes
identically.

\centerline{*}
 
The analysis is likely to extend to 
systems with $H_C$  replaced by 
$$H_W=-\Delta-b/r+W(r)$$ where $b$ may be zero and $W(r)=O(r^{-1-\epsilon})$ for large $r$ and is in $ L^\infty (\RR^3)$. Under these assumptions, $W(r)$ does
not participate in the asymptotics, to the orders relevant to the proofs.

\subsection{Further remarks on the asymptotics}\label{improvements}

\begin{Remark} 
\label{weakT3}{\rm 
A weaker statement than Theorem \ref{Rnasympt} suffices to
  complete the proof of Theorem~\ref{T3}. For instance, it suffices to show
  that for sufficiently large $j$, $|R_{k,j}| < 1$, where
$$
r^{l+1}v_{n_0-k_j} (r) = i^{k_j} r^{l}m_{k_j} (r) [1 + R_{k_j} (r)]$$
}\end{Remark}

\begin{Remark}
  {\rm Stronger results than those in Proposition \ref{P33.1} hold.
    Noting that for any integer $q\geqslant 0$ we have
$$ \| \mathcal{A}_{k_j+q}... 
\mathcal{A}_{k_j+2} \mathcal{A}_{k_j+1} \mathcal{A}_{k_j} 
[{\tilde h} - 1] \|_\infty \leqslant  
\prod_{q'=0}^\infty \left ( 1 + \frac{c_*}{(k_j+q')^2} \right ) 
\| {\tilde h} - 1 \|_\infty      
$$
while
$$ \mathcal{A}_{k_j+q}... 
\mathcal{A}_{k_j+2} \mathcal{A}_{k_j+1} \mathcal{A}_{k_j} [1] = 
1 + O(k_j^{-1}) $$
and the fact that 
$ \| \mathcal{H}_{k_j+q} {\tilde h} \|_\infty 
\leqslant  c_* k_j^{-2}$, it follows that the sequence 
${\tilde h}_{k}$, satisfying
$$
{\tilde h}_{k} = \mathcal{A}_k {\tilde h}_{k-1} + \mathcal{H}_k
{\tilde h}_{k+1}, $$ has the property $\lim_{k\to\infty} {\tilde h}_k
= 1$.  Indeed, this is in accordance with  the heuristic arguments
presented in \S\ref{Hr}.  While these results completely justify the formal
asymptotics, they  are not needed in the proofs and we omit the details.
}
\end{Remark}

\section{Acknowledgments} We thank R. D.  Costin, S. Goldstein, W. Schlag, A. Soffer and
C. Stucchio for very useful discussions.  We are very grateful to Kenji Yajima
for many useful comments and suggestions on earlier drafts of this paper. 

Work supported by in part by NSF
Grants DMS-0100495, DMS-0406193, DMS-0600369, DMS-0100490, DMS-0807266, 
 DMR 01-279-26 and
AFOSR grant AF-FA9550-04. O. C. and J. L. L. acknowledge the partial support
from IAS and IHES and S.T. acknowledges support by the EPSRC and the
Mathematics Institute at Imperial College during his 2005-2006 stay.  Any
opinions, findings, conclusions or recommendations expressed in this material
are those of the authors and do not necessarily reflect the views of the
National Science Foundation.
\printindex

\end{document}